\documentclass[english,reqno,12pt]{article}
\usepackage{amsmath}
\usepackage{amsthm}
\usepackage{amssymb}
\usepackage{amsfonts}
\usepackage{epic}
\usepackage{eepic}
 \usepackage{times}
\usepackage[matrix,arrow]{xy}
\usepackage{macros}
\usepackage{epsfig} 

\theoremstyle{plain}

\newtheorem{thm}{Theorem}

\newtheorem{propn}{Proposition}
\newtheorem{lem}{Lemma}

\newtheorem{conj}{Conjecture}

\newtheorem{defn}{Definition}

\theoremstyle{remark}
\newtheorem{rem}{Remark}

\newcommand{\sst}{\scriptscriptstyle}

\renewcommand{\1}{\one}
\renewcommand{\2}{\two}
\newcommand{\3}{\three}
\newcommand{\4}{{\mathfrak 4}}

\newcommand{\beq}{\begin{equation}}
\newcommand{\eeq}{\end{equation}}

\newcommand{\id}{{\rm id}}

\newcommand{\pa}{\partial}
\newcommand{\ot}{\otimes}
\newcommand{\ra}{\to}

\newcommand{\ti}{\times}

\newcommand{\fr}[2]{{\textstyle \frac{#1}{#2} }}

\newcommand{\fsl}{{\mathfrak s}{\mathfrak l}}

\newcommand{\RI}{{\mathbb I}}

\newcommand{\RE}{{\mathbb E}}

\newcommand{\RJ}{{\mathbb J}}

\newcommand{\df}{\equiv} 
\newcommand{\hz}{\hat{z}} 
 
\newcommand{\al}{\alpha}

\newcommand{\be}{\beta}
\newcommand{\ga}{c}
\newcommand{\Ga}{\Gamma}
\newcommand{\de}{\delta}

\newcommand{\ep}{\epsilon}
\newcommand{\la}{\lambda}
\newcommand{\om}{\omega}
\newcommand{\Om}{\Omega}
\newcommand{\si}{\sigma}

\newcommand{\vf}{\varphi}
\newcommand{\tpi}{\tilde{\pi}}

\newcommand{\CA}{{\mathcal A}}

\newcommand{\BFB}{{\rm\mathbf B}}
\newcommand{\BFC}{{\rm\mathbf C}}
\newcommand{\BFD}{{\rm\mathbf D}}
\newcommand{\BFO}{{\rm\mathbf O}}

\newcommand{\CC}{{\mathcal C}}
\newcommand{\CD}{{\mathcal D}}

\newcommand{\CF}{{\mathcal F}}
\newcommand{\CG}{{\mathcal G}}
\newcommand{\CH}{{\mathcal H}}

\newcommand{\CI}{{\mathcal I}}

\newcommand{\CK}{{\mathcal K}}
\newcommand{\CL}{{\mathcal L}}
\newcommand{\CM}{{\mathcal M}}

\newcommand{\CO}{{\mathcal O}}  
\newcommand{\CP}{{\mathcal P}}  
  
\newcommand{\CR}{{\mathcal R}}

\newcommand{\CT}{{\mathcal T}}
\newcommand{\TCT}{{\widetilde{\mathcal T}}}
\newcommand{\CU}{{\mathcal U}}
\newcommand{\CV}{{\mathcal V}}

\newcommand{\CZ}{{\mathcal Z}}

\newcommand{\SA}{{\mathsf A}}
\newcommand{\SB}{{\mathsf B}}
\newcommand{\SC}{{\mathsf C}}
\newcommand{\SD}{{\mathsf D}}
\newcommand{\SE}{{\mathsf E}}
\newcommand{\SF}{{\mathsf F}}
\newcommand{\SG}{{\mathsf G}}
\newcommand{\SH}{{\mathsf H}}
\newcommand{\SI}{{\mathsf I}}
\newcommand{\SJ}{{\mathsf J}}
\newcommand{\SK}{{\mathsf K}}

\newcommand{\SL}{{\mathsf L}}

\newcommand{\SM}{{\mathsf M}}

\newcommand{\SO}{{\mathsf O}}  
\newcommand{\SP}{{\mathsf P}}  
\newcommand{\SQ}{{\mathsf Q}}  
\newcommand{\SR}{{\mathsf R}}

\newcommand{\ST}{{\mathsf T}}
\newcommand{\SU}{{\mathsf U}}
\newcommand{\SV}{{\mathsf V}}
\newcommand{\SW}{{\mathsf W}}

\newcommand{\SZ}{{\mathsf Z}}

\newcommand{\fc}{{\mathfrak c}}

\newcommand{\ff}{{\mathfrak f}}
\newcommand{\fh}{{\mathfrak h}}

\newcommand{\fp}{{\mathfrak p}}
\newcommand{\fq}{{\mathfrak q}}

\newcommand{\fv}{{\mathfrak v}}

\newcommand{\ft}{{\mathfrak t}}
\newcommand{\fw}{{\mathfrak w}}

\newcommand{\csi}{\tau}
\newcommand{\hsi}{\si}
\newcommand{\cvf}{\theta}
\newcommand{\hvf}{\vf}

\newcommand{\sa}{{\mathsf a}}

\newcommand{\sd}{{\mathsf d}}

\renewcommand{\sf}{{\mathsf f}}
\newcommand{\sh}{{\mathsf h}}
\newcommand{\sll}{{\mathsf l}}

\newcommand{\sm}{{\mathsf m}}

\newcommand{\sq}{{\mathsf q}}
\newcommand{\spp}{{\mathsf p}}
\newcommand{\sr}{{\mathsf r}}
\newcommand{\mss}{{\mathsf s}}
\newcommand{\mst}{{\mathsf t}}
\newcommand{\su}{{\mathsf u}}
\newcommand{\sv}{{\mathsf v}}
\newcommand{\sw}{{\mathsf w}}
\newcommand{\sx}{{\mathsf x}}
\newcommand{\sy}{{\mathsf y}}
\newcommand{\sz}{{\mathsf z}}

\newcommand{\FA}{{\mathfrak A}}
\newcommand{\FB}{{\mathfrak B}}

\newcommand{\FL}{{\mathfrak L}}

\newcommand{\zero}{{\mathfrak 0}}
\newcommand{\0}{{\mathfrak 0}}
\newcommand{\one}{{\mathfrak 1}}
\newcommand{\two}{{\mathfrak 2}}
\newcommand{\three}{{\mathfrak 3}}

\newcommand{\BR}{{\mathbb R}}

\newcommand{\BJ}{{\mathbb J}}

\newcommand{\BC}{{\mathbb C}}

\newcommand{\BZ}{{\mathbb Z}}

\newcommand{\rf}[1]{(\ref{#1})}

\newcommand{\aufz}
{\begin{list}{$\bullet$}{\topsep0cm \itemsep0cm \parsep0cm}}
\newcommand{\eaufz}{\end{list}}
\begin{document}
\title{An analog of a  modular functor from\\ 
quantized Teichm\"uller theory}
\author{J. Teschner}
\address{Institut f\"ur theoretische Physik, 
Freie Universit\"at Berlin,\\ 
Arnimallee 14, 14195 Berlin, 
Germany, and:\\[1ex]
Max Planck Institut f\"ur Gravitationsphysik
(Albert Einstein Institut)\\
Am M\"uhlenberg 1,
D-14476 Potsdam,
Germany}
\maketitle
\begin{abstract}
We review the quantization of Teich\-m\"uller spaces as initiated 
by Kashaev and Checkov/Fock. It is shown that the quantized 
Teich\-m\"uller spaces have
factorization properties similar to those required 
in the definition of a modular functor.
\end{abstract}

\tableofcontents

\part{Introduction}

The program of the quantization of the 
Teich\-m\"uller spaces\index{quantization of Teichm\"uller spaces}  
$\CT(\Sigma)$
of Riemann surfaces $\Sigma$ which was 
started in \cite{Fo,CF} and independently in \cite{Ka1} \footnote{See  
\cite{FG,CP} for recent generalizations.}
is motivated by certain problems and conjectures from mathematical 
physics. One of the main aims of this program is to 
construct a one-parameter family of maps 
\begin{equation}\label{Tqfunct}
\Sigma \;\longrightarrow\; 
\big(~\CH_b^{\sst\CT}(\Sigma)~,~\BFO_b^{\sst\CT}(\Sigma)~,~\SM_b^{\sst\CT}(\Sigma)~\big)~, \qquad
\text{where}
\end{equation} 

(i) $b$ is a deformation parameter, related to the traditional $\hbar$ via
$b=\sqrt{\hbar}$,

(ii) $\Sigma$ is 
a two-dimensional topological surface possibly with boundary,

(iii) $\CH_b^{\sst\CT}(\Sigma)$ is a 
Hilbert-space (possibly infinite-dimensional),

(iv) $\BFO_b^{\sst\CT}(\Sigma)$ is an algebra of bounded 
operators on $\CH_b^{\sst\CT}(\Sigma)$ and,

(v) $\SM_b^{\sst\CT}(\Sigma)$ is a unitary projective representation
of the mapping class group\index{mapping class group 
representation}\index{representations of mapping class groups} 
of $\Sigma$ on $\CH_b^{\sst\CT}(\Sigma)$.\\[1ex]
The data $\big(\BFO_b^{\sst\CT}(\Sigma),\SM_b^{\sst\CT}(\Sigma)\big)$ are restricted by
the requirement that a suitably
defined limit of $\BFO_b^{\sst\CT}(\Sigma)$ for $b\ra 0$
should reproduce the commutative algebra of functions on the
Teich\-m\"uller space $\CT(\Sigma)$, whereas a natural limit $b\ra 0$
of the automorphisms
of  $\BFO_b^{\sst\CT}(\Sigma)$ which are induced by the 
representation $\SM_b^{\sst\CT}(\Sigma)$ should correspond to the 
natural action of the mapping class group\index{mapping class group} 
${\rm MC}(\Sigma)$ on $\CT(\Sigma)$.

\section{Motivation}

Motivation for studying this problem comes from 
mathematical physics. A conjecture of 
H. Verlinde \cite{V} can be formulated very schematically 
as the statement that
\begin{equation}\label{Tqisom}
\big(~\CH_b^{\sst\CT}(\Sigma)~,~\SM_b^{\sst\CT}(\Sigma)~\big)
\;\simeq\; 
\big(~\CH_c^{\rm\sst L}(\Sigma)~,~\SM_c^{\rm\sst L}(\Sigma)~\big)~, \qquad
\text{where}
\end{equation}
\begin{itemize}\itemsep -.8ex
\item[(i)] 
$\CH_c^{\rm\sst L}(\Sigma)$ has a definition in terms of the 
representation theory of the Virasoro algebra with central charge $c$
as the so-called space of conformal blocks\index{conformal blocks} associated to $\Sigma$, and\\
\item[(ii)] 
$\SM_c^{\rm\sst L}(\Sigma)$ is an 
action of the mapping class group ${\rm MC}(\Sigma)$ on 
$\CH_c^{\rm\sst L}(\Sigma)$ which is canonically associated
to the representation-theoretic definition of 
$\CH_c^{\rm\sst L}(\Sigma)$.
\end{itemize}
Part of the interest in the space $\CH_c^{\rm\sst L}(\Sigma)$
from the side of mathematical physics is due to the fact that 
the elements of $\CH_c^{\rm\sst L}(\Sigma)$ represent the 
basic building blocks in the so-called Liouville conformal 
field\index{Liouville theory}
theory \cite{TL1}. Deep connections between the
perturbative approach to quantum Liouville 
theory on the one hand and Teich\-m\"uller theory on the other hand have 
been exhibited by Takhtajan, Zograf and Teo
see e.g. \cite{TT} and references therein. 

This conjecture may be seen as a non-compact analog of
similar relations between the quantization of 
moduli spaces\index{moduli spaces} 
of flat connections on Riemann surfaces on the one hand,
and rational conformal field 
theories\index{rational conformal field theory}\index{conformal field theory}
on the other hand. 
For K being a compact group, the geometric
quantization of the moduli space $\CM_{\rm K}(\Sigma)$
of flat K-connections on a Riemann surface was performed 
in \cite{Hi2,AdPW}. Alternative approaches were based on more
explicit descriptions of the symplectic structure on 
$\CM_{\rm K}(\Sigma)$ \cite{FR,AGS,BR,AMR}.\footnote{The 
equivalence between the different quantization schemes has not been
discussed in detail so far. It boils down to the verification
that the monodromy representation of the KZ-connection constructed within
the geometric quantization framework of \cite{Hi2, AdPW} is equivalent
to the mapping class group representation defined in \cite{AGS,BR}.
It seems that e.g. combining the results of \cite{La} and \cite{BK2}
does the job.} In either case
one may schematically describe one of
the main results of these constructions as an assignment
\begin{equation}\label{Mqfunct}
\Sigma \;\longrightarrow\; 
\big(~\CH_k^{\sst\CM}(\Sigma)~,~\SM_k^{\sst\CM}(\Sigma)~\big)~, \qquad
\text{where}
\end{equation}

(i) $\CH_k^{\sst\CM}(\Sigma)$ is a finite-dimensional
vector space,

(ii) $\SM_k^{\sst\CM}(\Sigma)$ is a projective representation
of the mapping class group of $\Sigma$ on $\CH_k^{\sst\CM}(\Sigma)$.\\[1ex]
Part of the interest in these results was due to the close relations
between the representation $\SM_k^{\sst\CM}(\Sigma)$ and
the Reshetikhin-Turaev invariants\index{Reshetikhin-Turaev invariants} 
of three manifolds \cite{RT}. Another source
of interest were the relations to rational conformal field 
theory\index{rational conformal field theory}\index{conformal field theory},
which were predicted in \cite{Wi}, see
\cite{So} for a review of mathematical
approaches to the problem and further references. 
These relations may, again schematically,
be summarized as the existence of canonical isomorphisms
\begin{equation}\label{Mqisom}
\xymatrix{
&
\big(~\CH_k^{\sst\CM}(\Sigma)~,~\SM_k^{\sst\CM}(\Sigma)~\big)
\ar@{<->}[dl]_{\simeq}&\\
\big(~\CH_k^{\rm\sst W}(\Sigma)~,~\SM_k^{\rm\sst W}(\Sigma)~\big)
\ar@{<->}[rr]^{\simeq}&&
\big(~\CH_k^{\rm\sst RT}(\Sigma)~,~\SM_k^{\rm\sst RT}(\Sigma)~\big)
\ar@{<->}[ul]_{\simeq},\\
}
\end{equation}
where
\begin{itemize}\itemsep -.5ex
\item[(i)] 
$\CH_k^{\rm\sst W}(\Sigma)$ is the space of conformal 
blocks\index{conformal blocks} in the WZNW-model
associated to the compact group ${\rm K}$, which can be defined in terms of the
representation theory of the affine Lie algebra $\hat{\mathfrak g}_k$
with level $k$
associated to the Lie algebra $\mathfrak g$ of $K$,
\item[(ii)] $\SM_k^{\rm\sst W}(\Sigma)$ is the natural 
action of the mapping class group on $\CH_k^{\rm\sst W}(\Sigma)$, 
which can be defined by means of the monodromy representation of the
Knizhnik-Zamolodchikov connection\index{Knizhnik-Zamolodchikov connection},
\item[(iii)] $\CH_k^{\rm\sst RT}(\Sigma)$ is the space of 
invariants in certain tensor products of
representations of the quantum group $\CU_q(\mathfrak g)$,
\item[(iv)] $\SM_k^{\rm\sst RT}(\Sigma)$ is the mapping class group
representation on $\CH_k^{\rm\sst RT}(\Sigma)$ defined by
the construction of Reshetikhin-Turaev.
\end{itemize}
The quantization program \rf{Tqfunct} can be seen as a non-compact
analog of \rf{Mqfunct} in the following sense.
In \rf{Tqfunct} the role of the moduli space\index{moduli space} 
of flat connections 
$\CM_{\rm K} (\Sigma)$ is taken by the Teich\-m\"uller space $\CT(\Sigma)$,  
which can be identified with the component in the moduli space 
$\CM_{\rm G}(\Sigma)$ 
of flat $G=SL(2,\BR)$-connections that has maximal Euler class \cite{Hi1,Go2}.
Moreover, the natural symplectic structure on the moduli space 
of flat $SL(2,\BR)$-connections restricts to the Weil-Petersson 
symplectic form\index{Weil-Petersson form} 
on $\CT(\Sigma)$ \cite{Go1}. A quantization of 
the Teich\-m\"uller space\index{quantization of Teichm\"uller spaces} 
$\CT(\Sigma)$ may therefore be regarded 
as providing a quantization of a topological component in the 
moduli space $\CM_{\rm G}(\Sigma)$. 

We expect that any non-compact counterpart of the developments
mentioned above will be mathematically at least as rich as the
already known results associated to compact groups $\rm K$.
In particular we expect that certain analogs of the constructions
of Reshetikhin-Turaev and/or Turaev-Viro will capture information
on the geometry of hyperbolic three manifolds, similar and probably
related to the appearance of hyperbolic volumes in the 
asymptotic behavior of certain link invariants \cite{KaV}.

\section{Aims} 

A major step towards establishing H. Verlinde's  conjecture
\rf{Tqisom} is to show that the quantization of the 
Teich\-m\"uller spaces\index{quantization of Teichm\"uller spaces} 
\rf{Tqfunct} as initiated in 
\cite{Fo,CF}, \cite{Ka1} produces an analog of 
a {\it modular functor}\index{modular functor}. 
The basic data of a modular functor
are assignments such as \rf{Mqfunct}, which are required
to satisfy a natural set of axioms as discussed 
in Section \ref{stabMF}.
One of the most important
implications of the axioms of a modular functor  
are simple relations between the representations of the mapping class
groups associated to $\Sigma$ and $\Sigma^{\dagger c}$ respectively, where
$\Sigma^{\dagger c}$ is the surface that is obtained from
$\Sigma$ by cutting along a simple closed curve $c$. These relations
imply that the representation 
$\SM_k^{\sst\CM}(\Sigma)$ restricts to - and is generated by - the 
representations $\SM_k^{\sst\CM}(\Sigma')$
which are associated to those subsurfaces $\Sigma'$ that can be 
obtained from $\Sigma$ by cutting along a set of non-intersecting 
simple closed curves. This crucial {\it locality} property can be seen
as the hard core of the notion of a modular functor.

Within the formalisms introduced in 
\cite{Fo,CF}, \cite{Ka1}-\cite{Ka3} it is far from obvious 
that the quantization of Teich\-m\"uller spaces\index{quantization of Teichm\"uller spaces} 
constructed there
has such properties. To show that this is indeed the case is 
the main problem solved in this paper. The representation
$\SM_b^{\sst\CT}(\Sigma)$ constructed and investigated in 
\cite{Ka1}-\cite{Ka3} is obtained by exploiting the fact that
the mapping class group can be embedded into the so-called
Ptolemy groupoid\index{Ptolemy groupoid} 
associated to the transformations between different
triangulations of a Riemann surface $\Sigma$. A representation of the 
Ptolemy groupoid is constructed in \cite{Ka1}-\cite{Ka3}, which then
canonically induces a projective unitary representation of 
the mapping class group ${\rm MC}(\Sigma)$. The simplicity 
of the Ptolemy groupoid, which underlies the elegance of the 
constructions in \cite{Ka1}-\cite{Ka3} now turns out to cause 
a major problem from the point of view of our aims, since the 
above-mentioned locality properties implied in the 
notion of a modular functor are not transparently realized by
the Ptolemy groupoid. 

Essentially our task is therefore to go from 
triangulations to pants decompositions\index{pants decomposition}, 
which is the type of
decomposition of a Riemann surface $\Sigma$ that is 
naturally associated to the concept of a modular functor.
This requires to construct a change of representation
for $\CH_b^{\sst\CT}(\Sigma)$ from the one naturally
associated to triangulations of $\Sigma$ 
\cite{Ka1}-\cite{Ka3} to another one which is associated 
to pants decompositions. The main tool for doing this are the
geodesic length operators\index{geodesic length operators} 
introduced and studied 
in \cite{Fo,CF2} and \cite{Ka3,Ka4}, which are the observables
on the quantized Teich\-m\"uller spaces that are associated with the
geodesic length functions\footnote{See \cite{Wo1,Wo,Wo2} 
for some classical work on the symplectic nature of the Fenchel-Nielsen
coordinates which represents important background for our 
results.}
on the classical Teich\-m\"uller spaces.
The length operators associated to a maximal set of non-intersecting
simple closed curves turn out to furnish a set of commuting self-adjoint
operators, and the simultaneous diagonalization of these operators
defines the sought-for change of representation.

There is a natural groupoid associated with the transformations between
different pants decompositions. Of particular importance for us will be
a certain refinement of this groupoid which will be called the
modular groupoid\index{modular groupoid} ${\rm M}(\Sigma)$. The
modular groupoid ${\rm M}(\Sigma)$ has been
introduced for the study of rational conformal 
field theories\index{rational conformal field theory}\index{conformal 
field theory} 
by Moore and Seiberg in \cite{MS}, and it was further 
studied in particular in \cite{BK}. Constructing a modular 
functor\index{modular functor} is essentially 
equivalent to constructing a tower of representations of the 
modular groupoid\index{modular groupoid}. 
Our main aim in the present paper will be
to show that the quantization of Teich\-m\"uller spaces
allows one to construct
a tower of representations of the modular groupoid by 
unitary operators in a natural way.

\section{Overview}

This paper has three main parts. The first of these parts collects 
the necessary results from the 
``classical'' theory of Riemann surfaces. This includes
a review of two types of coordinate systems for the Teich\-m\"uller
spaces $\CT(\Sigma)$, one of which is associated to 
triangulations of $\Sigma$,
the other to pants decompositions\index{pants decomposition}. 
The coordinates associated 
to triangulations were first introduced by Penner\index{Penner coordinates} 
in \cite{P1}.
We will also need to discuss variants of these coordinates due 
to Fock \cite{Fo} and Kashaev \cite{Ka1} respectively. The 
changes of the underlying triangulation of $\Sigma$ 
generate a groupoid, the Ptolemy groupoid
${\rm Pt}(\Sigma)$, which has a useful representation 
in terms of generators and relations (Theorem 2).

The coordinates associated to pants decompositions are the
classical Fenchel-Nielsen coordinates\index{Fenchel-Nielsen coordinates}, which we review briefly 
in \S\ref{fencoor}. 
We furthermore explain
how the coordinates of Penner \cite{P1} and Kashaev \cite{Ka1}, which
were originally introduced to parameterize the Teich\-m\"uller spaces
of surfaces $\Sigma$ with punctures only, can be used to 
provide coordinates also for the case where the surface $\Sigma$
has holes represented by geodesics of finite length. 

The material in this part is mostly known, but it is scattered
over many places in the literature, and some basic results
were stated in the original references without a proof. We have 
therefore tried to give a reasonably self-contained and complete
presentation of the relevant material, providing proofs where
these are not available elsewhere.

The second part gives a largely self-contained
presentation of the foundations of the quantization of 
Teich\-m\"uller spaces. Our presentation is heavily inspired
by \cite{Ka1}-\cite{Ka3}, but we deviate from these references
in some important points. The treatment presented in this
paper seems to be the first complete and mathematically 
rigorous formulation of the quantum theory of the 
Teichm\"uller spaces.

The main aims of this paper are finally achieved in the third part.
We begin in \S\ref{stabMF}
by introducing the notion of a stable unitary 
modular functor\index{modular functor},
and by explaining why having a stable unitary modular functor 
is equivalent to having a tower of unitary projective representations
of the modular groupoid.

In \S\ref{modgroup} we will reformulate the main result of \cite{MS,BK}
concerning the description of ${\rm M}(\Sigma)$ in terms of 
generators and relations in a way that is convenient for us.

Of particular importance for us will be \S\ref{pt->MS},
where important first relations between 
certain subgroupoids of ${\rm M}(\Sigma)$ and 
${\rm Pt}(\Sigma)$ are observed. 

In \S\ref{lengthsec} we define the geodesic 
length operators\index{geodesic length operators} 
and establish their main properties. These results are of
independent interest since some important properties 
of the geodesic length operators had not been proven 
in full generality before.

A key step in our constructions is taken in \S\ref{lengthrep} by constructing a 
change of representation from the original one 
to a representation in which the length operators associated to 
a pants decomposition are simultaneously diagonalized. An important
feature of this construction is the fact that the unitary operator
which describes the change of representation factorizes into 
operators associated to the individual three holed spheres (trinions)
which appear in a pants decomposition.

In \S\ref{Mconstr} we construct the corresponding representation of the 
modular groupoid ${\rm M}(\Sigma)$. The operators
which represent ${\rm M}(\Sigma)$
are constructed out of compositions of the representatives for the
transformations in ${\rm Pt}(\Sigma)$. This makes it 
relatively easy to verify the relations of ${\rm M}(\Sigma)$, but the
price to pay is that some crucial locality properties are more 
difficult to prove.

\section{Outlook}

In a sequel \cite{TT2} to this paper we will calculate the 
matrix coefficients of the operators which generate
the representation of the modular groupoid explicitly. 
A close relation to the modular double\index{modular double} 
$\CD\CU_q(sl(2,\BR))$
of $\CU_q(sl(2,\BR))$ as
defined and studied in \cite{Fa,PT,BT}
will be found.

It should be noted 
that $(\CH_b^{\sst\CT}(\Sigma),
\SM_b^{\sst\CT}(\Sigma))$ will not satisfy all the 
usual axioms of a modular functor, which require, in particular,
that the vector space $\CV(\Sigma)$ assigned to each Riemann surface
should be finite-dimensional. Most importantly, however, the 
assignment $\Sigma\ra(\CH_b^{\sst\CT}(\Sigma),
\SM_b^{\sst\CT}(\Sigma))$ was up to now only constructed
for surfaces $\Sigma$ which have at least one boundary component.

What will allow us to overcome this 
unsatisfactory feature are the remarkable analytic properties that the
matrix coefficients of the operators which represent 
the modular groupoid will be shown to have. 
It turns out that the mapping class group representation 
$\SM_b^{\sst\CT}(\Sigma)$ assigned to a Riemann surface with
a boundary represented by geodesics of a fixed length
depends {\it analytically} on the values of these lengths.
The analytic properties of the matrix coefficients will furthermore
allow us to ``close a hole'' by taking a limit where
the length parameter assigned to this boundary component
approaches a certain {\it imaginary} value.
It will be shown in \cite{TT2} that the resulting mapping
class group representation is equivalent to the one on the 
surface which is obtained by gluing a disc into the
relevant boundary component.

Concerning the representation theoretic side of H. Verlinde's 
conjecture \rf{Tqisom} it should be mentioned that
a complete mathematical construction of 
$(\CH_c^{\rm\sst L}(\Sigma),\SM_c^{\rm\sst L}(\Sigma))$ is not
available so far, but 
nontrivial steps in the direction of constructing and describing 
$(\CH_c^{\rm\sst L}(\Sigma),\SM_c^{\rm\sst L}(\Sigma))$ precisely
have been taken in \cite{TL2} in the case of 
surfaces $\Sigma$ of genus zero. This includes in particular
the derivation of explicit formulae for a set of basic data which 
characterize the resulting representation
of the braid group uniquely.

The explicit computation of the matrix coefficients 
of the operators which generate
the representation of the modular groupoid carried out in
\cite{TT2} will therefore allow us to verify H. Verlinde's 
conjecture \rf{Tqisom} in the case of Riemann surfaces
of genus zero. 

We furthermore expect that it should be
possible to construct non-compact analogs of 
$(\CH_k^{\rm\sst RT}(\Sigma),\SM_k^{\rm\sst RT}(\Sigma))$
based the non-compact quantum
group $\CD\CU_q(sl(2,\BR))$, and thereby complete a non-compact analog 
of the triangle \rf{Mqisom}.\\[2ex]
{\bf Acknowledgements}\\[1ex]
The author would like to thank L. Chekhov, V. Fock and especially
R. Kashaev for useful discussions.

We furthermore gratefully 
acknowledge financial support from the DFG by a
Heisenberg fellowship, as well as the kind hospitality of 
the Humboldt University Berlin, the School of Natural Sciences
of the IAS Princeton, the
Physics Department of Chicago University
and the Caltech, Pasadena.

\newpage

\part{Coordinates for the Teichm\"uller spaces}

We will consider two-dimensional surfaces $\Sigma$ with genus $g\geq 0$
and $s\geq 1$ boundary components such that  
\[
M\; \df\; 2g-2+s\;>\; 0\,.
\]
On $\Sigma$ we will consider metrics of constant negative 
curvature $-1$. Our main interest will be 
the case where the boundary components can be represented
by geodesics of finite length. Such boundary components will also
be called holes in the following. However, to begin with we will
focus on the case where the boundary components are 
{\it punctures}, i.e. holes of vanishing geodesic circumference.

The space of deformations of the metrics of constant negative 
curvature is called the Teich\-m\"uller space $\CT(\Sigma)$. It will
be of basic importance for us to have useful systems of coordinates
for $\CT(\Sigma)$.

We will consider two classes of coordinate systems which are associated
to two types of graphs drawn on the Riemann surfaces respectively.
The first class of coordinates goes back to Penner and is associated
to triangulations of the Riemann surface or the corresponding dual 
graphs, the so-called fat graphs. We will also describe two useful
variants of the Penner coordinates\index{Penner coordinates} 
due to Kashaev and Fock respectively.

The second class of coordinates are the classical 
Fenchel-Nielsen length-twist coordinates\index{Fenchel-Nielsen coordinates}.
One may view them as being associated to a second type of graph on a 
surface $\Sigma$ called {\it marking} that determines 
in particular a decomposition
of the surface into three-holed spheres (trinions). 

In the following first part of this paper we shall describe these coordinate
systems in some detail, discuss the graphs on $\Sigma$ that
these coordinates are associated to, as well as the 
groupoids generated by the transformations between
different choices of these graphs.

\begin{figure}[b]
\centerline{\epsfxsize6cm\epsfbox{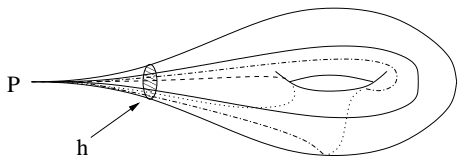}}
\caption{Triangulation of the once-punctured torus. }\label{triang}
\end{figure}

\section{The Penner coordinates\index{Penner coordinates}} \label{class}

\setcounter{equation}{0}

\subsection{Triangulations and fat graphs\index{fat graphs}}

Consider a fixed oriented topological surface $\Sigma$ of genus $g$  
with $s\geq 1$ punctures. An ideal triangulation $\tau$ of  $\Sigma$ is the
isotopy class of a collection of disjointly embedded arcs in $\Sigma$
running between the punctures such that $\tau$ decomposes $\Sigma$ into 
triangles. There are $2M$ triangles and $3M$ edges for any ideal
triangulation. As an example we have drawn a triangulation of
the once-punctured torus in Figure \ref{triang}.

\begin{figure}[htb]
\epsfxsize5cm
\centerline{\epsfbox{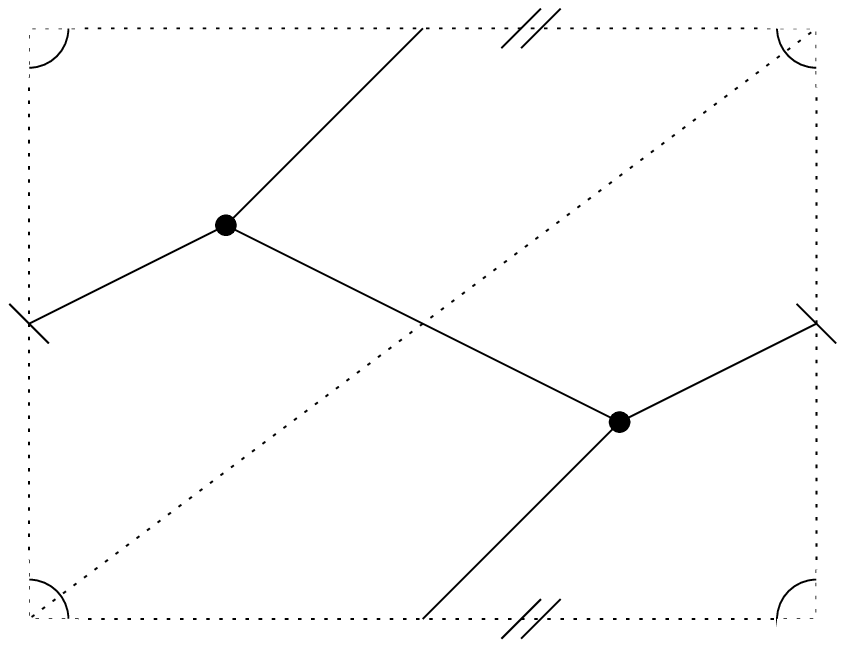}}
\caption{Another representation of the triangulation 
from Figure \ref{triang} and the dual  fat graph.}\label{fat}
\end{figure}

The graph dual to a triangulation is a trivalent fat graph i.e.
a trivalent graph embedded in the surface with fixed cyclic order 
of the edges incident to each vertex. An example for a fat graph 
on the once-punctured torus is depicted in Figure \ref{fat}.
The sets of vertices and edges of a fat graph $\vf$ will be denoted
$\vf_\0$ and $\vf_\1$ respectively. 

\begin{figure}[htb]
\epsfxsize4cm
\centerline{\epsfbox{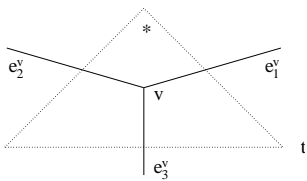}}
\caption{Graphical representation of the vertex $v$ dual
to a triangle $t$. The marked corner defines a
corresponding numbering of the edges that emanate at $v$.}\label{decor}
\end{figure}

An ideal triangulation is called {\it decorated} if the triangles
are numbered and if a corner is marked 
for each triangle \cite{Ka1}. 
The decoration of the triangle $t_v$ dual to a vertex $v\in\vf_\0$
can be used to fix a numbering convention
for the edges $e_i^v$, $i=1,2,3$ which emanate from $v$ 
as defined in Figure \ref{decor}.

\begin{rem}\label{convrem}
Decorated ideal triangulations are dual to decorated fat graphs, which means
that the vertices are numbered, and 
for each vertex $v\in\vf_\0$ one has chosen a distinguished 
edge $e_v\in\vf_\1$. As a convention we will assume that
fat graphs always carry such a decoration unless otherwise stated.
\end{rem}

\subsection{Penner coordinates\index{Penner coordinates}}

It turns out to be useful to consider a somewhat enlarged object
which keeps track of the choices of horocycles around each of the 
punctures, the so-called decorated Teich\-m\"uller space $\TCT(\Sigma)$.
$\TCT(\Sigma)$ is defined as a principal $\BR_+^s$-bundle over 
$\CT(\Sigma)$ by taking the s-tuple of horocycles around each of the
punctures as the fiber over a point of $\CT(\Sigma)$. The ordered s-tuple of 
hyperbolic lengths of the horocycles gives coordinates for the fibers.

Given any point $P$ in the decorated Teich\-m\"uller space $\TCT(\Sigma)$
and an ideal triangulation $\tau$ of $\Sigma$, Penner
assigns a coordinate value to each of the edges in $\tau$ by means
of the following construction. By means of Fuchsian 
uniformization\index{Fuchsian uniformization}
one may equip the surface $\Sigma$ with a unique hyperbolic metric $g$
associated to our chosen point $P\in \TCT(\Sigma)$. Let $\tau_\1$ 
be the set of edges of a triangulation $\tau$.
Each edge $e$ in $\tau_\1$
may be straightened to a geodesic for the 
hyperbolic metric $g$. The coordinate $l_e(P)$ is defined as the 
hyperbolic length of the segment of $e$ that lies between the 
two horocycles surrounding the punctures that $e$ connects,
taken with positive sign if the two horocycles are disjoint, with 
negative sign otherwise. We are going to consider the tuple
$(l_e)_{e\in\tau_\1}$ as a vector in the vector space $\BR^{\tau_\1}$ 
of dimension $3M$.

\begin{thm} ${}$ {\rm Penner \cite{P1}\cite{P2}}\\
(a)
For any fixed ideal triangulation $\tau$ of $\Sigma$, the 
function 
\[
l:\TCT(\Sigma)\ra \BR^{\tau_\1},\qquad P\ra (l_e(P))_{e\in\tau_\1}
\]
is a homeomorphism.\\
(b) The pull-back of the Weil-Petersson two-form $\omega$ 
on $\CT(\Sigma)$
is given by the expression
\[
\omega\;=\;-\sum_{t\in\tau_2}\bigl(dl_{e_1(t)}\land dl_{e_2(t)}+
 dl_{e_2(t)}\land dl_{e_3(t)}+dl_{e_3(t)}\land dl_{e_1(t)}\big),
\]
where the summation is extended over the set $\tau_2$ of triangles
of $\tau$, and $e_i(t)$, $i=1,2,3$ are the edges bounding the triangle
$t$, labelled in the counter-clockwise sense. 
\end{thm}
The Teich\-m\"uller space $\CT(\Sigma)$ itself can finally be described as the
space of orbits in $\TCT(\Sigma)$ under the following
symmetry. Choose a number $d(\fp)$ for each 
puncture $\fp$. Let the action of the symmetry be defined by
\begin{equation}\label{gaugesymm}
l_e' \;\df\; l_e + d(\fp)+d(\fp')
\end{equation} 
if the edge $e$ connects the punctures $\fp$ and $\fp'$.

\subsection{Fock coordinates\index{Fock coordinates}}\label{Fockvars}

\begin{figure}[t]
\epsfxsize3cm
\centerline{\epsfbox{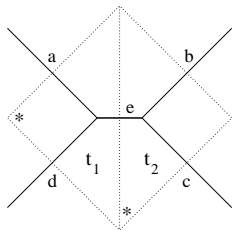}}
\caption{Two adjacent triangles and the dual fat graph.}\label{labels}
\end{figure}

There is a useful variant of the Penner coordinates which was introduced
by V. Fock in \cite{Fo}. In terms of the Penner coordinates one may
define the Fock coordinates in terms of certain 
cross-ratios\index{cross-ratio}.
Given a quadrilateral formed by two adjacent 
triangles we are going to keep the labelling of the edges introduced in 
Figure \ref{labels}. Let  
\begin{equation}\label{penner-fock}
z_e\;=\;l_a+l_c-l_b-l_d.
\end{equation}
The dependence of the Penner coordinates on the choice of
horocycles drops out in the Fock coordinates. However,
the variables $z_e$ assigned to the $3M$ edges in a triangulation
are not all independent. To describe the relations that
they satisfy it is convenient to think
of the Fock coordinates 
as being assigned to the edges of the fat graph dual to the
given triangulation. 
Each closed curve $c$ on $\Sigma$ is homotopic to a
unique path $g_c$ on the fat graph which has minimal length 
w.r.t. the metric defined by assigning each edge of $\vf$ the length one.
Such paths will also be called graph geodesics in the following.
The path $g_c$ may be described
by a sequence of edges $e_1^c,\dots,e_{n_{c}}^c\in\vf_\1$.
To a closed curve $c$ let us  associate
\begin{equation}\label{Fockconstrdef}
f_{\vf,c}\;\df\;\sum_{i=1}^{n_c} \,z_{e_i^c}\, .
\end{equation}
The definition \rf{penner-fock} then implies the relations
\begin{equation}\label{Fockconstr}
f_{\vf,c}\;=\;0
\end{equation}
for any closed curve $c$ that encircles one puncture only.
If one uses the equations \rf{Fockconstr} to express $s$ of the 
variables $z_{e}$ in terms of the others one obtains
a set of coordinates for $\CT(\Sigma)$.

On the $3M$-dimensional space $F_{\vf}$ that is spanned
by the coordinate functions $z_e(P)$ one may define
a Poisson bracket  $\Omega_{\rm WP}$ which is such that 
\begin{itemize}
\item[(i)]
the elements $f_a$, $a=1,\dots,s$
span the set $C_{\vf}$ of all $\fc\in F_{\vf}$ such that
\begin{equation}
\Omega_{\rm WP}(\fc,\fw)=0\quad {\rm for\;\, all}\;\,
\fw \in F_{\vf}, 
\end{equation}
\item[(ii)] the Poisson bracket which is induced by
$\Omega_{\rm WP}$ on the quotient $F_{\vf}/C_{\vf}$ coincides
with the Poisson bracket which corresponds to the 
Weil-Petersson symplectic form on $\CT(\Sigma)$.
\end{itemize}
There exists a rather simple description for this Poisson-bracket \cite{Fo}:
 \begin{equation}\label{poisson}
\Omega_{\rm WP}(z_{e}^{},z_{e'}^{})\;=\;n_{e,e'}, \quad{\rm where}\quad
n_{e,e'}\in\{-2,-1,0,1,2\}.
\end{equation}
The value of $n_{e,e'}$ depends on how edges $e$ and $e'$ are embedded
into a given fat graph. If $e$ and $e'$ don't have a common 
vertex at their ends, or if one of $e$, $e'$ starts and
ends at the same vertex then $n_{e,e'}=0$. In the case that 
$e$ and $e'$ meet at two vertices one has 
$n_{e,e'}=2$ (resp. $n_{e,e'}=-2$) 
if $e'$ is the first edge to the right\footnote{w.r.t.
to the orientation induced by the embedding of the fat-graph 
into the surface} (resp. left) of $e$ at both
vertices, and $n_{e,e'}=0$ otherwise. 
In all the remaining cases $n_{e,e'}=1$ (resp. $n_{e,e'}=-1$) 
if $e'$ is the first edge to the right (resp. left) of $e$ at the common
vertex.

\begin{figure}[h]
\centering
\setlength{\unitlength}{0.0004in}
\begingroup\makeatletter\ifx\SetFigFont\undefined
\def\x#1#2#3#4#5#6#7\relax{\def\x{#1#2#3#4#5#6}}%
\expandafter\x\fmtname xxxxxx\relax \def\y{splain}%
\ifx\x\y   
\gdef\SetFigFont#1#2#3{%
  \ifnum #1<17\tiny\else \ifnum #1<20\small\else
  \ifnum #1<24\normalsize\else \ifnum #1<29\large\else
  \ifnum #1<34\Large\else \ifnum #1<41\LARGE\else
     \huge\fi\fi\fi\fi\fi\fi
  \csname #3\endcsname}%
\else
\gdef\SetFigFont#1#2#3{\begingroup
  \count@#1\relax \ifnum 25<\count@\count@25\fi
  \def\x{\endgroup\@setsize\SetFigFont{#2pt}}%
  \expandafter\x
    \csname \romannumeral\the\count@ pt\expandafter\endcsname
    \csname @\romannumeral\the\count@ pt\endcsname
  \csname #3\endcsname}%
\fi
\fi\endgroup
{\renewcommand{\dashlinestretch}{30}
\begin{picture}(7244,3032)(0,-10)
\put(1822.000,270.000){\arc{1800.000}{3.1416}{6.2832}}
\put(4522.000,270.000){\arc{3600.000}{3.1416}{6.2832}}
\put(3352.000,270.000){\arc{1260.000}{3.1416}{6.2832}}
\put(5152.000,224.100){\arc{2341.800}{3.1808}{6.2440}}
\put(3622.000,307.500){\arc{5400.521}{3.1277}{6.2971}}
\thicklines
\path(22,270)(7222,270)
\put(2497,1755){\makebox(0,0)[lb]{\smash{{{\SetFigFont{10}{14.4}{rm}$t_i$}}}}}
\put(3892,1260){\makebox(0,0)[lb]{\smash{{{\SetFigFont{10}{14.4}{rm}$t_j$}}}}}
\put(3307,675){\makebox(0,0)[lb]{\smash{{{\SetFigFont{10}{14.4}{rm}$c$}}}}}
\put(5152,1100){\makebox(0,0)[lb]{\smash{{{\SetFigFont{10}{14.4}{rm}$b$}}}}}
\put(3622,2790){\makebox(0,0)[lb]{\smash{{{\SetFigFont{10}{14.4}{rm}$a$}}}}}
\put(832,0){\makebox(0,0)[lb]{\smash{{{\SetFigFont{10}{14.4}{rm}$x_1$}}}}}
\put(2632,0){\makebox(0,0)[lb]{\smash{{{\SetFigFont{10}{14.4}{rm}$x_2$}}}}}
\put(3892,0){\makebox(0,0)[lb]{\smash{{{\SetFigFont{10}{14.4}{rm}$x_3$}}}}}
\put(6232,0){\makebox(0,0)[lb]{\smash{{{\SetFigFont{10}{14.4}{rm}$x_4$}}}}}
\put(1732,915){\makebox(0,0)[lb]{\smash{{{\SetFigFont{10}{14.4}{rm}$d$}}}}}
\put(4477,1825){\makebox(0,0)[lb]{\smash{{{\SetFigFont{10}{14.4}{rm}
$e$}}}}}
\end{picture}}
\caption{Representation of the triangles $t_i$ and $t_j$ 
in the upper half plane.}\label{uhpfig}
\end{figure}
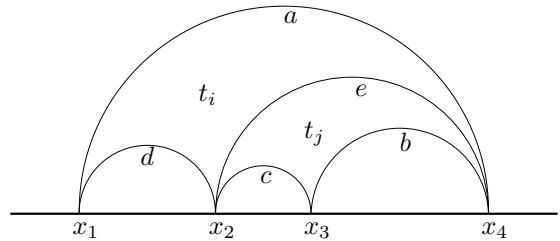

The coordinates $z_e$ also have a nice geometrical meaning \cite{Fo}.
In the Fuchsian uniformization\index{Fuchsian uniformization} 
the two triangles that share the 
common edge $e$ will be mapped into ideal hyperbolic triangles in the
upper half plane. The edges are then represented by half-circles,
and the corners will be at points $x_1,\dots,x_4$ on the real line, 
see Figure \ref{uhpfig}.

We then have 
\begin{equation}
\exp(z_e)\;=\;\frac{(x_4-x_1)(x_3-x_2)}{(x_4-x_3)(x_2-x_1)}.
\end{equation}
By means of M\"obius-transformations $x_i\ra \frac{ax_i+b}{cx_i+d}$
one may map the corners of one of the two triangles to 
$-1$, $0$ and $\infty$ respectively. The variable 
$z_e$, being expressed in terms of the M\"obius-invariant 
cross-ratio\index{cross-ratio}
therefore parameterizes the different ways of gluing two 
ideal hyperbolic triangles along a common edge modulo 
M\"obius-transformations. Given the variables $z_e$ one 
may reconstruct the
Riemann surface as represented in the Fuchsian uniformization 
by successively mapping ideal hyperbolic triangles
into the upper half-plane, glued along the edges $e$ in the way
prescribed by the given value $z_e$ \cite{Fo}.

\section{The Ptolemy groupoids\index{Ptolemy groupoid}}
\setcounter{equation}{0}

In the previous section we had associated coordinate systems to 
fat graphs on a surface $\Sigma$: 
A change of graph will of course induce a change of coordinates. 
The groupoid generated by the moves between
different fat graphs\index{fat graphs} will be the subject 
of the present 
section. 

\subsection{Groupoids vs. complexes}

The groupoids that we will be interested in can be 
conveniently described as 2-dimensional connected CW complexes $\CG$.
The set of vertices $\CG_{\zero}$ 
of $\CG$ will be represented by certain sets of graphs, within this section
called fat graphs. The (directed) edges $E\in\CG_{\1}$
that connect these vertices correspond to the 
generators (``elementary moves'') of the groupoid, while the 
faces $F\in\CG_{\2}$ of $\CG$ yield the 
relations. 

The groupoid ${\rm G}$ associated to the  
2-dimensional connected CW complexes $\CG$ will then 
simply be the path groupoid of $\CG$, which has
the vertices in $\CG_{\zero}$ as objects and
the homotopy classes of edge paths between two vertices as morphisms.
The homotopy class of paths 
leading from vertex $V_1\in\CG_\0$ 
to vertex $V_2\in\CG_\0$ will be denoted by $[V_2,V_1]$.
Similarly we will sometimes use the notation $[W_E,V_E]$ 
for the element of ${\rm G}$ which corresponds to 
an edge $E\in\CG_\1$.

A path $\pi$ which represents an element in the homotopy class $[W,V]$
may be represented by a {\it chain} of edges $E\in\CG_{\1}$, i.e.
an ordered sequence $(E_{\pi,n(\pi)},\dots,E_{\pi,\1})$, 
$E_{\pi,\jmath}\in\CG_{\1}$ for $\jmath=1,\dots,n(\pi)$ 
such that $E_{\pi,\jmath}\in[V_{\pi,\jmath+1},V_{\pi,\jmath}]$ for 
$\jmath=1,\dots,n(\pi)-1$,
and $V_{\pi,\1}=V$, $V_{\pi,n(\pi)}=W$.
We will also use the suggestive notation
$E_{n}\circ E_{n-1}\circ\dots\circ E_1$ to denote a chain.

\subsection{Change of the triangulation} \label{Ptolemy}

In the case of the Ptolemy groupoids ${\rm Pt}(\Sigma)$
we will consider a complex $\CP t(\Sigma)$, where 
the set $\CG_\0=\CP t_\0(\Sigma)$ is defined to be the set of fat graphs 
on $\Sigma$. 
Let us furthermore define $\CP t_\1(\Sigma)$ 
to consist of the following elementary moves.
\begin{itemize}\item[(i)]
{\it Permutation $(vw)$:} Exchanges the labels of the vertices $v$ and 
$w$.
\item[(ii)]
{\it Rotation $\rho_v:$} See Figure \ref{ft}.
\begin{figure}[htb]
\epsfxsize5cm
\centerline{\epsfbox{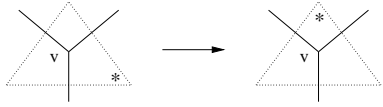}}
\caption{Transformation $\rho_v$ changes the
  marked corner of the triangle dual to a vertex $v\in\vf_\0$.}\label{ft}
\end{figure}
\item[(iii)]
{\it Flip $\om_{vw}:$} See Figure \ref{flip}.
\begin{figure}[htb]\epsfxsize7cm
\centerline{\epsfbox{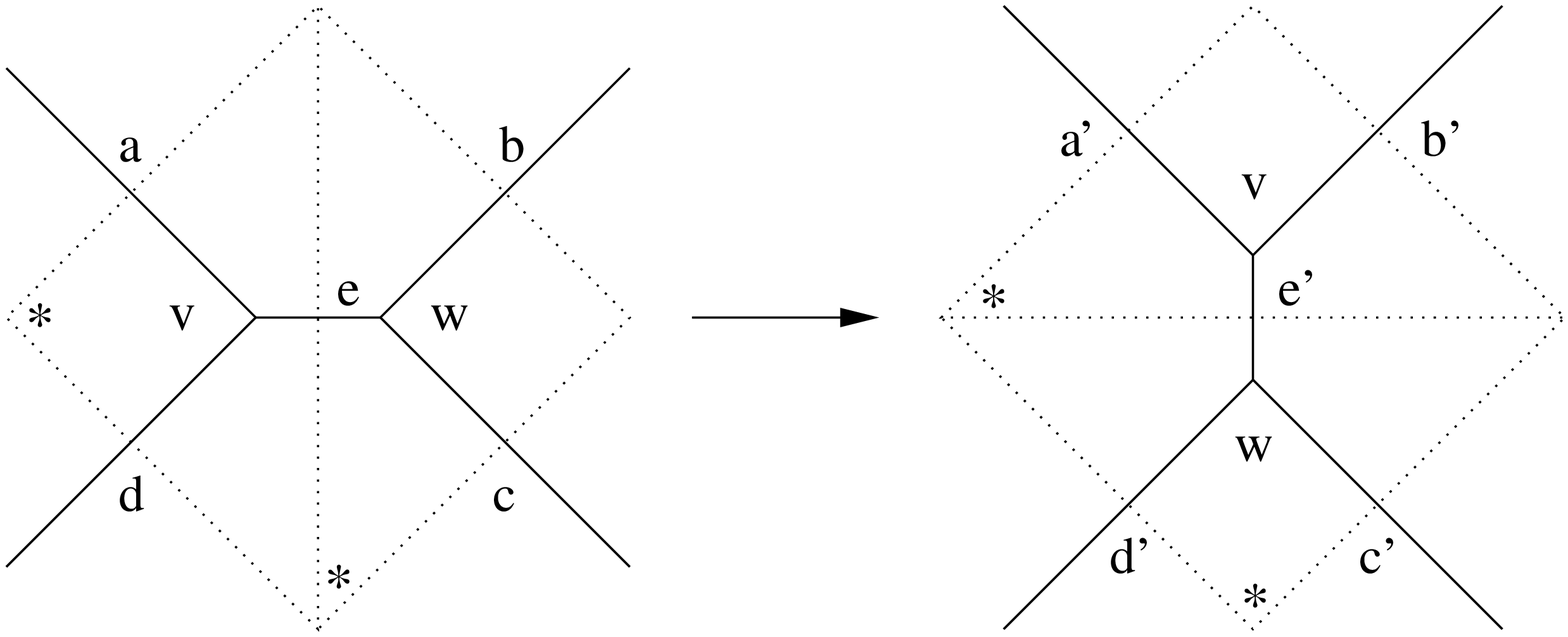}}
\caption{The flip transformation $\omega_{vw}$ 
changes the diagonal in the quadrilateral formed by 
the two adjacent triangles $t_v$ and $t_w$.}\label{flip}
\end{figure}
\end{itemize}

\begin{propn} $\;\frac{\quad}{}$
\cite[Proposition 7.1]{P1}\cite{Ka3} $\;\frac{\quad}{}$\\
The complex $\CP t(\Sigma)$ is connected, i.e.
for a given surface $\Sigma$, {\it any} two 
fat graphs $\vf$ and $\vf'$ can be connected by a
chain of elementary transformations. 
\end{propn}

Validity of the following relations in ${\rm Pt}(\Sigma)$ 
can easily be verified
pictorially \cite{Ka3}. 
\begin{equation}\label{sPtrels}\begin{aligned}
{}&  \rho_v\circ\rho_v\circ\rho_v=id,\\
& \omega_{u_\1u_\2}\circ\omega_{u_\3u_\4}=
\omega_{u_\3u_\4}\circ\omega_{u_\1u_\2},\quad u_r\neq u_s\;\;{\rm for}
\;\; r\neq s,\\
&  \omega_{vw}\circ\omega_{uw}\circ\omega_{uv}=\omega_{uv}\circ\omega_{vw},
\\
&  (\rho_v^{-1}\times\rho_w)\circ\omega_{vw}=\omega_{wv}\circ(\rho_v^{-1}
\times\rho_w),\\
& \omega_{wv}\circ\rho_v\circ\omega_{vw}=(vw)\circ(\rho_v\times\rho_w).
\end{aligned}\end{equation}

\begin{thm} \label{Ptrelcompl}
The complex $\CP t(\Sigma)$ is simply connected, i.e.
any relation between the generators $(vw)$, $\rho_v$ and $\om_{vw}$
of the Ptolemy groupoid is a consequence of the
relations \rf{sPtrels} together with the relations 
of the permutation group. 
\end{thm}

The proof of this theorem is explained in Appendix \ref{Ptrelproof}.

One of the 
main virtues of the Penner coordinates\index{Penner coordinates}
is that the corresponding
change of coordinates can be described rather simply.

\begin{lem} $\;\frac{\quad}{}$
{\rm Lemma A.1a of \cite{P2}}$\;\;\frac{\quad}{}$\\
Let $\tau'$ be 
the triangulation obtained by applying the flip of Figure \ref{flip}
to a pair of adjacent triangles in a given initial triangulation $\tau$, 
and denote $e$ and $e'$ the diagonal edge before and after the flip.
The coordinates associated to $\tau$ and $\tau'$ will then agree
for each edge that the two triangulations have in common, and
\begin{equation}
\la_{e'}\;=\;\frac{1}{\la_e}(\la_a\la_c+\la_b\la_d),\quad
\la_f\df\sqrt{2}\exp(\fr{1}{2}l_f),\;\;\forall f\in\vf_\1\,, 
\end{equation}
where we have labelled the edges according to Figure \ref{flip}.
\end{lem}
The corresponding transformation of the Fock variables\index{Fock coordinates}
is also easy to describe:
\begin{equation}\label{flipfovarclass}\begin{aligned}
e^{-z_{a'}}=& e^{-z_a}(1+e^{-z_e})\\
e^{+z_{d'}}=& e^{+z_d}(1+e^{+z_e})
\end{aligned} 
\qquad z_{e'}=-z_e\qquad\begin{aligned}
e^{-z_{b'}}=& e^{+z_b}(1+e^{+z_e})\\
e^{+z_{c'}}=& e^{-z_c}(1+e^{-z_e})
\end{aligned}
\end{equation}

\subsection{The representation of the mapping class 
group\index{mapping class group}}\label{MCGclass}

The mapping class group ${\rm MC}(\Sigma)$ of the topological surface 
$\Sigma$ is the group of isotopy classes of orientation-preserving 
diffeomorphisms of $\Sigma$. An element $\mu\in {\rm MC}(\Sigma)$
of the mapping class group will map a given (decorated) 
triangulation $\tau$
into another one, $\mu.\tau$. The fact 
that any two triangulations can be mapped into each other by a 
composition of the elementary transformations introduced in 
subsection \ref{Ptolemy} therefore leads to an embedding of
the mapping class group into the Ptolemy groupoid:
\begin{equation}\label{MCimbed}
\Phi_{\tau}: {\rm MC}(\Sigma)\;\ra\; {\rm Pt}(\Sigma),\qquad
\Phi_{\tau}(f)\;\df\;[\mu.\tau,\tau].
\end{equation}
More precisely, $\Phi_{\tau}$ induces a homomorphism 
${\rm MC}(\Sigma)\;\ra\; {\rm Pt}(\Sigma)$ in the sense that
\[ \Phi_{\tau}(\mu_\2\circ \mu_\1)\;=\;\Phi_{\mu_\1.\tau}(\mu_\2)\circ
\Phi_{\tau}(\mu_\1)\;\;
\text{for any}\;\;\mu_\2,\mu_\1\in {\rm MC}(\Sigma),
\]
which embeds ${\rm MC}(\Sigma)$ {\it injectively} into ${\rm Pt}(\Sigma)$
\cite[Theorem 1.3]{P3}.

\section{Teich\-m\"uller space as the phase space of a constrained system}
\label{constr}
\setcounter{equation}{0}

As a preparation for the description of the quantum Teich\-m\"uller spaces
it will be useful to parameterize the Teich\-m\"uller spaces by means of 
variables assigned to the {\it triangles} instead of the {\it edges} of 
a triangulation \cite{Ka1}.  In the following section 
we shall elaborate upon the results and constructions in \cite{Ka1},
strengthening them somewhat.

\subsection{Kashaev's coordinates\index{Kashaev coordinates}}

Assume
given a fat graph $\vf$ with set of vertices $\vf_{\0}$.
For each vertex $v\in\vf_{\0}$ one may introduce a pair of
variables $(q_v,p_v)$ according to the following rule. Let us label the 
edges that emanate from the vertex $v$ by $e^v_i$, $i=1,2,3$
according to  Figure \ref{decor}. 
We will
denote the Penner coordinates associated to the 
edges $e^v_i$ by $l^v_i$, $i=1,2,3$. Let us then 
define the pair of
variables $(q_v,p_v)$ as
\begin{equation}\label{trivars}
\big(\,q_v\,,\,p_v\,\big)\;=\;
\big(\,l^v_3-l^v_2\,,\,
l^v_1-l^v_2\,\big).
\end{equation}
Following Kashaev \cite{Ka1} we will consider the vector space
$V_{\vf}\simeq \BR^{4M}$
obtained by regarding the variables $q_v$, $p_v$ as the components
$q_v(\fv)$, $p_v(\fv)$ of vectors $\fv\in V_{\vf}$.
The space of linear coordinate functions on $V_{\vf}$ 
will be called the Kashaev space $W_{\vf}$. On $W_{\vf}$ 
we will consider the Poisson bracket $\Omega_{\vf}$  defined by
\begin{equation}\label{omegadef}
\Omega_{\vf}(p_v,q_w)\,=\,\de_{vw},\quad
\Omega_{\vf}(q_v,q_w)\,=\,0,\quad
\Omega_{\vf}(p_v,p_w)\,=\,0.
\end{equation}
The assignment \rf{trivars} associates a 
vector $\fv(P)$ in a subspace $T_{\vf}\subset V_{\vf}$ 
to each point $P\in\widetilde{\CT}(\Sigma)$. 
Kashaev has observed that
the subspace $T_{\vf}$ can be characterized by a suitable
set of linear forms $h_c\in W_{\vf}$ (``constraints''). 

\subsubsection{The constraints}\label{constraints}

To define the linear forms $h_c$ let us introduce an embedding  
of the first homology $H_1(\Sigma,\BR)$ into $W_{\vf}$ as follows.
Each graph geodesic  $g_c$ which 
represents an element $[\ga]\in H_1(\Sigma,\BR)$ may 
be described by an ordered sequence of vertices
$v_i\in\vf_\0$, and edges $e_{i}\in\vf_\1$,
$i=0,\dots,n$, where $v_0=v_n$, $e_0=e_n$, and we assume 
that $v_{i-1}$, $v_{i}$ 
are connected by the single edge $e_i$. We will define $\omega_i=1$
if the arcs connecting $e_i$ and $e_{i+1}$ turn around the 
vertex $v_i$ in the counterclockwise sense, $\omega_i=-1$ otherwise.
The edges emanating from $v_i$ will be numbered
$e^i_{j}$, $j=1,2,3$ according to the convention introduced in Figure
\ref{decor}. 
To each $[c]\in H_1(\Sigma,\BR)$ we will assign
\begin{equation}\label{constrdef}
h_{\ga}\df\sum_{i=1}^{n} u_i,\qquad 
u_i:\;=\;\omega_i\left\{ \begin{array}{ll} 
-q_{v_i} &\;\text{if $\{e_{i},e_{i+1}\}=\{e_3^i,e_1^i\}$,}\\ 
\phantom{-}p_{v_i}   &\;\text{if $\{e_{i},e_{i+1}\}=\{e_2^i,e_3^i\}$,}\\ 
q_{v_i}-p_{v_i} & \; \text{if $\{e_{i},e_{i+1}\}=\{e_1^i,e_2^i\}$.}
\end{array}\right. 
\end{equation}
$h_{\ga}$ is 
independent of the choice of representative $\ga$
within the class $[c]$. Let $C_{\vf}$ be the
subspace in $W_{\vf}$ that is spanned by the $h_{{c}}$, 
$[\ga]\in H_1(\Sigma,\BR)$.

\begin{lem} $\frac{\quad}{}$ \cite{Ka1} $\frac{\quad}{}$ \begin{itemize}
\item[(i)] The mapping $H_1(\Sigma,\BR)\ni [c]\mapsto h_{{c}}\in C_\vf$
is an isomorphism of vector spaces.
\item[(ii)]
The restriction of $\,\Omega_{\vf}$ to $C_{\vf}$ coincides with the
intersection form ${\rm I}$ on $H_1(\Sigma,\BR)$,
\[
\Omega_{\vf}(h_{c_\1},h_{c_\2})\;=\;{\rm I}(c_\1,c_\2).
\]
\item[(iii)]
The linear forms
$h_{\ga}$, $[\ga]\in H_1(\Sigma,\BR)$ vanish identically
on the subspace $T_{\vf}$.
\end{itemize}
\end{lem}
The equations $h_{\ga}(\fv)=0$, 
$[\ga]\in H_1(\Sigma,\BR)$ 
characterize the image of 
$\widetilde{\CT}(\Sigma)$ within $V_{\vf}$.
It is useful to recall that $H_1(\Sigma,\BR)$ splits as 
$H_1(\Sigma,\BR)=H_1(\Sigma_{\rm\sst cl},\BR)\oplus B(\Sigma)$,
where $B(\Sigma)$ is the $s-1$-dimensional 
subspace spanned by the homology classes
associated to the punctures of $\Sigma$, and 
$\Sigma_{\rm\sst cl}$ is the compact Riemann surface
which is obtained by ``filling'' the punctures of $\Sigma$.
The corresponding splitting of $C_{\vf}$ will be written as
$C_{\vf}=H_{\vf}\oplus B_{\vf}$.

\subsubsection{Change of fat graph\index{fat graphs}}

In order to describe the change of Kashaev 
variables\index{Kashaev coordinates} 
induced by a change of fat graph let us, following 
\cite{Ka1}, define the following two
transformations associated to the elementary moves 
$\omega_{vw}$ and $\rho_v$ respectively.
\begin{align}
A_v&: \quad(~q_v~,~p_v~) \,\mapsto \, (~p_v-q_v~,~-q_v~),\label{Aclassdef}\\
T_{vw}&:\quad\left\{\begin{aligned}
(~U_v~,~V_v~)& \,\mapsto \,\big(~U_vU_w~,~U_vV_w+V_v~\big),\label{Tdef}\\
(~U_w~,~V_w~)& \,\mapsto \,
\big(~U_wV_v(U_vV_w+V_v)^{-1}~,~V_w(U_vV_w+V_v)^{-1}~\big),
\end{aligned}\right.
\end{align}
where we have set $U_v\df e^{q_v}$ and $V_v\df e^{p_v}$ for all
$v\in\vf_\0$.
\begin{lem}\label{Ptcanon} The maps $A_v:W_\vf\ra W_{\rho_v\circ \vf}$
and $T_{vw}:W_\vf\ra W_{\omega_{vw}\circ \vf}$
defined in \rf{Aclassdef} and \rf{Tdef} respectively are 
canonical, i.e. they preserve the Poisson structure $\Omega_\vf$.
\end{lem}
The proof is again straightforward. 
Lemma \ref{Ptcanon} 
implies in particular that the mapping class 
group\index{mapping class group} acts on $W_{\vf}$ 
by {\it canonical} transformations.

\subsection{The structure of the Kashaev space $W_{\vf}$}

\subsubsection{Fock variables\index{Fock coordinates} 
vs. Kashaev's variables\index{Kashaev coordinates}}\label{S:kafock}

There is a canonical way to reconstruct the Fock-variables in terms of 
Kashaev's variables which is found
by combining equations \rf{penner-fock} and \rf{trivars}. 
The result may be formulated as follows. Let $v,w\in\vf_\0$
be the vertices that are connected by the edge $e\in\vf_\1$,
and let $e_i^v$, $i=1,2,3$ be the  edges introduced in Figure \ref{decor}.
\begin{equation}\label{kafock}
\hz_e=\hz_{e,v}+\hz_{e,w},\qquad \hz_{e,v}=\left\{ 
\begin{aligned}{}
p_v & \;\;{\rm if}\;\,e=e^v_1,\\
-q_v& \;\;{\rm if}\;\,e=e^v_2,\\
q_v-p_v & \;\;{\rm if}\;\,e=e^v_3.
\end{aligned}\right.
\end{equation}
The definition \rf{kafock} defines a linear map
$I_{\vf}:F_{\vf}\to \widehat{F}_\vf
\subset W_{\vf}$. It will be useful
to describe the properties of this map a bit more precisely.
\begin{lem}\label{Fockembed}
\begin{align*}{\rm(i)}\qquad& \hz_e(\fv(P))\;=\;z_e(P) \quad\forall 
e\in\vf_\1,\;\,
\forall P\in \TCT(\Sigma),\\
{\rm(ii)}\qquad&\Om_\vf(\hz_e,\hz_{f})\;=\;\Om_{\rm WP}(z_e,z_f)\quad\forall e,f\in\vf_\1,\\
{\rm(iii)}\qquad&\Om_\vf(\hz_e,h_{{c}})\;=\;0 \quad\forall e\in\vf_\1,\;\,
\forall c\in H_1(\Sigma,\BR),\\
{\rm(iv)} \qquad& \hat{f}_{c}\;\df\;I_{\vf}(f_c)
\;=\;h_{c},\quad \forall [c]\in B(\Sigma).
\end{align*}
\end{lem} 
\begin{proof} Straightforward verifications.
\end{proof}

It is also useful to remark
that the transformation of the Fock variables $\hat{z}_e$, $e\in\vf_\1$ 
that is induced by \rf{Aclassdef}, \rf{Tdef} coincides with
\rf{flipfovarclass}.

\subsubsection{Splitting of $W_\vf$}

The linear forms $h_{{c}}\in B^{}_{\vf}$ turn out to be 
the Hamiltonian generators for the symmetry \rf{gaugesymm} \cite{Ka1}.
It is therefore natural to consider the subspace $M_{\vf}\subset W_{\vf}$
which is spanned by the Hamiltonian vector fields that are generated by the
linear forms $h_{{c}}\in B_\vf$, as well as $N_{\vf}\df
M_{\vf}\oplus B_{\vf}$.

\begin{propn}\label{kashdecomp}
There exists a canonical transformation establishing the
isomorphism of Poisson vector spaces
$W_{\vf} \,\simeq\,T_{\vf} \oplus N_{\vf}\oplus H_\vf$,
such that \begin{itemize}\item[(i)]
$T_{\vf} \simeq \CT'(\Sigma)$ is the space of linear functions
on the Teich\-m\"uller space $\CT(\Sigma)$,
\item[(ii)] The restriction of $\Omega_\vf$ to $T_\vf$ coincides 
with the Poisson bracket induced by the Weil-Petersson symplectic form.
\end{itemize} 
\end{propn}
\begin{proof}
As a warmup it may be instructive to 
count dimensions: We have ${\rm dim}(W_\vf)=8g-8+4s$ and
${\rm dim}(C_\vf)={\rm dim}(H_1(\Sigma,\BR))=2g+s-1$. In order to determine
${\rm dim}(N_\vf)$ let us choose a canonical basis for 
$H_1(\Sigma,\BR)$, represented by curves $\al_1,\dots,\al_g,
\be_1,\dots,\be_g,\gamma_1,\dots,\gamma_{s-1}$ such that the only
nontrivial intersection pairings are ${\rm I}(\al_i,\be_j)=\de_{ij}$.
$N_{\vf}$ is spanned by the images of the classes
$[\gamma_1],\dots,[\gamma_{s-1}]$, together with the 
Hamiltonian vector fields that they generate.
It follows that ${\rm dim}(N_\vf)=2s-2$.

The main point that remains to be demonstrated is the 
existence of a decomposition of $\widehat{F}_{\vf}$ as
the direct sum 
\begin{equation} \label{decompFphi}
\widehat{F}_{\vf}=T_{\vf}\oplus B_{\vf}\quad {\rm such~~that}\quad
\Omega_{\vf}(t,h)=0\quad \forall~~ 
t\in T_\vf, \quad\forall~~h\in N_{\vf}.
\end{equation}
Thanks to Lemma \ref{Ptcanon} we may choose a convenient fat graph
to this aim. Let us pick a basis $\FB(\Sigma)$
for $B(\Sigma)$ represented by 
the curves which encircle $s-1$ of the $s$ punctures.
It is easy to see that we can always construct a
fat graph $\vf'$ such that the elements of $\FB(\Sigma)$
are represented by singles edges in $\vf_\1'$. These edges end in
a unique vertex $v(c)$ for each $c\in \FB(\Sigma)$.
It is clear that the expression for $h_c$ only involves the variables 
$(q_{v(c)},p_{v(c)})$ for all $c\in \FB(\Sigma)$. 
It follows from \rf{omegadef} that the Hamiltonian 
vector field generated by $h_c$ can likewise be
expressed in terms of 
$(q_{v(c)},p_{v(c)})$ only. The existence of the 
sought-for decomposition \rf{decompFphi} is obvious in 
this case. 

The result is carried over to the general
case with the help of Lemma \ref{Ptcanon}. 
It is clear that the subspace $T_{\vf}\subset\widehat{F}_{\vf}$  
is defined uniquely by the condition \rf{decompFphi}.
$T_\vf$ may then also be described as the quotient 
of $\widehat{F}_{\vf}$ by the conditions $h_c=0$ for all
$[c]\in B(\Sigma)$. It therefore 
follows from our discussion in \S\ref{Fockvars} that
$T_\vf$ is canonically isomorphic to $\CT'(\Sigma)$, 
the vector space of linear functions on $\CT(\Sigma)$.

To complete the proof it remains to observe 
that we have
\begin{equation}
\begin{aligned}{\rm (i)}\quad &
\Omega_{\vf}(t,h)=0\quad\forall~~ 
t\in T_\vf, \quad\forall~~h\in H_{\vf}, \\
{\rm (ii)}\quad &
\Omega_{\vf}(h,n)=0\quad\forall~~ 
h\in H_\vf, \quad\forall~~n\in N_{\vf}.
\end{aligned}
\end{equation}
(i) follows directly from Lemma \ref{Fockembed}, 
whereas part (ii) can easily be verified by considering the 
fat graph $\vf'$ above.
\end{proof}

\section{The Fenchel-Nielsen coordinates\index{Fenchel-Nielsen coordinates}}
\setcounter{equation}{0}

We will now be interested in the case of Riemann surfaces
with a boundary that is represented by a collection of $s>0$
geodesics. 
Another standard set of coordinates for the Teich\-m\"uller spaces
is associated to the decomposition of a Riemann surface into trinions
(three-holed spheres). We are now going to review the definition of these
coordinates. Different sets of  Fenchel-Nielsen coordinates will be
associated to different {\it markings}\index{marking} of the Riemann surface 
in a way which is analogous to the relation between the Penner
coordinates and fat graphs.

Let us denote by $S_\3$ the sphere with three holes (trinion). As a concrete model 
we may e.g. choose
\begin{equation}\label{standardSthree}
S_\3\;\df\;\{z\in\BC;|z|\geq\ep,|1-z|\geq\ep,|z|<1/\ep\}.
\end{equation}
Any trinion with a smooth boundary is 
diffeomorphic to $S_\3$.

\begin{defn}\label{markdefn}
A marking $\si$ of a surface $\Sigma$ consists of the following data.
\begin{itemize}
\item[(i)] A {\it cut system} $\CC_\si$, which is a
set $\CC_\si=\{c_1,\dots,c_{3g-3+s}\}$
of simple non-intersecting oriented closed 
curves $c_i$ on $\Sigma$. Cutting $\Sigma$ along the 
curves in $\CC$ decomposes the surface into a collection $\CP_\si$ 
of trinions. 
\item[(ii)] A choice of a trivalent 
graph $\Gamma_{T}$ with a single vertex $v_{\sst T}$
in each trinion $T\in\CP_\si$  such that the graphs on the different trinions 
glue to a connected graph $\Gamma_\si$ on $\Sigma$.
\item[(iii)] A choice of a distinguished boundary component for each 
trinion  $T\in\CP_\si$.
\end{itemize}
These data will be considered up to isotopy. 
\end{defn}
An example for the graphs $\Gamma_{T}$ is depicted in Figure \ref{marking}.

\begin{figure}[htb]
\epsfxsize3cm
\centerline{\epsfbox{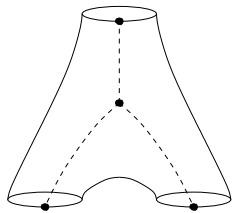}}
\caption{A trinion equipped with a marking graph}\label{markinggraph}
\end{figure}

\subsection{Definition of the Fenchel-Nielsen coordinates}\label{fencoor}

The basic observation underlying the definition of the 
Fenchel-Nielsen twist coordinates is the fact that for each triple
$(l_1,l_2,l_3)$ of positive real numbers there is a {\it unique}
metric of constant curvature $-1$ on the three-holed sphere (trinion)
such that the boundary components are geodesics\index{geodesic length} 
with lengths
$l_i$, $i=1,2,3$. A trinion with its metric of constant curvature $-1$
will be called hyperbolic trinion.
There furthermore exist three distinguished geodesics 
on each hyperbolic trinion
that connect the boundary components pairwise.

Let 
us assume  that the geodesic $c$ separates two 
trinions ${T}_{a}$ and ${T}_{b}$. Pick boundary components
$c_{a}$ and $c_{b}$ of ${T}_{a}$ and ${T}_{b}$ respectively
by starting at $c$, following the marking
graphs, and turning left at the vertices. As mentioned above,
there exist distinguished geodesics on ${T}_{a}$ and ${T}_{b}$ 
that connect $c$ with $c_{a}$ and $c_{b}$ respectively. 
Let $\de_c$ be the signed geodesic distance between the end-points
of these geodesics on $c$, and let
\begin{equation}
\theta_c\;=\;2\pi\frac{\de_c}{l_c}
\end{equation}
be the corresponding twist-angle. In a similar way one may define 
$\theta$ in the case that cutting along $c$ opens a handle.

Given a cut system 
$\CC=\{c_1,\dots,c_{\kappa}\}$, $\kappa=3g-3+s$, 
we thereby associate to each 
Riemann surface $\Sigma$ a tuple 
$(l_1,\dots,l_{\kappa};\theta_1,\dots,\theta_{\kappa})$ of real numbers.
It can be shown (see e.g. \cite{IT}) that the Riemann 
surface $\Sigma$ is characterized uniquely by the 
tuple $(l_1,\dots,l_{\kappa};e^{i\theta_1},\dots,e^{i\theta_{\kappa}})\in
(\BR^+)^{\kappa}\ti (S^1)^{\kappa}$. 
In order to describe the Teich\-m\"uller space $\CT(\Sigma)$ of {\it deformations}
of $\Sigma$ it suffices to allow for 
arbitrary {\it real} values of the twist angles $\theta_i$. 
Points in $\CT(\Sigma)$ are then parametrized by tuples
$(l_1,\dots,l_{\kappa};\theta_1,\dots,\theta_{\kappa})\in
(\BR^+)^{\kappa}\ti \BR^{\kappa}$.

\begin{rem}
The marking graph $\Ga_\si$ 
allows one to distinguish systems of 
Fenchel-Nielsen coordinates which are related to each other by 
Dehn twists\index{Dehn twist},
$\theta_c'=\theta_c+2\pi k_c$, $k\in\BZ$, $c\in\CC$. 
To use the markings for the parametrization of
different systems of Fenchel-Nielsen coordinates is 
then closely analogous to using
fat graphs for the specification of systems of 
Penner coordinates.

The definition of the Fenchel-Nielsen coordinates does not 
use the choice of a distinguished boundary component for 
each trinion. The latter has been included into the 
definition \ref{markdefn} for later convenience only.
\end{rem}

\subsection{Symplectic structure}

Let us furthermore notice that the 
Weil-Petersson symplectic form\index{Weil-Petersson form} 
becomes particularly simple in terms of
the Fenchel-Nielsen coordinates:
\begin{thm} \cite{Wo}$\;\frac{\quad}{}$
\begin{equation}
\omega_{\rm WP}^{}\;=\;
\sum_{i=1}^{\kappa}\,d\tau_i\land dl_i,\quad\tau_i=\frac{1}{2\pi}l_i\theta_i.
\end{equation}
\end{thm}
The content of the theorem may also be paraphrased as follows:
\begin{itemize}
\item[(i)] The geodesic length functions associated to non-intersecting 
closed curves Poisson-commute.
\item[(ii)] The Hamiltonian flows generated by the 
geodesic length functions coincide with the Fenchel-Nielsen twist flows.
\end{itemize}

\subsection{Geodesic lengths\index{geodesic length} 
from the Penner coordinates\index{Penner coordinates}}
\label{Fockleng}

A nice feature of the Fock coordinates is that they 
lead to a particularly simple way to reconstruct the Fuchsian
group\index{Fuchsian uniformization} 
corresponding to the point $P$ in Teich\-m\"uller space that is
parametrized by the variables $z_e(P)$, $e\in\vf_\1$. Assume given a 
graph geodesic $g_c$ on the fat graph homotopic to a simple closed
curve $c$ on $\Sigma$. Let the edges be labelled $e_i$, $i=1,\dots,r$
according to the order in which they appear on $g_c$, 
and define $\si_i$ to be $1$ if the path turns left\footnote{w.r.t.
to the orientation induced by the embedding of the fat-graph 
into the surface} 
at the vertex
that connects edges $e_i$ and $e_{i+1}$, and to be equal to $-1$ 
otherwise. The generator ${\rm X}({\ga})$ of the Fuchsian group 
that corresponds to $\ga $ is then constructed as follows \cite{Fo}.  
\begin{equation}\label{fuchsgen}
{\rm X}_{\ga}\;=\;{\rm V}^{\si_r}{\rm E}(z_{e_r})\dots {\rm V}^{\si_1}
{\rm E}(z_{e_1}),
\end{equation}
where the matrices ${\rm E}(z)$ and ${\rm V}$ are defined respectively by
\begin{equation}
{\rm E}(z)\;=\;\left(\begin{array}{cc} 0 & +e^{+\frac{z}{2}}\\
-e^{-\frac{z}{2}} & 0 \end{array}\right),\qquad
V\;=\;\left(\begin{array}{rr} 1 & 1 \\ -1 & 0 \end{array}\right).
\end{equation}
Given the generator ${\rm X}_{\ga}$ 
of the Fuchsian group one
may then calculate the hyperbolic length of the closed geodesic
isotopic to $c$ via
\begin{equation}\label{glength}
2\cosh\big(\fr{1}{2}l_{\ga}\big)\;=\;|{\rm tr}({\rm X}_{\ga})|.
\end{equation}

The proof of \rf{fuchsgen} was omitted in \cite{Fo}. 
We are therefore now going to explain how to verify the validity
of equation \rf{fuchsgen}.

It was remarked in the Subsection \S\ref{Fockvars} that for given values of the
coordinates $z_e$ one may construct a uniformized representation of the
corresponding Riemann surface by successively mapping ideal hyperbolic
triangles into the upper half plane which are glued according to the 
values $z_e$. Iterating this procedure ad infinitum one generates
a tessellation of the upper half plane by ideal hyperbolic
triangles. Let us now consider 
a generator ${\rm X}(\ga)$ of the Fuchsian group, represented 
on the upper half plane by a M\"obius 
transformation $M_{{\rm X}(\ga)}$, where 
$
M_{{\rm X}}(u)\; \df \;\frac{au+b}{cu+d}\;\;\text{if}\;\;
X=\left(\begin{smallmatrix} a & b \\ c & d\end{smallmatrix}\right).
$
The element $\ga\in\pi_1(\Sigma)$ 
may then be represented
by an open path on the upper half plane which leads from a chosen
base point $u$ to its image under $M_{{\rm X}}$. 

But one may 
equivalently represent the motion along the path by standing still
at the base point and moving the tessellation around by means of 
M\"obius-transformations. More precisely, let us assume 
that our base point $u$ is located within the 
ideal hyperbolic triangle $t_0$ with corners at $-1,0,\infty$, and that
the path $\CP_{\ga}$ 
representing our chosen element $\ga\in\pi_1(\Sigma)$ crosses
the edges $e_i$ $i=1,\dots, r$ in the order of the labelling. 
We may assume that the edge $e_1$ connects the points $0$
and $e^{z_{e_1}}$. After having crossed edge $e_1$ one would have 
left the triangle $t_0$ into the triangle $t_1$
with corners at $0,e^{z_{e_1}},\infty$. 
The M\"obius-transformation
$M_1$ corresponding to ${{\rm E}(z_{e_1})}$ 
brings one back into $t_0$: It can be checked
that it leaves the set of corners $\{-1,0, e^{z_{e_1}},\infty\}$ 
on the two 
triangles glued along $e_1$ invariant, but exchanges the two triangles.
To continue along the path $\CP_{\ga}$ in this fashion we now need
to map the next edge $e_2$ 
to be crossed to the edge going from $0$ to $\infty$
before we can apply the M\"obius-transformation corresponding to
${\rm E}(z_{e_2})$ in the same manner as before.
This is precisely what the M\"obius-transformation 
$M_{{\rm V}^{\si_1}}$ does: It simply rotates the edges of our 
fundamental triangle $t_0$. Moreover, 
the triangle $t_2$ that would be reached when leaving $t_0$
through $e_2$ will be mapped by $M_{{\rm V}^{\si_1}}$ into
the ideal hyperbolic triangle with corners $0, e^{z_{e_2}},\infty$. 
The fact that $e^{z_{e_2}}$ is indeed the position that the 
corner of $t_2$ is mapped into by $M_{{\rm V}^{\si_1}}$ follows
from the prescription for gluing $t_1$ and $t_2$ along $e_2$
in terms of $z_{e_2}$ and the fact that $M_{{\rm V}^{\si_1}}$
preserves cross-ratios.

By continuing in this fashion one generates the M\"obius 
transformation $M_{{\rm X}(\ga)}$ that evidently maps the original
tessellation representing the chosen point $P$ in $\CT(\Sigma)$ into
another one that is equally good as a representation for $P$.
By assumption, $M_{{\rm X}(\ga)}$ represents a closed path on the
considered fat graph. This means that the points that are
mapped into each other by $M_{{\rm X}(\ga)}$ are to be identified 
as different representatives for the same 
points on the surface corresponding to our point $P\in\CT(\Sigma)$.

\section{Coordinates for surfaces with holes of 
finite size}\label{coordholes}
\setcounter{equation}{0}

In the present paper we are mainly interested in the case of Riemann
surfaces which have a boundary $\pa\Sigma$
represented by $s$ geodesics of finite length. 
We therefore need to discuss 
how to introduce analogs of the previously described 
coordinate systems for $\CT(\Sigma)$ for the cases of interest here.
When considering surfaces with holes of finite size one has to
choose if one wants to keep the geodesic lengths of the 
boundary components variable, or if one wants to consider 
surfaces $\Sigma_\Lambda$ which have fixed boundary length 
given by the tuple $\Lambda=(l_\1,\dots,l_s)\in\BR_+^s$.
We shall find the first option often more convenient to work with.
Passing to a representation in which the boundary lengths are fixed will
then be almost trivial.

\subsection{Useful fat graphs\index{fat graphs} 
on surfaces with holes of finite size}
\label{normaltype}

Riemann surfaces $\Sigma$ with $s$ holes can always be represented 
by considering a Riemann surface $\Sigma^{\rm e}$ with $s$ pairs
of punctures, from which $\Sigma$ is obtained by cutting 
$\Sigma^{\rm e}$ along the geodesics $b_1,\dots,b_s$
that encircle the pairs of punctures. This simple observation
allows us to use the coordinates discussed 
previously in order to define coordinates for the Teich\-m\"uller
spaces of surfaces with $s$ holes. In order to spell out more 
precisely how to do this, let us first introduce a convenient 
class of fat graphs.

\begin{figure}[htb]
\epsfxsize5cm
\centerline{\epsfbox{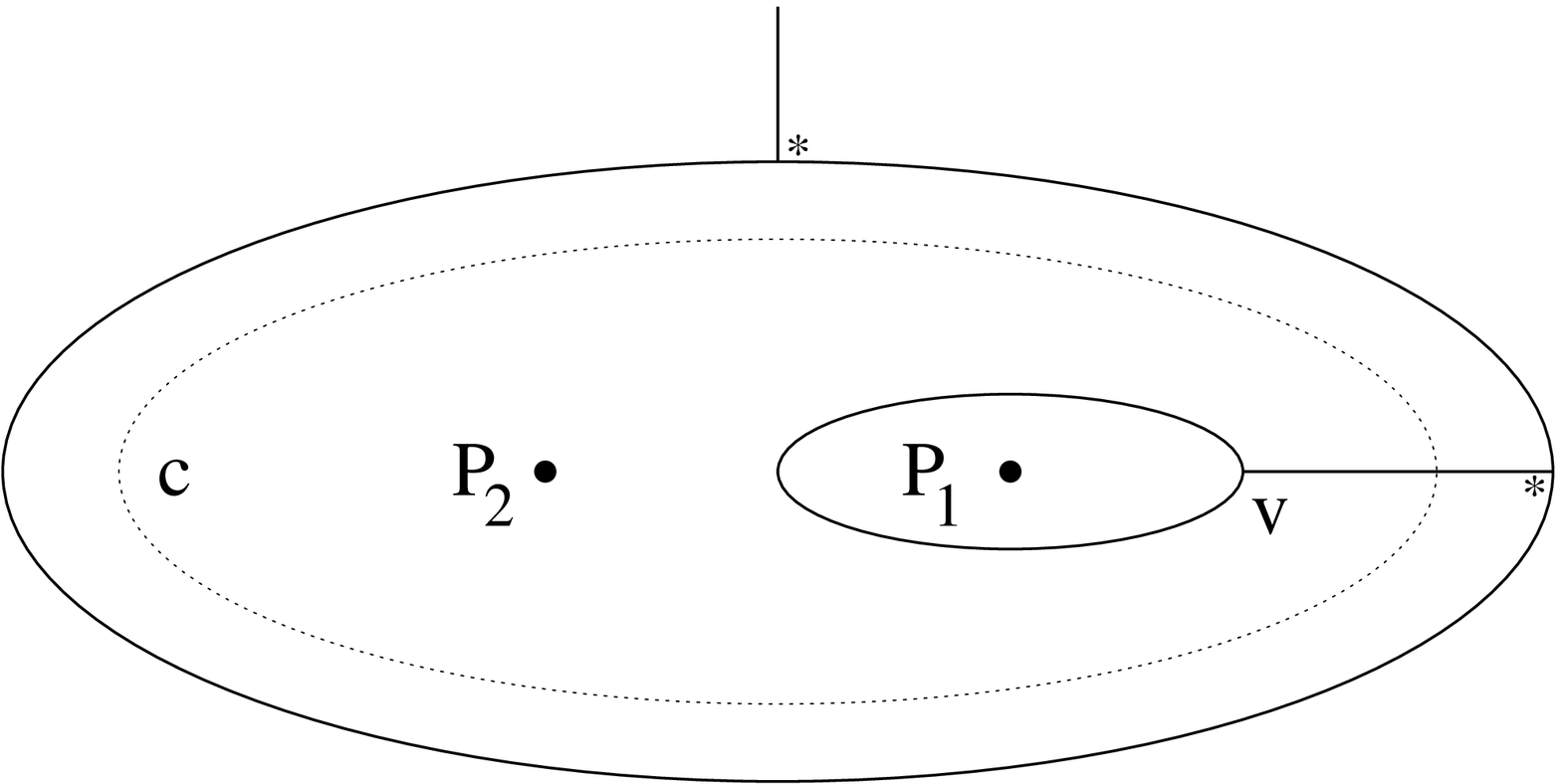}
\hspace{1cm}\epsfxsize5cm\epsfbox{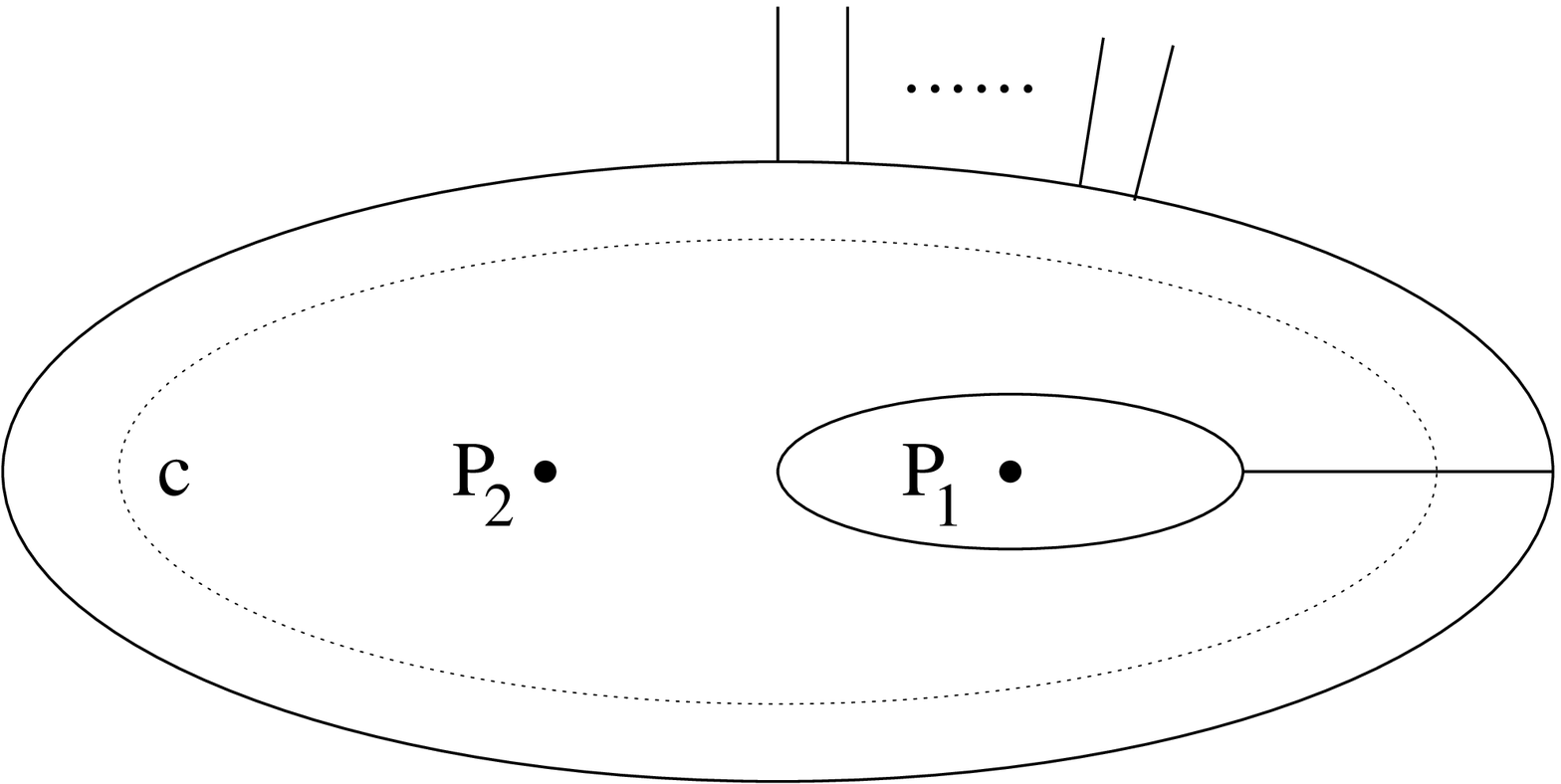}
}
\caption{Simple fat graphs in a neighborhood of the disc defined by the 
geodesic $c$ which encircles two punctures $P_1$ and $P_2$.}
\label{standard}
\end{figure}

Let us consider a pair of punctures  $(P_1,P_2)$. Let $c$ be a geodesic
such that cutting $\Sigma$ along $c$ produces two 
connected components one of which is a two-punctured disc $D$ with
punctures $P_1$ and $P_2$. A given fat graph $\vf^{\rm e}$ will
be said to have standard form near $D$  
if there exists a neighborhood of 
the disc $D$ in which $\vf^{\rm e}$ is homotopic to the fat graph
depicted on the left half of Figure \ref{standard}. 

For surfaces $\Sigma^{\rm e}$ with $s$ pairs of punctures
there exist fat graphs $\vf^{\rm e}$
which are of standard form in the neighborhood of $s-1$ discs $D_i$.
The simplest possible form of a fat graph around the remaining two
punctures is indicated on the right of Figure \ref{standard}. 
A fat graph $\vf^{\rm e}$ on a surface with $2s$ punctures
will be said to have standard form if it
has standard form near $s-1$  discs $D_i$, and if it has the 
form depicted on the right of Figure \ref{standard} in a neighborhood
of the remaining disc.

We will finally say that a fat graph $\vf$ on a Riemann surface 
$\Sigma$ with $s$ geodesic boundaries has standard form if 
$(\Sigma,\vf)$ can be obtained from a pair 
$(\Sigma^{\rm e},\vf^{\rm e})$ consisting of a 
$2s$-punctured surface $\Sigma$ and a fat graph $\vf^{\rm e}$ of standard
form by cutting $\Sigma^{\rm e}$
along $s$ geodesics $b_1,\dots,b_s$, each of which 
encircles a pair of punctures. The embedding $\Sigma\hookrightarrow 
\Sigma^{\rm e}$ furthermore induces an embedding 
${\rm MC}(\Sigma)\hookrightarrow {\rm MC}(\Sigma^{\rm e})$ of the
respective mapping class groups. The 
subgroup of ${\rm MC}(\Sigma^{\rm e})$ which is generated by
the diffeomorphisms that are supported on $\Sigma\subset\Sigma^{\rm e}$
preserves the set of fat graphs which have standard form.

\subsection{Kashaev type coordinates\index{Kashaev coordinates}}\label{Kashbound}

If we only use fat graphs of standard form, it becomes easy
to adapt the previously discussed systems of coordinates to 
the case of interest in the rest of this paper. We may in particular
consider the space $W_{\vf}$ of Kashaev variables associated to the 
fat graph $\vf$. A subspace $T_\vf$ of $W_{\vf}$ can again be defined
by means of the decomposition \rf{decompFphi}. It is furthermore
convenient to introduce the set $\vf_\1'$ which only contains the 
edges of $\vf$ that do not end in boundary components of $\Sigma$.

\begin{lem} \label{TTlem} \begin{itemize}
\item[(i)]
We have $T_{\vf}^{}\simeq T_{\vf^{\rm e}}$.
\item[(ii)] 
The Fock coordinates 
$\{z_e;e\in\vf_\1'\}$ form a set of coordinates for $T_\vf$.
\end{itemize}
\end{lem}
\begin{proof}
Let us first note that the linear form 
$h_1$ which is  via \rf{constrdef} associated to 
puncture $P_1$ can be expressed in terms of the variables $(p_v,q_v)$
associated to the vertex $v$ in Figure \ref{standard}
{\it only}. This means that both $p_v$ and $q_v$ 
are contained in $N_{\vf^{\rm e}}$. Instead of the linear form
$h_2$ associated to puncture $P_2$ we may 
consider $h_{c}=h_1+h_2$, which can be expressed exclusively
in terms of the variables associated to the vertices contained
in $\Sigma$. Part (i) of the lemma follows easily from these observations.

In order to verify part (ii) one may again consider $\vf^{\rm e}$.
When writing the 
relations $f_c=0$, $c\in B(\Sigma)$ in terms of the Fock variables
$z_e$, $e\in\vf^{\rm e}_\1$ one will always find contributions
containing the $z_e$, $e\in\vf^{\rm e}_\1\setminus\vf_\1'$. It is then
easy to convince oneself that the  relations $\hat{f}_c=0$ may be used to
express the $z_e$, $e\in\vf^{\rm e}_\1\setminus\vf_\1'$ in terms
of the $z_e$, $e\in\vf_\1'$. After this is done, all relations
$f_c=0$, $c\in B(\Sigma)$ are satisfied, the variables
$z_e$, $e\in\vf_\1'$ are therefore unconstrained.
\end{proof}  

However, in this case
the relation between the Teich\-m\"uller space 
$\CT(\Sigma)$ and 
$T_\vf$ is slightly more complicated. In order to describe this relation
let us consider the spaces ${\rm Fun}(\CT(\Sigma))$ 
and ${\rm Fun}(T_\vf)$ of smooth functions 
on $\CT(\Sigma)$ and $T_\vf$ respectively. These spaces carry 
canonical Poisson brackets $\{\,.\,,\,.\,\}_{\rm \sst WP}$ and
$\{\,.\,,\,.\,\}_{\vf}$
uniquely defined by the bilinear 
forms $\Omega_{\rm \sst WP}$ and $\Omega_{\vf}$ respectively. 

\begin{propn}\label{hatTT}
We may represent ${\rm Fun}(\CT(\Sigma))$ as the subspace of 
${\rm Fun}(T_\vf)$ which is defined by the conditions
\begin{equation}
\{\,F\,,\,l_i\,\}_{\vf}\;=\;0, \qquad i=1,\dots,s,\;\;F\in{\rm Fun}(T_\vf),
\end{equation}
where $l_i$ is the length function which is associated to the i-th
boundary component via equations \rf{fuchsgen} and \rf{glength}.
\end{propn}

It should be noted that the length functions associated to the 
boundary components are contained in ${\rm Fun}(\CT(\Sigma))$.

\newpage
\part{Quantization of the Teich\-m\"uller 
spaces\index{quantization of Teichm\"uller spaces}}

Our aim is to construct certain classes of infinite-dimensional 
representations of the mapping class groups ${\rm MC}(\Sigma)$. One possible
approach to this problem is to ``quantize'' Poisson
manifolds like the Teichm\"uller spaces 
on which ${\rm MC}(\Sigma)$ acts as a group of symmetries. 

Our construction will proceed in two main  steps. Quantization of the
Kashaev spaces $W_\vf$ leads to a rather elegant 
construction of projective unitary  representations
of the mapping class groups \cite{Ka1}. However, these representations 
turn out to be reducible. The second step
will therefore be to identify distinguished subrepresentations
within the representations coming from the quantization of the
Kashaev spaces as the mapping class group representations which
are naturally associated to the quantization of the Teichm\"uller spaces.
A direct construction of the latter is not known, 
which is why this somewhat indirect construction seems to be most
efficient at the moment.

\section{Quantization of the Teichm\"uller spaces}\label{Quantsimple}

\subsection{Canonical quantization}\label{CanQ}

Quantization of a Poisson manifold $\CP$ 
means ``deforming'' the space of functions on $\CP$ 
into a one-parameter ($\hbar$) family of noncommutative
algebras $\CP_{\hbar}$  in such a way that 
the deformed 
product $f\ast_\hbar^{} g$ satisfies
\begin{equation}
f\ast_\hbar^{} g\,=\,fg+\hbar \Om(f,g)+\CO(\hbar^2),
\end{equation} 
where $fg$ is the ordinary commutative product of functions on $\CP$ and
$\Om(f,g)$ is the Poisson bracket on $\CP$. 
If the Poisson manifold $\CP$ has a group $G$ of symmetries
it is natural to demand that these symmetries
are preserved by quantization in the sense that 
any $g\in G$ is realized as an automorphism 
$f\ra\sa_g(f)$ of $\CP_{\hbar}$. 

Representations of the group $G$ can be constructed by studying
representations of the algebra $\CP_\hbar$ by
operators $\SO(f)$ on a Hilbert space $\CH$,
\[
\SO(f)\cdot\SO(g)\,=\,\SO(f\ast_\hbar^{} g)\,.
\] 
The Hilbert space $\CH$
will then typically come equipped 
with a unitary projective 
representation of the group $G$ of symmetries by operators
$\SU_g$ such that the automorphisms $\SA_g(\SO(f))\equiv\SO(\sa_g(f))$
are realized as
\[
\SA_g(\SO)\,=\,\SU_g\cdot\SO\cdot\SU_g^{-1}\,.
\]

Quantization is particularly simple if there exist coordinate
functions 
$q_1,\dots,q_N$ and $p_1,\dots,p_N$ defined globally
on $\CP$ such that the Poisson bracket 
takes the form
\begin{equation}\label{canform}
\Omega(p_v,q_w)\,=\,\de_{vw},\quad
\Omega(q_v,q_w)\,=\,0,\quad
\Omega(p_v,p_w)\,=\,0,
\end{equation}
for $v,w\in\{1,\dots,N\}$. One may then define $\CP_\hbar$ in such a way
that the relations
\begin{equation}
i[p_v,q_w]\,=\,\hbar\de_{vw},\quad
[q_v,q_w]\,=\,0,\quad
[p_v,p_w]\,=\,0,
\end{equation}
hold for $[f,g]\equiv f\ast_\hbar^{}g-g\ast_\hbar^{}f$. There exists 
a standard representation $\SO$ of these commutation relations
on (dense subspaces of) 
the Hilbert space $\CH=L^2(\BR^N)$ of square-integrable functions 
$\Psi(\fq)$,
$\fq=(q_1,\dots,q_N)$, which is generated by pairs of
operators 
\[
\spp_v\,\equiv\,\frac{1}{2\pi\hbar}  \SO(p_v),\qquad
\sq_v\,\equiv\,\SO(q_v)\,,\qquad v=1,\dots,N\,,
\]
that are defined respectively by
\begin{equation}\label{spsqdef}
\sq_v\Psi(\fq)\,\equiv\,q_v\Psi(\fq)\,,\qquad
\spp_v\Psi(\fq)\,\equiv\,\frac{1}{2\pi i}\frac{\pa}{\pa q_v}\Psi(\fq)\,.
\end{equation}
This simple example for the quantization of Poisson manifolds 
is often referred to as ``canonical quantization''.

\begin{rem} \label{irredrem}
The representation that is constructed in this way is irreducible
in the following sense. If $\SO$ is a bounded  operator on $L^2(\BR^N)$
which commutes\footnote{Commutativity $[\SO,\SA]=0$ with a self-adjoint 
unbounded operator $\SA$ is, by convention, understood 
in the sense of commutativity with the spectral projections of $\SA$.
For the reader's convenience, we have collected the 
relevant operator-theoretical results in 
Appendix B.} 
with the operators $\spp_v$ and $\sq_v$
for all $v=1,\dots,N$ then $\SO=\chi$, the operator of multiplication 
with the complex number $\chi$.
\end{rem}

Following the discussion in Sections \ref{class},
it seems natural to define the quantized Teichm\"uller
spaces as the noncommutative algebra $\CT_\hbar(\Sigma)$ with
generators $z_e$ and relations 
 \begin{align}
{\rm(i)}\qquad&i\,[\;z_{e_\2}\, , \, z_{e_\1}\,]\;=\;\hbar\, 
\Om_{\rm\sst WP}^{}(\,z_{e_\2}\, , \, z_{e_\1}\,)\label{Qfockrel}\\
{\rm(ii)}\qquad& f_{\vf,c}\,=\,0\,, \;\;
\forall\;c\in B(\Sigma).\label{Qfockconstr}
\end{align}
The space $\CH(\Sigma)$ will be defined as an irreducible 
representation of the commutation relations \rf{Qfockrel}
which satisfies the additional conditions \rf{Qfockconstr}.

\subsection{Quantization of the Kashaev space $W_\vf$}\label{quantsect}


For each given fat graph $\vf$ let us define the 
Hilbert space $\CK(\vf)$ as the space of square integrable functions
$\Psi(\fq)$ of the Kashaev-variables $\fq=(q_1,\dots,q_{2M})$.
On $\CK(\vf)$ we shall consider the basic operators $\spp_v$, $\sq_v$
defined in \rf{spsqdef} for 
$v=1,\dots,2M$. The noncommutative algebras of operators
which are generated by the operators $\sq_v$, $\spp_w$ with the
commutation relations
\begin{equation}
[\,\spp_v\,,\sq_w\,]\,=\,(2\pi i)^{-1}{\de_{v,w}}\,
\end{equation}
may be considered as representing quantized algebras of 
functions on 
the Kashaev space $W_\vf$.

\begin{rem}
It is worth noting that the decoration of the triangles
is used to define the concrete realization of the space 
$\CK(\vf)$ as a space of square integrable functions.
\end{rem}

When quantizing the Kashaev spaces $W_\vf$ we get more than we
ultimately want. In order to see this,
let us  introduce quantum analogs of the  coordinate
functions $h_c$ and $\hat{z}_e$ respectively, i.e. self-adjoint operators
$\hat\sh_c$ and $\hat\sz_e$ on $\CH(\vf)$ which are defined 
by formulae very similar 
to \rf{constrdef} and \rf{kafock} respectively, normalized in such a way that
the following commutation relation hold:
 \begin{align}
{\rm(i)}\qquad& i[\,\hat\sz_{e_\1}\, , \, \hat\sz_{e_\2}\,]\;=\;\hbar\, 
\Omega_{\rm\sst WP}(z_{e_\1},z_{e_\2}),\label{zz}\\
{\rm (ii)}\qquad &i[\,\hat\sh_{c_\1}\, ,\, \hat\sh_{c_\2}\,]\,=\,
 \hbar \,{\rm I}(c_\1,c_\2).
\label{hh}\\
{\rm(iii)}\qquad&[\,\hat\sz_e\, , \, \hat\sh_{{c}}\,]\;=\;0 \quad\forall c\in 
H_1(\Sigma,\BR).\label{zh} 
\end{align}
We observe that we {\it do} have a representation, henceforth denoted 
$\SZ_\vf$, of the 
algebra \rf{Qfockrel}, but this representation is neither 
irreducible, nor does it fulfill the additional relations
\rf{Qfockconstr}. 
The latter point becomes most clear if 
one introduces the operators $\hat{\sf}_c$, $c\in B(\Sigma)$
associated to the relations \rf{Qfockconstr} which are defined
by replacing $z_e\ra \hat{\sz}_e$ in \rf{Fockconstrdef}.
We may then observe that  
\begin{equation}\label{sfshrel}
\hat{\sf}_c=\sh_c\quad{\rm for}\;\; c\in B(\Sigma)\,, 
\end{equation} 
which is verified in the same way as statement (iv) in
Lemma \ref{Fockembed}.

\subsection{Reduction to the quantized Teichm\"uller spaces}\label{Qredsubsec}

In order to see that the representations
$\SZ_{\vf}$ ``contain'' irreducible  representations of 
the quantized Teichm\"uller spaces $\CT_\hbar(\Sigma)$, 
let us  consider the noncommutative algebra
$\CZ_\hbar(\Sigma)$ 
with generators $z_e$ and the only relations \rf{Qfockrel}. 
One should observe that
the $f_c$ generate the center of the algebra $\CZ_\hbar$. Irreducible 
unitary representations $\SZ_\ff$ 
of this algebra are parametrized by linear functions $\ff:B(\Sigma)\ra\BR$.
Such representations are such that the operators 
${\sf}_c=\SZ_\ff(f_c)$ are realized
as the operators of multiplication with the  real numbers $\ff(c)$, 
$c\in B(\Sigma)$. 

The representations can be constructed concretely by forming linear 
combinations $t_k^{}$ and $t_{k'}^{\sst\vee}$, $k=1,\dots,3g-3+s$
of the $z_e$  which \begin{itemize}
\item[(i)] are mutually linearly independent,
and linearly independent of the $f_c$, 
\item[(ii)] and which satisfy the commutation relations 
\begin{equation}
2\pi \,
\big[\,t_k^{}\,,\,t_{k'}^{\sst\vee}\,\big]\,=\,i\,\de_{kk'}\,.
\end{equation}
\end{itemize}
Canonical quantization realizes $\SZ_\ff$ on the Hilbert space 
$\CH_{\ff}(\vf)\simeq L^2(\BR^{3g-3+s})$ which consists of square-integrable
functions $\Phi(\ft)$, $\ft=(t_1,\dots,t_{3g-3+s})$. 
The operators $\sz_{\ff,e}\equiv \SZ_\ff(z_e)$ are then realized 
as linear combination of the operators 
\[
\mst_k^{}\Phi_\ff^{}(\ft)=t_k^{}\Phi_\ff^{}(\ft),\quad
\mst_{k}^{\sst\vee}\Phi_\ff^{}(\ft)=\frac{1}{2\pi i}\frac{\pa}{\pa t_k} 
\Phi_\ff^{}(\ft),\quad
\sf_c^{}\Phi_\ff^{}(\ft)=\ff(c)\Phi_\ff^{}(\ft),
\]
where $k=1,\dots,3g-3+s$ and $c\in B(\Sigma)$.


Our aim is to describe how the representation $\SZ_\vf$ decomposes
into the representations $\SZ_\ff$.
In order to do this, 
let us introduce the representation
\[
{\SZ}'_\vf\,\equiv\int_{B'(\Sigma)}^{\oplus}d\ff \;\SZ_{\ff}\quad
\text{on the space}\quad
{\CH}_z(\vf)\equiv\int_{B'(\Sigma)}^{\oplus}d\ff \;\CH_\ff(\vf)\,.
\]
The space ${\CH}_z(\vf)$ is spanned by
square-integrable families 
$\Phi\equiv(\Phi_\ff)_{\ff\in B'(\Sigma)}$ of 
functions $\Phi_\ff\in\CH_\ff(\vf)$ which are associated to the 
linear functions $\ff:B(\Sigma)\ra \BR$ 
in the dual $B'(\Sigma)\simeq\BR^{s-1}$ up to a set of measure zero.
The 
representatives 
$\sz_e\equiv\SZ_\vf(z_e)$ are defined as follows
\begin{equation}
\sz_e\,\Phi\,\equiv\,\big(\,\sz_{\ff,e}\Phi_{\ff}\,\big)_{\ff\in B'(\Sigma)\;.}
\end{equation}

\begin{propn} \label{Qkashdecomp} The decomposition of 
the representation $\SZ_\vf$ into irreducible 
representations of $\CZ_\hbar$ may be written as follows:
\begin{equation}
\SZ_\vf\,\simeq\,\left(\,\int_{B'(\Sigma)}^{\oplus}d\ff \;\SZ_{\ff}\right)
\ot 1_{\CH_h(\vf)}^{}\,,
\end{equation}
where the space $\CH_h(\vf)$ is isomorphic to $L^2(\BR^{g})$.
There exists a unitary operator $\SI_\vf$,
$\SI_\vf:\CK(\vf)\ra \CH_z(\vf)\ot\CH_h(\vf)$ such that
\begin{equation}\label{ZHmap} 
\SI_\vf^{}\cdot\hat{\sz}_e^{}\cdot\SI_\vf^{-1}=\sz_e^{}\ot 1\quad{\rm and}
\quad\SI_\vf^{}\cdot\hat{\sh}_c^{}\cdot\SI_\vf^{-1}=1\ot\sh_c^{}\,,
\end{equation}
for any $e\in\vf_\1$ and $c\in H_1(\Sigma_{\rm\sst cl},\BR)$, respectively.
\end{propn}
\begin{proof}
Let us recall 
the direct sum decomposition
\begin{equation}\label{kashdecomp2}
W_{\vf} \,\simeq\,T_{\vf} \oplus N_{\vf}\oplus H_\vf\,.
\end{equation}
To each of the three spaces $T_\vf$, $N_\vf$, $H_\vf$ 
one may choose coordinates which bring the Poisson bracket
to the canonical form \rf{canform}. 
The  corresponding operators
\begin{align*}
&(\,\hat\mst_k^{}\,,\,\hat\mst_k^{\sst \vee}\,)\,, \quad k=1,\dots,3g-3+s\,,\\
&(\,\hat\sf_l^{}\,,\,\hat\sf_l^{\sst\vee}\,)\,, \quad l=1,\dots,s-1\,,\\
&(\,\hat\sh_m^{}\,,\,\hat\sh_m^{\sst\vee}\,)\,, \quad m=1,\dots,g\,,
\end{align*}
can be constructed as linear combinations of the $(\spp_v,\sq_v)$, 
$v\in\vf_\1$, in such a way that the only nontrivial commutation 
relations are
\begin{equation}\label{ccrfactor}
\begin{aligned}
 &2\pi \,
\big[\,\hat\mst_k^{}\,,\hat\mst_{k'}^{\sst\vee}\,\big]\,=\,ib^2\de_{kk'}, 
 \quad k,k'=1,\dots,3g-3+s\,,\\
&2\pi \,\big[\,\hat\sf_l^{}\,,\,\hat\sf_{l'}^{\sst\vee}\,\big]\,=\,
ib^2\de_{ll'}, \quad l,l'=1,\dots,s-1\,,\\
&2\pi \,\big[\,\hat\sh_m^{}\,,\,\hat\sh_{m'}^{\sst\vee}\,\big]\,=\,ib^2
\de_{mm'}, \quad m,m'=1,\dots,g\,.
\end{aligned}
\end{equation}
It will be convenient to form the following vectors with
$4M$ operator-valued components:
\begin{align*}
\sv &=(\dots ,\sq_v,\dots,\spp_w,\dots)\,,\\
\hat{\sv} &=(\dots ,\hat\mst_k^{},\dots,\hat\sf_l^{},\dots,
\hat\sh_m^{},\dots,\hat\mst_k^{\sst\vee},\dots,\hat\sf_l^{\sst\vee},\dots,
\hat\sh_m^{\sst\vee},\dots)\,.
\end{align*}
The linear change of variables $\hat{\sv}=\hat{\sv}(\sv)$ 
can then be represented by 
a symplectic $(4M\times 4M)$-matrix $J_\vf$,
\begin{equation}\label{Jphi}
\hat{\sv}\,=\,J_{\vf}\sv\,,\qquad J_\vf\in Sp(2M,\BR).
\end{equation}

On the other hand let us note that $\CH_z(\vf)\ot\CH_h(\vf)$ is canonically
isomorphic to $L^2(\BR^{2M})$ via
\begin{equation}
\SK:\Phi\ot\psi\ra \Psi,\qquad \Psi(\ft,\ff,\fh)\,\equiv\,
\Phi_\ff(\ft)\psi(\fh)\,.
\end{equation}
The corresponding  representation of the commutation
relations \rf{ccrfactor} on the Hilbert space $L^2(\BR^{2M})$ 
is obtained by renaming the operators $q_v,p_v$ as follows
\begin{align*}
(\dots ,\sq_v,\dots, &\spp_w,\dots)\\
    &\equiv\,(\dots ,\mst_k^{},\dots,\sf_l^{},\dots,
\sh_m^{},\dots, \mst_k^{\sst\vee},\dots,\sf_l^{\sst\vee},\dots,
\sh_m^{\sst\vee},\dots)\,.
\end{align*}
It follows from the Stone - von Neumann uniqueness theorem for the 
representation of the commutation relations \rf{ccrfactor} that
these two representation must be related by a unitary transformation.
This transformation may be characterized more precisely as follows.

\begin{lem} $\quad$ \label{canqlem}
\begin{itemize}
\item[a)] 
To each $\gamma\in Sp(2M,\BR)$ there exists a unitary 
operator $\SJ_\gamma$ on $L^2(\BR^{2M})$ such that 
\begin{equation}\label{Jtrsf}
\SJ_{\gamma}^{}\cdot \sv\cdot \SJ_\gamma^{-1}\,=\,\gamma\sv\,.
\end{equation}
The operators $\SJ_\gamma$ generate a projective unitary 
representation of $Sp(2M,\BR)$.
\item[b)]
The operators $\SJ_\gamma$ can be represented in  the form 
\begin{equation}
\SJ_\gamma\,=\,\exp\big(i J_{\gamma}(\sv)\big)\,,
\end{equation}
where $J_{\gamma}(\sv)$ is a quadratic expression in the operators
$\sv$. 
\end{itemize}
\end{lem}
\begin{proof} Part a) is a classical result of I.E. Segal, \cite{Se}.
Part b) follows easily from the observation that 
the quadratic functions of the operators $\sv$
generate a representation
of the Lie algebra of $Sp(2M,\BR)$ which satisfies 
the infinitesimal version of \rf{Jtrsf}, 
see e.g. \cite{GS} for
more details. 
\end{proof}

One may therefore find an operator $\SJ_\vf$ on $L^2(\BR^{2M})$
which represents the transformation \rf{Jphi} in the sense that
\begin{equation}
\SJ_{\vf}^{}\cdot \sv_{}\cdot \SJ_\vf^{-1}\,=\,
\hat{\sv}_{}(\sv)\,=\,J_{\vf}^{}\sv\,.
\end{equation}
The sought-for isomorphism $\SI_\vf$  can finally be constructed
as $\SI_\vf^{}=\SK^{-1}\cdot\SJ_{\vf}^{-1}$.
\end{proof}

\begin{rem} 
It is worth noting that the definition of $\CH_\ff(\vf)$ depends
only on the combinatorial structure of the fat graph $\vf$, not on the
way it is embedded into the  Riemann surface $\Sigma$. It follows
that the isomorphism $\CH_\ff(\mu.\vf)\simeq \CH_\ff(\vf)$, 
$\mu\in {\rm MC}(\Sigma)$ is {\it canonical}.
\end{rem}

\section{Representations of the mapping class groups}

The representations of the mapping class group associated to the 
quantized Teichm\"uller spaces will be obtained by means of a
very general construction which produces representations
of the group $G$ of symmetries of a two-dimensional CW complex $\CG$ 
out of representations of the edge path groupoid of $\CG$. 
We will first describe this construction, before
we discuss how to construct representations of the Ptolemy groupoid
on the quantized Teichm\"uller spaces. The latter will then induce the 
sought-for representation of ${\rm MC}(\Sigma)$.

\subsection{Projective unitary representations of 
groupoids\index{representations of groupoids}}\label{projgroup}

Let us recall that a groupoid ${\rm G}$
is a category such that all morphisms are invertible.
The objects of ${\rm G}$ will here be denoted by letters
$U,V,W,...$. Anticipating that the groupoids ${\rm G}$ 
we will be interested in are path groupoids of some topological
space we will use the notation 
$[W,V]\df {\rm Hom}_{\rm G}(V,W)$. The elements of $[W,V]$ will
also be called "paths".  

\begin{defn}\label{groupoidrepdef}
A unitary projective
representation of the groupoid ${\rm G}$ consists of the following data:
\begin{itemize}
\item[(i)] A  Hilbert space $\CH(V)$ associated to each object 
$V\in{\rm Ob}({\rm G})$,
\item[(ii)] a map $\su$ which associates to 
each path $\pi\in[W,V]$ in ${\rm G}$
a unitary operator \[
\su(\pi)~~:~~\CH(V)~\ra~\CH(W)~,\]
\item[(iii)] a family of maps $\zeta_V$, $V\in{\rm Ob}({\rm G})$
which associate to each {\rm closed} path
$\pi\in[V,V]$ a number  $\zeta_V(\pi)\in\BC$ with $|\zeta(\pi)|=1$.
\end{itemize}
These data are required to satisfy the relations
 \begin{align}\label{grpoidrels1}{}& 
 \begin{aligned}
 &{\rm a)} \quad
 \su(\pi_\2\circ\pi_\1)\;=\;\su(\pi_\2)\cdot\su(\pi_\1)\,,\\
 &{\rm c)} \quad\zeta_V(\pi_\2\circ\pi_\1)\;=\;\zeta_V(\pi_\2)\zeta_V(\pi_\1)\,,\\
&{\rm e)}  \quad
 \su(\pi)\;=\;\zeta_V^{}(\pi)\;\;{\rm if}\;\;\pi\in[V,V]\,,
 \end{aligned}
 \qquad \quad
 \begin{aligned}
 &{\rm b)} \quad
 \su(\pi^{-1})\;=\;\su^{\dagger}(\pi)\,,\\
&{\rm d)} \quad
 \zeta_V(\pi^{-1})\;=\;(\zeta_V(\pi))^*\,,\\
&{\rm f)} \quad
 \su(\id)\;=\;1,\\
 \end{aligned}\end{align}
where we use the notation $\zeta_V^{}(\pi)$ also to denote the operator
which multiplies each vector of $\CH_V$ by the 
 number $\zeta_V(\pi)$.
\end{defn}

The groupoids of interest will be the path groupoids 
${\rm G}$ of two-dimensional CW complexes $\CG$. The set of objects
is given by the set of vertices $\CG_\0$, whereas the set of morphisms
coincides with the set of paths in the complex $\CG$. 
Since each path $\pi$ may be represented as a chain 
$E_{\pi,n(\pi)}\circ\dots E_{\pi,2}\circ E_{\pi,1}$ 
of edges in $\CG_\1$
it is clear that a projective unitary representation of the path groupoid 
${\rm G}$ of a two-dimensional CW complex $\CG$ is characterized completely
by specifying the images $\su(E)$ for $E\in\CG_\1$. 
Existence of the family of maps $\zeta_V$ such that 
relation e) is fulfilled represents a rather  nontrivial 
constraint that the operators $\su(E)$, $E\in\CG_\1$ have to
satisfy. Of course it suffices to satisfy these constraints for the
2-cells $\pi\in\CG_\2$.


\subsection{Representations of symmetries of a groupoid} \label{symmG}

The group of symmetries $G$  of a two-dimensional CW complex $\CG$ 
is the group of all invertible mappings 
\[
\mu\;:\;\left\{
\begin{aligned}
 \CG_\0\ni V\;&\longrightarrow\;\mu.V\in\CG_\0\,,\\
 \CG_\1\ni E\;&\longrightarrow\;\mu.E\in\CG_\1\,.
\end{aligned}\right.
\]
There is an associated action on the edge paths in the complex $\CG$,
\[
[W,V]\ni \pi\;\ra\; \mu.\pi\in[\mu. W,\mu. V],
\]
which is such that 
\begin{equation}
\mu(\pi_\2\circ \pi_\1)\,=\,\mu(\pi_\2)\circ\mu(\pi_\1), \qquad
\mu(\pi^{-1})\,=\,(\mu(\pi))^{-1}.
\end{equation}
We will assume that
we are given 
a unitary projective representation of ${\rm G}$ which is compatible 
with the symmetry $G$ in the sense that $\CH(V)$ is canonically
isomorphic with $\CH(\mu.V)$, $\CH(V)\simeq\CH(\mu.V)$. 
We are going to show that the given representation of the groupoid 
${\rm G}$ canonically induces a representation of 
its group $G$
of symmetries. 

Let us fix a base point $V\in{\rm Ob}({\rm G})$
and assume having chosen a path
$\pi_{\sst V}(\mu)\in[\mu. V,V]$ for each $\mu\in G$.
Let then $\SR(\mu):\CH(V)\ra\CH(V)\simeq\CH(\mu. V)$ be defined by
\begin{equation}
\SR_V(\mu)\;=\;\su(\pi_{\sst V}(\mu))\,.
\end{equation} 
We are going to assume that the paths $\pi_{\sst \mu_\1.V}(\mu_\2)$ 
are the translates of $\pi_{\sst V}(\mu_\2)$ under $\mu_\1$, i.e. that
$\pi_{\sst \mu_\1.V}(\mu_\2)=\mu_1.\pi_{\sst V}(\mu_\2)$. 
It follows that
\[
\SR_{\mu_\1.V}(\mu_\2)=\SR_{V}(\mu_\2)\,.
\]
The operators $\SR_V(\mu)$ satisfy the relations
\begin{align}
{}& \SR_{V}(\mu_\2)\cdot\SR_{V}(\mu_\1)
\,=\,\vartheta_V^{}(\mu_\2,\mu_\1)\,\SR_V(\mu_\2\circ\mu_\1),
\label{projgrouplaw}\\
& \vartheta_V^{}(\mu_\2,\mu_\1)\,=\,\zeta_V
\big(\pi_{V}^{-1}(\mu_\2\circ\mu_\1)\circ
\pi_{\mu_\1. V}(\mu_\2)\circ\pi_{V}(\mu_\1)\big).
\end{align}

We may next observe that the apparent dependence 
on the base point $V\in\CG_\0$ 
is inessential. Let $V,W\in\CG_\0$, and let us 
pick a path $\pi_{\sst W,V}\in
[W,V]$. For an 
operator $\SO_V:\CH(V)\ra
\CH(V)$ we will define 
\begin{equation}\label{basechange}
\SA_{[W,V]}\big(\SO_V\big)\,=\,
\su(\pi_{\sst W,V})\cdot
\SO_V(\mu)\cdot\su^\dagger(\pi_{\sst W,V}).
\end{equation}
It is easy to convince oneself that $\SA_{W,V}\big(\SO_V(\mu)\big)$ 
does not depend on the choice of a path $\pi_{\sst W,V}\in
[W,V]$. We furthermore have
\begin{equation}\label{Atrsf}
\SA_{W,V}\big(\SR_V(\mu)\big)\,=\,\SR_{W}(\mu).
\end{equation}
It easily follows that
$\vartheta_V^{}(\mu_\2,\mu_\1)$ does not depend on $V$, i.e.
$\vartheta_V^{}(\mu_\2,\mu_\1)\equiv\vartheta(\mu_\2,\mu_\1)$.

To summarize:
The operators $\SR_{ V}(\mu)$ generate a projective
unitary representation $\SR_{ V}(G)$ of $G$ on $\CH$, 
\begin{align}
{}& \SR_{V}(\mu_\2)\cdot\SR_{V}(\mu_\1)
\,=\,\vartheta(\mu_\2,\mu_\1)\,\SR_{V}(\mu_\2\circ\mu_\1)\,.
\label{projgrouplaw2}
\end{align}
The operators $\SA_{[W,V]}$ express the unitary equivalence of the
representations $\SR_{ V}$ associated to the different $V\in\CG_\0$,
which allows us to regard 
\[
\SR\,\equiv\,
\Big[\,
\big(\SR_{V}^{}\big)_{V\in\CG_\0}^{}
\,,\,
\big(\SA_{E}^{}\big)_{E\in\CG_\1}^{}
\,\Big]
\]
as the representation of $G$ canonically associated to the given
representation of the groupoid ${\rm G}$.

\begin{rem}
There is of course some ambiguity in the construction, coming from
the choice of a representative $\pi_{\sst V}(\mu)\in[\mu. V,V]$.
However, it is clearly natural to consider two representations
$\sr$, $\sr'$ as equivalent if the generators $\sr_{\sst V}(\mu)$
and $\sr_{\sst V}'(\mu)$ differ from each other just by 
multiplication with a (possibly $\mu$-dependent) central element.
The cocycle $\vartheta$ of the representation $\sr$ will differ
from the cocycle $\vartheta'$ of $\sr'$ by a coboundary.
\end{rem}

\subsection{The projective representation of the 
Ptolemy groupoid on $\CK(\vf)$ \index{Ptolemy groupoid}}
\label{qPtolemy}

Following \cite{Ka3} closely we shall define a 
projective representation of the Ptolemy groupoid in 
terms of the following set of 
unitary operators on $\CK(\vf)$ 
\begin{equation}\begin{aligned}
\SA_v\;\equiv\;& e^{\frac{\pi i}{3}}
e^{-\pi i (\spp_v+\sq_v)^2}e^{-3\pi i \sq_v^2}\\
\ST_{vw}\;\equiv\;& e_b(\sq_v+\spp_w-\sq_w)e^{-2\pi i\spp_v\sq_w} ,
\end{aligned}\qquad
\text{where}\;\;v,w\in\vf_{\zero}\,.
\label{qPtgens}\end{equation}
The special function $e_b(U)$ can be defined in the strip 
$|\Im z|<|\Im c_b|$, $c_b\equiv i(b+b^{-1})/2$ by means of the 
integral representation
\begin{equation}
\log e_b(z)\;\equiv\;\frac{1}{4}
\int\limits_{i0-\infty}^{i0+\infty}\frac{dw}{w}
\frac{e^{-2{\mathsf i}zw}}{\sinh(bw)
\sinh(b^{-1}w)}.
\end{equation}
We refer to Appendix \ref{specapp} for more details on this 
remarkable special function. 
These operators are unitary for $(1-|b|)\Im b=0$. They satisfy the 
following relations \cite{Ka3}
\begin{align}
{\rm (i)}& 
\qquad\ST_{vw}\ST_{uw}\ST_{uv}\;=\;\ST_{uv}\ST_{vw},\label{pentrel}\\
{\rm (ii)}& \qquad\SA_{v}\ST_{uv}\SA_{u}\; =\; 
\SA_{u}\ST_{vu}\SA_{v},
\label{symrel}\\
{\rm (iii)}& \qquad\ST_{vu}\SA_{u}\ST_{uv}
\;=\;\zeta\SA_{u}\SA_{v}\SP_{uv},
\label{invrel}\\
{\rm (iv)}& \qquad\SA_u^3\;=\;\id,\label{cuberel}
\end{align}
where $\zeta=e^{\pi i c_b^2/3}$, $c_b\df\frac{i}{2}(b+b^{-1})$.
The relations 
\rf{pentrel} to \rf{cuberel} allow us to define a projective 
representation of the Ptolemy groupoid as follows.  
\begin{itemize}
\item Assume that $\om_{uv}\in[\vf',\vf]$.
To $\om_{uv}$ let us associate
the operator
\[ 
\su(\omega_{uv})\,\df\,\ST_{uv}\;:\;\CK(\vf)\ni \fv\;\,\ra\;\,
\ST_{uv}\fv\in\CK(\vf').
\]
\item For each fat graph $\vf$ and vertices $u,v\in \vf_{\zero}$
let us define the following operators
\[ \begin{aligned}
{}& \SA_{u}^{\vf}\;:\;\CK(\vf)\ni \fv\;\,\ra\;\,
\SA_{u}^{}\fv\in\CK(\rho_{u}\circ\vf).\\
{}& \SP_{uv}^{\vf}\;:\;\CK(\vf)\ni \fv\;\,\ra\;\,
\SP_{uv}^{}\fv\in\CK((uv)\circ\vf).
\end{aligned}
\]
\end{itemize}
It follows immediately from \rf{pentrel}-\rf{cuberel}
that the operators $\ST_{uv}$, $\SA_{u}$ 
and $\SP_{uv}$ can be used to generate
a unitary projective representation of the Ptolemy groupoid in 
$\CK(\vf)\simeq L^2(\BR^{2M})$.

\subsection{Reduction to the quantized Teichm\"uller spaces}\label{reduction}

\begin{thm}\label{thmred}
The isomorphism $\SI_\vf$ maps the operators $\su(\pi)$ which 
represent the Ptolemy groupoid on $\CK(\vf)$ to operators
of the form $\su'(\pi)=\SV_z(\pi)\ot\SV_h(\pi)$, where 
$\SV_z(\pi)\equiv(\SV_\ff(\pi))_{\ff\in B'(\Sigma)}$  
is a family of unitary operators $\SV_\ff(\pi)$ on 
$\CH_\ff(\vf)$.

For each fixed  $\ff\in B'(\Sigma)$ one may use the
operators $\SV_\ff(\pi)$ to generate
a unitary projective representation of
the Ptolemy groupoid.
\end{thm}

\begin{proof}
To begin with, let us note that
each path $\pi\in\rho\df[\vf',\vf]$ canonically defines a map 
$c\mapsto c'\equiv \FA(c)$
for each $c\in H_1(\Sigma,\BR)$. This map is defined in an obvious 
way for the elementary moves depicted in Figures \ref{ft}, \ref{flip} 
if we require that $c'$ coincides with $c$ outside the 
triangles depicted in these figures.
\begin{lem} $\quad$ \label{homact}\begin{itemize}
\item[a)] The map $c\mapsto \FA(c)$ preserves the symplectic (intersection) 
form on $H_1(\Sigma,\BR)$.
\item[b)] The operator $\su(\pi)$
maps 
$
\su(\pi)\cdot\sh_{\vf,c}\cdot(\su(\pi))^{-1}=\sh_{\vf',\FA(c)}.
$
\end{itemize} 
\end{lem}

\begin{proof}
Direct verifications.
\end{proof}

\begin{propn} \label{firstred} For each path $\pi\in[\vf', \vf]$ there 
exists an operator $\SH(\pi)$ on $\CK(\vf)$ such that 
the operators $\SV(\pi)$ on
$L^2(\BR^{2M})$ defined by 
\[ 
\SV(\pi) \,=\,\SH(\pi)\cdot\SU(\pi),\qquad
\SU(\pi)\,\equiv\,\SJ_{\vf'}^{-1}\cdot\su(\pi)\cdot\SJ_{\vf}^{}\,,
\]
\begin{itemize}
\item[(i)] commute with all operators ${\sh}_m^{}$, 
$\sh_m^{\sst\vee}$, $m=1,\dots,g$ and
$\sf_l$, $l=1,\dots,s-1$,
\item[(ii)] generate a unitary 
projective representation of ${\rm Pt}(\Sigma)$
on $L^2(\BR^{2M})$.
\end{itemize}
The operators $\SH(\pi)$ can  be represented
in the form 
\begin{equation}\label{SHpirepr}
\SH(\pi)\,\equiv\,
\exp(i H_{\pi}(\sh_{\vf'}))\,,
\end{equation}
where $H_{\pi}(\sh_{\vf'})$ is a quadratic function
of the $2g$ operators ${\sh}_1',\dots, {\sh}_g',
\sh_1'\!{}^{\sst^\vee}_{},\dots, 
\sh_g'\!{}^{\sst^\vee}_{}$ on $\CK(\vf')$.
\end{propn}

\begin{proof}
The existence of operators $\SH(\pi)$ of the 
form \rf{SHpirepr} which are such that statement (i) of the
Proposition is verified 
follows directly from Lemma \ref{canqlem} if one takes into
account that the transformation 
$\sh_{\vf',c}\ra\sh_{\vf',c'}$ is represented by an element 
of $Sp(g,\BR)$
according to part a) of Lemma \ref{homact}.

In order to prove statement (ii) of the proposition,
we mainly need to check that the operators $\SV(\pi)$ 
satisfy the relations of the Ptolemy groupoid. 
Let us consider a closed path $\pi\in[\vf,\vf]$ which decomposes
into a chain of edges as $\pi=\pi_n\circ\dots\circ\pi_{1}$.
\[
\SV(\pi)\,\equiv\,
\SV(\pi_n)\cdots \SV(\pi_1)
\]

On the one hand one may observe that
$\SV(\pi)$ can be factorized as 
\begin{equation} \label{SV} \SV(\pi)\,=\,\SH(\pi)\cdot 
\SU(\pi_n)\cdots \SU(\pi_1)\,=\,\SH(\pi)\zeta_\vf(\pi)\,.
\end{equation}
The operator $\SH(\pi)$ in \rf{SV} can be represented as follows:
\[
\SH(\pi)\,=\,\SH(\pi_n)\cdot
\big[\SU(\pi_{n}')\cdot \SH(\pi_{n-1})\cdot \SU(\pi_{n}')^{\dagger}\big]
\cdots 
\big[\SU(\pi_{1}')\cdot \SH(\pi_1)\cdot \SU(\pi_{1}')^{\dagger}\big]\,,
\]
where $\pi_{n-k}'\equiv\pi_n\circ\dots\circ\pi_{n-k}$. It follows from 
Lemma \ref{homact} together with \rf{SHpirepr} that
\[
\SU(\pi_{j+1}')\cdot 
\SH(\pi_{j})\cdot \SU(\pi_{j+1}')^{\dagger} 
\,=\,\exp( i H_{\pi_{j}}(\sh_{\vf}))\,.
\]
Taking into account equation \rf{ZHmap} we conclude that
\begin{equation}\label{SHfactor}
\SK^{-1}\cdot\SH(\pi)\cdot\SK\,=\, 1\ot \SH_h(\pi)\,.
\end{equation}

On the other hand let us note that $\SV(\pi)$ 
commutes with all operators ${\sh}_m^{}$, 
$\sh_m^{\sst\vee}$, $m=1,\dots,g$.
Equation \rf{SV}  implies that the same is true for
$\SH(\pi)$. However, the representation of the 
operators ${\sh}_m^{}$, 
$\sh_m^{\sst\vee}$, $m=1,\dots,g$ on $\CH_h(\vf)$ is irreducible
(see Remark \ref{irredrem}). This  allows us 
to conclude that $\SH(\pi)=\eta_{\vf}(\pi)\in \BC$, $|\eta_\vf(\pi)|=1$. 
Inserting this into \rf{SV} proves our claim.
\end{proof} 
It follows from statement (i) in Proposition \ref{firstred} that 
\begin{equation}
\SK^{-1}\cdot\SV(\pi)\cdot\SK\,=\, \SV_z(\pi)\ot 1\,.
\end{equation}
The task remains to describe the operators $\SV_z(\pi)$ more precisely.


\begin{propn}\label{secondred}
There exists a family of unitary operators 
$\SV_\ff(\pi):\CH_\ff(\vf)\ra \CH_\ff(\vf')$, $\ff\in B'(\Sigma)$ which
represent the action of the operator 
$\SV_z(\pi): \CH_z(\vf)\ra \CH_z(\vf')$.
\end{propn}
\begin{proof} To begin with, let us observe that it follows from 
\rf{qPtgens} that
the operators $\su(E)$ associated to the edges $E=\rho_v$, 
$E=\om_{vw}$ in ${\rm Pt}_\1(\Sigma)$
can all be factorized as
$\su(E)=\SQ(E)\cdot G(\hat{\sz}_e)$, where
$\SQ(E)$ is of the form $\SQ(E)=\exp(iQ_E(\sv))$
for a quadratic function $Q_E$,
and $G\equiv 1$ if $E=\rho_v$ and $G(z)=e_b(z)$ if $E=\om_{vw}$.
It follows that the corresponding operator
$\SV(E)$ defined in Proposition \ref{firstred} can be factorized as
\begin{equation}
\SV(E)\,=\,\SQ'(E)\cdot
G(\sz_{e}),\;\;{\rm where}\;\;\SQ'(E)\,\equiv\,
\SH(E)\cdot\SJ_{\vf'}^{-1}\cdot\SQ(E)\cdot\SJ_{\vf}^{}.
\end{equation}
The operator $\SQ'(E)$ is a product of four operators 
$\SJ_{\gamma_k}$, $\gamma_k\in Sp(2M,\BR)$ for $k=1,2,3,4$. If follows 
from Lemma \ref{canqlem} that it can be 
represented in the form 
\[
\SQ'(E)\,=\,\exp\big(i J_{E}(\sv_{})\big)\,,
\]
for some expression $J_{E}(\sv_{})$ which is quadratic
in $\sv_{}$. Note that the operators $G(\sz_{e})$ and $\SV(E)$ 
commute with all operators $\sh_m$, $\sh_m^{\sst\vee}$, $m=1,\dots,g$ and
$\sf_l$, $l=1,\dots,s-1$. It follows that the same is true for
$\SQ'(E)$, which implies that 
$
J_E(\sv_{})\equiv J_E(\sz)
$ 
depends only on the vector $\sz\equiv(\sz_e)_{e\in\vf_\1}$. 
$\SV(E)$ is therefore of the form
$
\SV(E)=\exp\big( i J_{E}(\sz)\big)\cdot
G(\sz_e)\,.
$
Our claim follows easily, 
$\SW_\ff(E)=\exp\big( i J_{E}(\sz_{\ff,e})\big)\cdot
G\big(\sz_{\ff,e})\big)$ does the job.
\end{proof}
Theorem \ref{thmred}
follows by combining Propositions  \ref{firstred} and \ref{secondred}.
\end{proof}

We are finally in the position to define more precisely
what we will regard as the quantized Teichm\"uller spaces. 
To this aim let us note that the Hilbert spaces 
$\CH_\0(\vf)$ associated to the origin $\0$ in $B'(\Sigma)$ 
form irreducible representations of the relations \rf{Qfockrel},
\rf{Qfockconstr}. Funktions of the operators
$\sz_{\0,e}$ generate the algebras $\BFB(\CH_\0(\vf))$ 
of all bounded operators
on $\CH_\0(\vf)$, which suggests to interprete 
$\BFB(\CH_\0(\vf))$ as 
particular representations 
of the quantized algebras 
$\CT_\hbar(\Sigma)$
of functions on the Teichm\"uller spaces.

The operators $\SV_\0(\pi)$ generate a unitary projective representation
of the Ptolemy groupoid which allows us to regard two operators
$\SO_{\vf_\2}\in\BFB(\CH_\0(\vf_\2))$ and 
$\SO_{\vf_\1}\in\BFB(\CH_\0(\vf_\1))$ as equivalent, 
$\SO_{\vf_\2}\sim\SO_{\vf_\1}$, iff
\begin{equation}
\SO_\2\,=\,\SV_\0(\pi)\cdot\SO_\1\cdot
(\SV_\0(\pi))^{-1},\quad
\pi\in[\vf_\2,\vf_\1]\,.
\end{equation}

Let furthermore $\SM\SC_\vf$ be 
the unitary projective 
representation of the mapping class group ${\rm MC}(\Sigma)$ 
which is generated from the operators $\SV_\0(\pi)$ 
 by means of the 
construction in Subsection \ref{symmG}.

\begin{defn} $\quad$ \begin{itemize}
\item[(i)] We define the algebra $\CT_\hbar(\Sigma)$
as the algebra generated by the
families $\CO\equiv(\SO_\vf)_{\vf\in\CP t_\0(\Sigma)}$ of bounded
operators $\SO_\vf$ on $\CH_\0(\vf)$ such that
$\SO_{\vf_\2}\sim\SO_{\vf_\1}$ for all 
$\vf_\2,\vf_\1\in\CP t_\0(\Sigma)$.
The algebra $\CT_\hbar(\Sigma)$ will
be called the quantized algebra of functions on the Teichm\"uller
spaces.

\item[(ii)] 
Let  $\CM\CC_\hbar(\Sigma)$ be the subalgebra
of $\CT_\hbar(\Sigma)$ generated by the families 
$\CM\CC(\mu)\equiv(\SM\SC_\vf(\mu))_{\vf\in \CP t_\0(\Sigma)}$
for all $\mu\in{\rm MC}(\Sigma)$.
\end{itemize}
\end{defn}

\newpage

\newpage
\part{A stable modular functor\index{modular functor} from the quantum 
Teichm\"uller spaces}

Let us recall that systems
of Fenchel-Nielsen coordinates\index{Fenchel-Nielsen coordinates} 
are naturally associated
to markings of the surfaces $\Sigma$.
The transformations between the different markings\index{marking} of a
Riemann surface $\Sigma$ generate yet another groupoid, called the stable
modular groupoid\index{modular groupoid}. Given that the quantum version of the
changes between the Penner coordinates
associated to different fat graphs was represented by a 
unitary projective
representation of the Ptolemy groupoid it is natural to expect that
the quantization of the Fenchel-Nielsen
coordinates should similarly come with a unitary projective
representation of the modular groupoid. 

Pants decompositions\index{pants decomposition} have one 
big advantage over ideal triangulations:
The gluing operation allows us to build ``larger'' surfaces from 
simple pieces of the same type, namely hyperbolic surfaces with 
holes. It is natural to require that the unitary projective
representations of the modular groupoid\index{modular groupoid} 
assigned to surfaces $\Sigma$
should correspondingly be
organized in a ``tower-like'' fashion: They should allow 
restriction to, and should be generated by the representations
assigned to the surfaces $\Sigma'$ which are obtained from
$\Sigma$ by cutting along simple closed curves on $\Sigma$.

Our aim in the rest of this paper will be to show that
such a structure can be constructed from the quantized
Teichm\"uller spaces considered in the previous part 
of this paper. This is of great importance
since 
\begin{equation*}\boxed{
\begin{aligned}
& \text{\bf having a tower of 
projective unitary representations
of the stable modular}\\
& \text{\bf  groupoid is equivalent to
having a stable unitary modular
functor.}
\end{aligned}}
\end{equation*}
The notion of a stable unitary modular functor
will be introduced in the next section.
One main difference to the more conventional 
(two-dimensional) modular functors 
as defined e.g. in \cite{Tu,BK2} is that
one restricts attention to Riemann surfaces
$\Sigma$ of genus $g$ and with $n$ parametrized boundary components
which are stable in the sense that $2g-2+n>0$.

We will then explain why having a stable unitary modular functor 
is equivalent 
to having a tower of
projective unitary representations
of the modular groupoid before we take up the task to
actually construct the latter from the quantization of the
Teichm\"uller spaces as described previously.

\section{The notion of a stable unitary modular functor}\label{stabMF}
\setcounter{equation}{0}

Given that the usual definitions 
of a modular functor take several pages to fully write them down
\cite{Tu,BK2},
we shall only briefly explain the most important
features. The missing details will not differ much 
from the definitions discussed in \cite{Tu,BK2}.


%
%

\subsection{Rigged Riemann surfaces\index{Rigged Riemann surfaces}}

We will consider compact oriented
surfaces $\Sigma$ with boundary 
$\pa\Sigma=\coprod_{\be\in A(\Sigma)} b_\beta$, where 
$A(\Sigma)\equiv\pi_0(\pa\Sigma)$ is the set of connected components of
$\pa\Sigma$.
A surface $\Sigma$ is called an extended surface if one has
chosen
orientation-preserving homeomorphisms 
$p_\be:b_\be\ra S^1$
for each connected component $b_\be$ of the
boundary. To be concrete, let $S^1=\{z\in\BC;|z|=1\}$. 
An e-surface $\Sigma$ of genus $g$ and with $n$ boundary circles 
is called stable if $2g-2+n>0$.

We will use the terminology
rigged Riemann surface, or r-surface for short, for triples 
\newcommand{\RSigma}{\widehat{\Sigma}}
$\RSigma=(\Sigma,y,\fc)$,
where 
\begin{itemize}\item[$\triangleright$]
$\Sigma$ is a stable extended surface,
\item[$\triangleright$] $y$ is a Lagrangian subspace of 
$H_1(\Sigma,\BR)$, and
\item[$\triangleright$] $\fc:A(\Sigma)\ra\CL$ is  
a coloring of the boundary of $\Sigma$ by
elements of a set $\CL$. 
\end{itemize}

Given an r-surface $\RSigma$ 
and given $\be,\be'\in A(\Sigma)$ such that $\fc(\be)=\fc(\be')$
we can define a new r-surface $\RSigma'\equiv\sqcup_{\be\be'}\RSigma=(\Sigma',y',\fc')$ which is 
obtained from $\RSigma$ by gluing the boundary components
$b_\be$ and $b_{\be'}$. \begin{itemize}\item[$\triangleright$]
The surface $\Sigma'\equiv\sqcup_{\be\be'}\Sigma$
is defined by identifying all points $p\in b_\be$ with
$(p_{\be'}^{-1}\circ p_\be^o)(p)\in b_{\be'}$, where $p_\be^o:b_\be\ra S^1$
is defined by $p_\be^o(p)=-\overline{p_\be(p)}$. There is a 
corresponding projection $P_{\be\be'}:\Sigma\ra\Sigma'\equiv\sqcup_{\be\be'}\Sigma$ 
which maps $b_\be$, $b_{\be'}$ to the same simple closed curve on $\Sigma'$.
\item[$\triangleright$]
The Lagrangian subspace $y'$ is given by the image of
$H_1(\Sigma,\BR)$ under the projection $P_{\be\be'}$, 
\item[$\triangleright$] The
coloring $\fc'$ is obtained from $\fc$ by putting $\fc'(\al')=\fc(\al)$ if 
$P_{\be\be'}(b_\al)=b_{\al'}\in \pa\Sigma'$.
\end{itemize}



\subsection{Stable unitary modular functors}\label{stablesubsec}

Let $\CL$ now be a space with a measure $d\nu$. A stable modular functor 
with central charge is the following collection of data.
\begin{itemize}
\item[$\triangleright$] Assignment $\RSigma=(\Sigma,y,\fc)\ra 
\CH(\Sigma,y,\fc)$, where
\begin{itemize}
\item $\RSigma$ is an r-surface,
\item $\CH(\Sigma,y,\fc)$ is a separable Hilbert space,
\end{itemize}
\item[$\triangleright$] {\it Mapping class group:} Assignment
\[
[f]\;\longrightarrow\;\;\big(\,
\SU_{[f]}: \CH(\Sigma,y,\fc)\ra \CH(\Sigma',y'_f,\fc'_f)\,\big),
\]
\begin{itemize}
\item $[f]$  is the isotopy class of a homeomorphism 
$f:\Sigma\ra\Sigma'$,
\item $\SU_{[f]}$ is a unitary operator,
\item $y'_f$ 
is the Lagrangian subspace of $H_1(\Sigma',\BR)$ determined from $y$ via $f$,
\item $\fc'_f:A(\Sigma')\ra\CL$: the coloring of boundary components
of $\Sigma'$ induced from $\fc$ via $f$.
\end{itemize}
\item[$\triangleright$] {\it Disjoint union:} There exist 
unitary operators
\[ 
\SG_{\2\1}:\CH(\Sigma_\2\sqcup\Sigma_\1,y_\2\oplus y_\1,\fc_\2\sqcup\fc_\1)\,
\overset{\sim}{\longrightarrow}\,
\CH(\Sigma_\2,y_\2,\fc_\2)\otimes \CH(\Sigma_\1,y_\1,\fc_\1).
\]
\item[$\triangleright$] {\it Gluing:} Let $(\Sigma',y')$ 
be obtained from $(\Sigma,y)$ by gluing of 
two boundary components $\al,\be$. There then exists a unitary operator
\[ 
\SG_{\al\be}\;:\;\int_{\CL}d\nu(s)\;
\CH\big(\Sigma,y,\fc_{s\lfloor\al\be}\big)\,
\overset{\sim}{\longrightarrow}\,
\CH(\Sigma',y',\fc'),\]
where the 
coloring 
$\fc_{s\lfloor\al\be}:A(\Sigma)\ra\CL$ is defined from $\fc'$ via
\[\begin{aligned}
& \fc_{s\lfloor \al\be}(c)=\fc'(P_{\al\be}(c))\quad{\rm if} \quad\ga\in 
A(\Sigma)\setminus\{\al,\be\},\\ 
& \fc_{s\lfloor \al\be}(c)=s\quad{\rm if} \quad\ga\in \{\al,\be\}.
\end{aligned}
\]
\end{itemize}
These data are required to satisfy the following ``obvious'' consistency
and compatibility conditions:
\begin{itemize}
\item[$$] {\bf Multiplicativity:} For all 
homeomorphisms $f:\Sigma_1\ra\Sigma_2$,
$g :\Sigma_2\ra\Sigma_3$ there exists $\zeta(f,g)\in S^1$ such thay
we have
\begin{equation}\label{multcond}
\SU_{[f]}\SU_{[g]}\,=\,\zeta(f,g)\SU_{[f\circ g]}.
\end{equation}
$\zeta(f,g)$ has to satisfy the condition
$\zeta(f,g)\zeta(f\circ g,h)=\zeta(g,h)\zeta(f,g\circ h)$.
\item[$$] {\bf Functoriality:} The gluing isomorphisms and the
disjoint union isomorphisms are functorial in $\Sigma$. \footnote{One 
is considering
the category with objects r-surfaces, and morphisms isotopy classes 
of homeomorphisms of r-surfaces, equipped additionally with the
gluing and disjoint union operations.}

\item[$$] {\bf Compatibility:} The gluing isomorphisms and the
disjoint union isomorphisms are mutually compatible.
\item[$$] {\bf Symmetry of Gluing:} $G_{\al\be}=G_{\be\al}$.
\end{itemize}

It would take us several pages to write out all conditions
in full detail, we therefore refer to \cite{Tu,BK2} for more details.
However, it seems that the following two ``naturality'' requirements
represent a key to the understanding of the notion of the
modular functor:

\begin{quote}
{\bf Naturality:}\\[1ex]
{\bf a)} Let $f_\1:\Sigma_\1^{}\ra\Sigma_\1'$, 
$f_\2:\Sigma_\2^{}\ra\Sigma_\2'$ be r-homeomorphisms.
We then have
\[
\SG_{\2\1}^{}\cdot \SU_{[f_\2\sqcup f_\1]}\cdot
\SG_{\2\1}^{\dagger}
\,\equiv\,\SU_{[f_\2]}\ot\SU_{[f_\1]}\,.
\]
{\bf b)} If the r-homeomorphism $f:\Sigma_\1\ra\Sigma_\2$ induces 
an r-homeomorphism $f':\Sigma_\1'\ra\Sigma_\2'$ of the surfaces 
$\Sigma_1'$, $\Sigma_2'$ obtained from $\Sigma_1$, $\Sigma_2$
by the gluing construction, we have
\[
\SG_{\al\be}^{}\cdot \left(
\int_{\CL}d\nu(s)\;\SU_{[f]}\big(\Sigma,y,\fc_{s\lfloor\al\be}\big)
\right)\cdot \SG_{\al\be}^{\dagger}
\,=\,
\SU_{[f']}(\Sigma',y',\fc'\big).
\]
\end{quote}
These requirements make clear how  the mapping class group 
representations on the spaces $\CH(\Sigma',y',\fc')$
restrict to and are generated by the representations
assigned to the surfaces $\Sigma$ which are obtained from
$\Sigma'$ by cutting along simple closed curves on $\Sigma'$.

\begin{rem} The standard definitions of modular functors
assume that the Hilbert  spaces
$\CH(\Sigma,y,\fc)$ are {\it finite-dimensional}.
They are therefore not suitable for 
{\it nonrational} conformal field theories.
Our definition should be seen as a first step towards
the definition of  analogs of modular functors which are
associated to nonrational conformal field 
theories\index{rational conformal field theory}\index{conformal field theory}
in a way 
similar to the connections between rational conformal field theories
and modular functors mentioned in the introduction.

However, there is one important ingredient of the usual definition that
does not have an obvious counterpart in our framework. In the more
standard definitions of modular functors it is required that there
exists a distinguished
element $s_\0$ in $\CL$ which has the property that coloring a boundary 
component $\be$ with $s_\0$ is equivalent to ``closing'' this boundary component. 
More precisely, let $\fc_{s\lfloor\be}$ be a coloring of the boundary components
of an extended surface $\Sigma$ which assigns $s\in\CL$ to the 
component with label $\be$,
and let $\RSigma_{\check{\beta}}=
(\Sigma_{\check{\beta}},y_{\check{\beta}},
\fc_{\check{\beta}})$ be the r-surface 
obtained from $\RSigma=(\Sigma,y,\fc)$ by gluing a disc
to $b_{\beta}$. The more
standard definitions of modular functors assume or imply existence of an
element $s_\0$ in $\CL$ such that 
\[
\CH(\Sigma,y,\fc_{s_\0\lfloor\be})\simeq \CH(\Sigma_{\check{\beta}},y_{\check{\beta}},
\fc_{\check{\beta}}),\quad
\SU_{[f]}(\Sigma,y,\fc_{s_\0\lfloor\be})\simeq \SU_{[f]}
(\Sigma_{\check{\beta}},y_{\check{\beta}},
\fc_{\check{\beta}}).
\]
This yields additional relations between the mapping class group representations
assigned to surfaces with different numbers of boundary components. 

In the case of the quantized Teichm\"uller spaces it ultimately turns
out that an analog of the element $s_\0$ in $\CL$ exists only if one
considers the analytic continuation of $\CH(\Sigma,y,\fc)$
and $\SU_{[f]}(\Sigma,y,\fc)$ with respect to the boundary labels $\fc(\be)\in\CL$.
This fact, and the corresponding improvement of our definition of a stable
modular functor will be elaborated upon elsewhere.
\end{rem}

\subsection{Representations
of the modular groupoid
versus modular functors} \label{towerdef}

It turns out to be very useful to reformulate the notion
of a stable unitary modular functor in  terms of generators and relations as
follows. Any surface $\Sigma'$ can be glued from a 
surface $\Sigma_\0=\coprod_{p\in\si_\0}T_p$ which is a disjoint union of
trinions. The 
different ways of doing this can be parametrized\footnote{Note 
that markings with the same 
cut system will yield equivalent representations for $\Sigma'$.
This redundancy 
will give a useful book-keeping device when the representation of 
the mapping class group is considered.}
by markings $\si$. 
The gluing construction determines a canonical 
Lagrangian subspace $y'_\si$ of $H_1(\Sigma',\BR)$
from the tautological Lagrangian subspace
$y_\0\equiv H_1(\Sigma_\0,\BR)$.
By 
iterating the gluing and disjoint union isomorphisms one 
defines unitary operators $\SG({\si},\fc'):
\CH(\Sigma',y'_\si,\fc')\ra\CH(\si,\fc')$, where  
\begin{equation}\label{CHtriniondecomp}
\CH(\si,\fc')\,\equiv\,
\int_{\mathfrak L}d\nu_\si^{}({\rm S})\;\bigotimes_{p\in\si_\0}
\CH(T_p^{},\fc_{\si,p}^{{\rm\sst S}})\,.
\end{equation}
We have used the following notation:
\begin{itemize}
\item The integration is extended over the set ${\mathfrak L}$ 
of all colorings ${\rm S}:\CC_\si\ni c\ra s_c\in\CL$ equipped with the
canonical product measure $d\nu_\si^{}({\rm S})$ obtained from $d\nu$ 
by choosing any
numbering of the elements of $\CC_\si$.
\item $\fc_{\si,p}^{{\rm\sst S}}$ 
is the coloring of the 
boundary components of $T_p$ which is defined by
assigning
\[ \fc_{\si,p}^{{\rm\sst S}}(\be)\,=\,\left\{\begin{aligned}
& \fc'(\be') \;\;{\rm if}\;\;P_{\si,p}^{}(\be)=\be'\in A(\Sigma'),\\
& s_c \;\;{\rm if}\;\; P_{\si,p}^{}(\be)=c\in\CC_\si.
\end{aligned}\right.
\]
$P_{\si,p}^{}(\be)$ is the embedding $T_p\hookrightarrow\Sigma'$ defined from $\si$
by the gluing construction.
\end{itemize}

Unitary mappings between the different spaces $\CH(\si,\fc)$ 
arise in two
ways: First, one may have different markings $\si_\2$, $\si_\1$ such that
the Lagrangian subspaces defined
by the gluing construction coincide, $y_{\si_\1}=y_{\si_\2}$.
In this case one gets unitary operators
$\SF_{\si_\2\si_\1}(\fc):\CH(\si_\1,\fc)\ra\CH(\si_\2,\fc)$ from the 
composition 
\begin{equation}
\SF_{\si_\2\si_\1}(\fc)\,\equiv\,\SG({\si_\2},\fc)\cdot 
(\SG({\si_\1},\fc))^{\dagger}.
\end{equation}
Secondly, one has the mappings $\SU_{[f]}$ which may map
between spaces $\CH(\Sigma,y_\si,\fc)$
and $\CH\big(\Sigma',(y_{\si})'_f,\fc'_f\big)$. In the case that
$\mu:\Sigma\ra\Sigma$ 
represents an element of the mapping class group it is natural
to define operators $\SV_{\si\mu}(\fc):\CH(\si,\fc)
\ra\CH(\si'_\mu,\fc'_\mu)$ by
\begin{equation}
\SV_{\si\mu}(\fc)\,=\,\SG(\si'_\mu,\fc'_\mu)\cdot \SU_{[\mu]}\cdot
(\SG({\si},\fc))^{\dagger},
\end{equation}
where 
$\si'_\mu$ is the image of the marking $\si$ under $\mu$.

It turns out - as will be reviewed in the next section - that there
exists a set $\CM_\1(\Sigma)$ of elementary moves between 
the different markings $\si$ such that
any two markings $\si_\2,\si_\1$ can be connected by 
paths $\pi$ which are composed out of the elementary moves. 
There furthermore exists a set $\CM_\2(\Sigma)$ of relations
which makes the resulting two-dimensional CW  complex $\CM(\Sigma)$
simply connected.
The correponding path groupoid is called the modular groupoid and denoted
by ${\rm M}(\Sigma)$.
 
Our construction of the operators $\SF_{\si_\2\si_\1}(\fc)$, 
$\SV_{\si\mu}(\fc)$ assigns unitary operators to each of the 
elementary moves . This yields unitary operators 
$\SU(\pi,\fc):
\CH(\si_\1,\fc)\ra
\CH(\si_\2,\fc_\pi)$ for each path $\pi\in[\si_\2,\si_\1]$, where
$\fc_\pi$ is the coloring of boundary components which is obtained from $\fc$
by tracking the relabeling of boundary components defined by $\pi$.
Due to the fact that
multiplicativity holds only projectively, cf. eqn. \rf{multcond}, 
one will find that there exists $\zeta(\pi)\in S^1$ such 
that $\SU(\pi)=\zeta(\pi)$ for all {\it closed} paths, starting and
ending at the same marking $\si$. In other words, a stable unitary modular
functor canonically defines projective unitary representations
of the stable modular groupoid.

Conversely, the definition of 
a {\it tower of projective unitary representations of the modular
groupoids} ${\rm M}(\Sigma)$ involves the following data.
\begin{itemize}
\item[$\triangleright$] Assignment $(\si,\fc)\ra \CH(\si,\fc)$,
where $\CH(\si,\fc)$ is constructed out of spaces $\CH(S_\3,\fc_\3)$
assigned to trinions as in \rf{CHtriniondecomp}.
\item[$\triangleright$]  Assignment 
\[
\pi\in[\si_\2,\si_\1]\longrightarrow \big(\,\SU(\pi,\fc)\,:\,
\CH(\si_\1,\fc)\ra
\CH(\si_\2,\fc_\pi)\,\big).
\]
\end{itemize}
It is required that the operators $\SU(\pi,\fc)$ generate unitary projective
representations of the modular groupoids associated to the surfaces $\Sigma$.
In order to formulate the additional requirements which turn this collection
of representations of ${\rm M}(\Sigma)$ into a {\it tower} let us 
note that markings can be glued in a 
natural way. We will use 
the notation $\sqcup_{\al\be}\si$ for the 
marking obtained by gluing boundary components
$\al$ and $\beta$. Gluing and disjoint union 
now become realized in a trivial manner, 
\begin{align}
& \CH(\sigma_\2\sqcup\sigma_\1,\fc_\2\sqcup\fc_\1)\,=\,
\CH(\sigma_\2,\fc_\2)\otimes \CH(\sigma_\1,\fc_\1),\\
& 
\CH(\sqcup_{\al\be}\sigma,\fc')=\int_{\CL}d\nu(s)\;
\CH\big(\sigma,\fc_{s\lfloor\al\be}\big)\,.
\label{factor2}\end{align}
The data specified above are then required to satisfy the following
naturality conditions. 
\begin{quote}
{\bf Naturality:}\\[1ex]
{\bf a)} Let $\pi=\pi_\2\sqcup\pi_\1\in[\si_\2',\si_\2]
\sqcup[\si_\1',\si_\1]$
We then have
\begin{equation}
\SU(\pi_\2\sqcup\pi_\1,\fc_\2\sqcup\fc_\1)
\,\equiv\,\SU(\pi_\2,\fc_\2)\ot\SU(\pi_\1,\fc_\1)\,.
\end{equation}\label{natur1}
{\bf b)} Given a path $\pi\in[\si_\2,\si_\1]$ in 
$\CM(\Sigma)$ 
let $\sqcup_{\al\be}\pi\in [\sqcup_{\al\be}\si_\2,
\sqcup_{\al\be}\si_\1]$ be the corresponding path 
in $\CM(\sqcup_{\al\be}\Sigma)$ defined by the 
gluing construction. We then have
\begin{equation}\label{natur2}
\SU(\sqcup_{\al\be}\pi,\fc') \,=\,
\int_{\CL}d\nu(s)\;\SU(\pi,\fc_{s\lfloor\al\be}\big).
\end{equation}
\end{quote}

From a tower of projective unitary representations of the modular
groupoids one can reconstruct a stable 
unitary modular functor as follows. The system of isomorphisms
$\SU(\pi,\fc):\CH(\si_\1,\fc)\ra \CH(\si_\2,\fc)$ allows us to 
identify the Hilbert spaces associated to different markings $\si$, 
thereby defining $\CH(\Sigma,\fc)$.
The representation of the modular groupoid defines a representation
of the mapping class group on $\CH(\Sigma,\fc)$ via the 
construction in subsection \ref{projgroup}. We refer to \cite{BK2}
for more details.

It is important to note (see next section for a detailed discussion) 
that the elements of the set $\CM_\1(\Sigma)$ of 
elementary moves only change the
markings within subsurfaces of genus zero with three or four holes,
or within subsurfaces of genus one with one hole.
The corresponding operators will be called Moore-Seiberg data.
They characterize a tower of 
representations of the modular groupoid completely.
Let us furthermore note that 
the faces/relations that one needs to 
define a two-dimensional CW complex $\CM(\Sigma)$ which has
as set of edges $\CM_\1(\Sigma)$ turn out to involve
only subsurfaces of genus zero with three to five holes,
and subsurfaces of genus one with one or two holes. 
This means that one only needs to verify a finite number
of relations to show that a given set 
of Moore-Seiberg data defines a stable unitary modular functor.

\section{The modular groupoid\index{modular groupoid}}\label{modgroup}
\setcounter{equation}{0}

The modular groupoid is the groupoid generated from the 
natural transformations relating the different markings of a
Riemann surface  $\Sigma$
\cite{MS}\cite{BK}. 
A complete set of generators 
and relations has been determined in \cite{MS,BK,FuG}. 

\begin{rem}
The CW complex $\CM(\Sigma)$ defined below will be a subcomplex
of the complex denoted $\CM^{\rm max}(\Sigma)$ 
in \cite{BK}, since that reference
allowed for cut systems that yield connected components 
with less than three boundary components, whereas we
will exclusively consider
cut systems that yield connected components 
with exactly three holes.
All other deviations from \cite{BK} are due to slightly different 
conventions in the definition of the generators.
\end{rem}

\subsection{Notations and conventions}


We have a unique curve $c(\si,e)\in\CC'$ associated to each 
edge $e\in\si_\1$ of the marking graph $\Gamma_\si$.
The trinions $T_p\in\CP_\si$ are in one-to-one correspondence 
with the vertices $p\in\si_\0$ of $\Gamma_{\si}$. 

\begin{figure}[htb]
\epsfxsize8cm
\centerline{\epsfbox{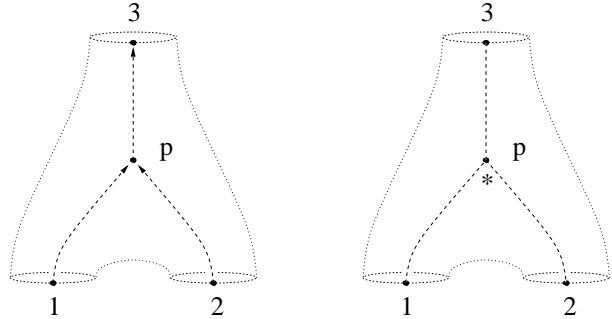}}
\caption{Two representations for the decoration on a marking graph}\label{marking}
\end{figure}

The choice of a distinguished boundary component $c_p$ for each trinion 
$T_p$,  $p\in\si_\0$ will be called the
{\it decoration} of the marking graph $\Gamma_{\si}$. 
The distinguished boundary component $c_p$ will be called
{\it outgoing}, the other two boundary components of $T_p$ {\it incoming}. 
Two useful
graphical representation for the 
decoration are depicted in Figure \ref{marking}. 

\subsection{Generators}\label{MS-gens}

\begin{figure}[htb]
\epsfxsize6cm
\centerline{\epsfbox{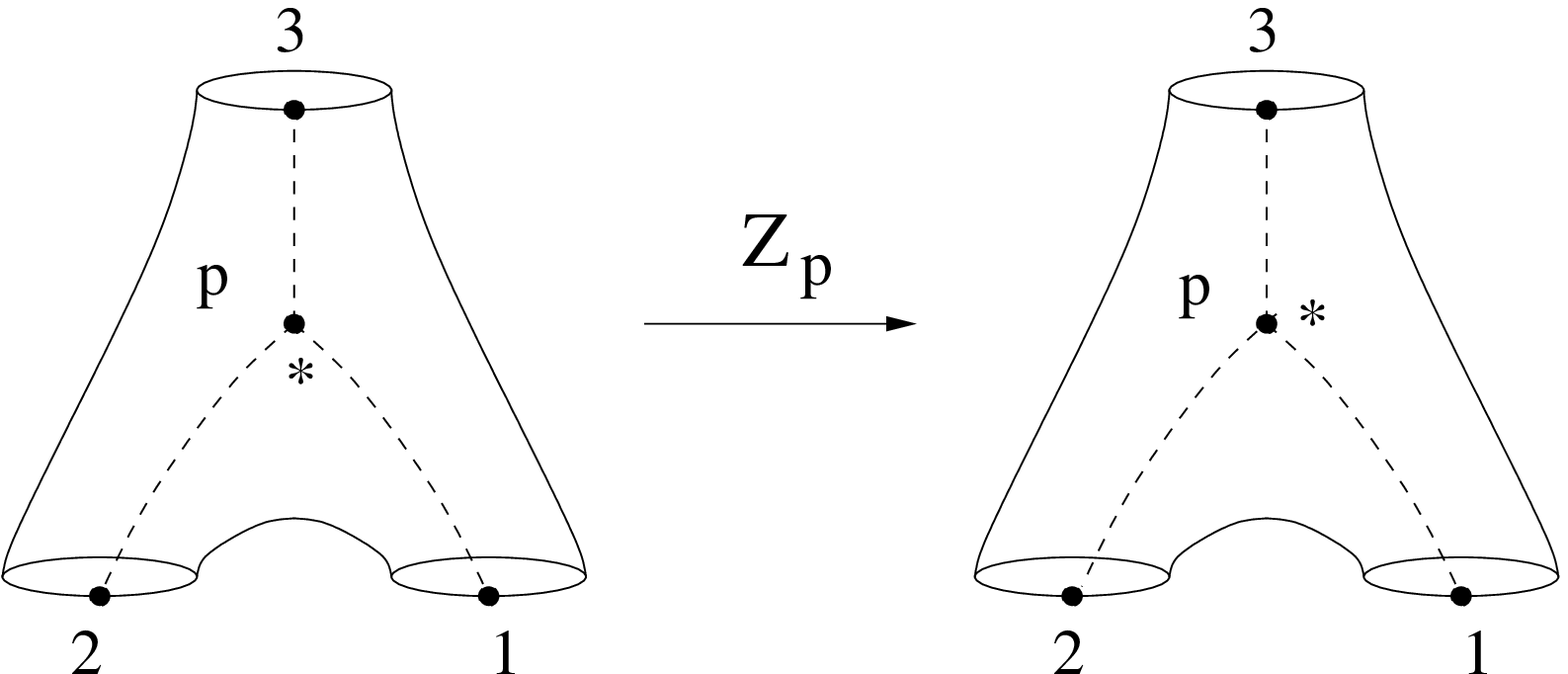}}
\caption{The Z-move}\label{zmove}\vspace{.3cm}
\epsfxsize6cm
\centerline{\epsfbox{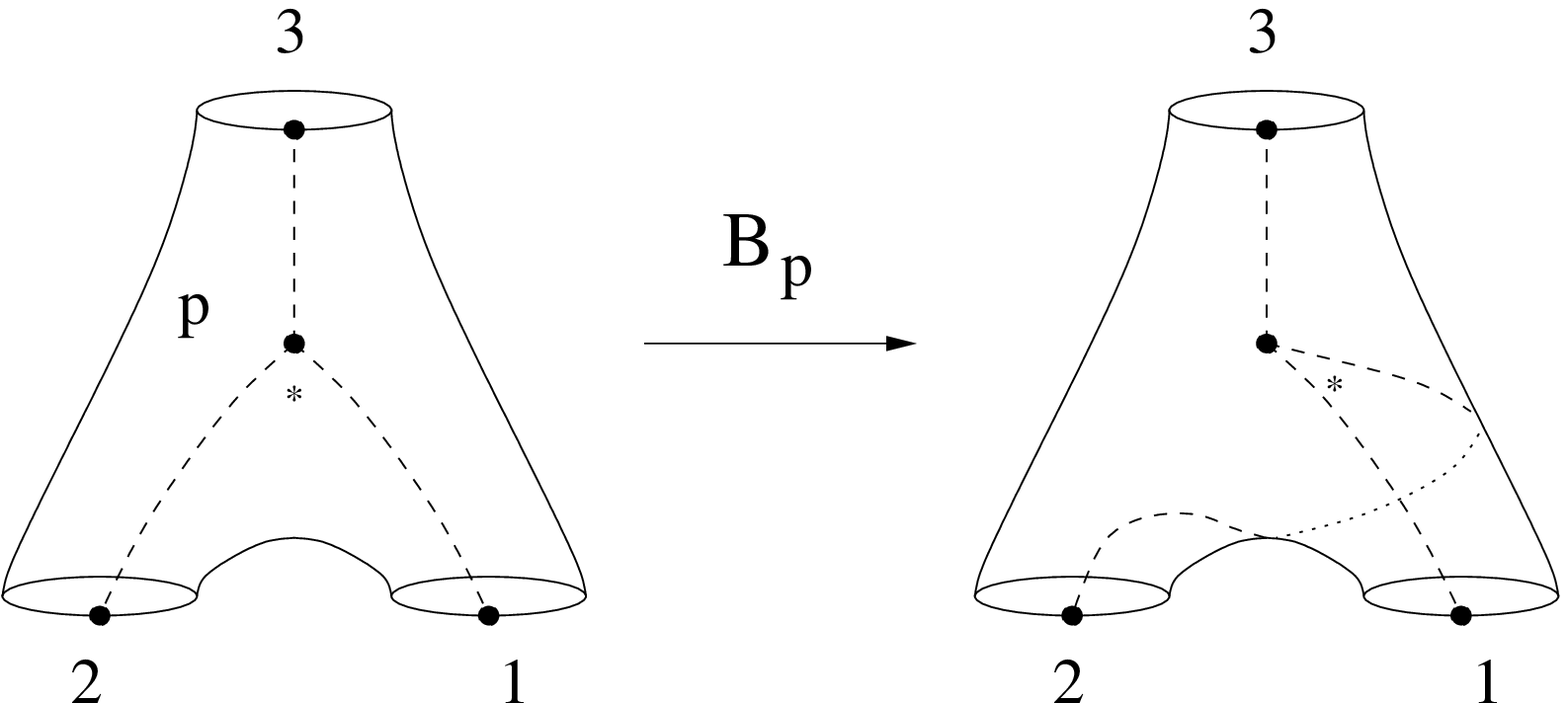}}
\caption{The B-move}\label{bmove}
\end{figure}

\begin{figure}[htb]
\epsfxsize7cm
\centerline{\epsfbox{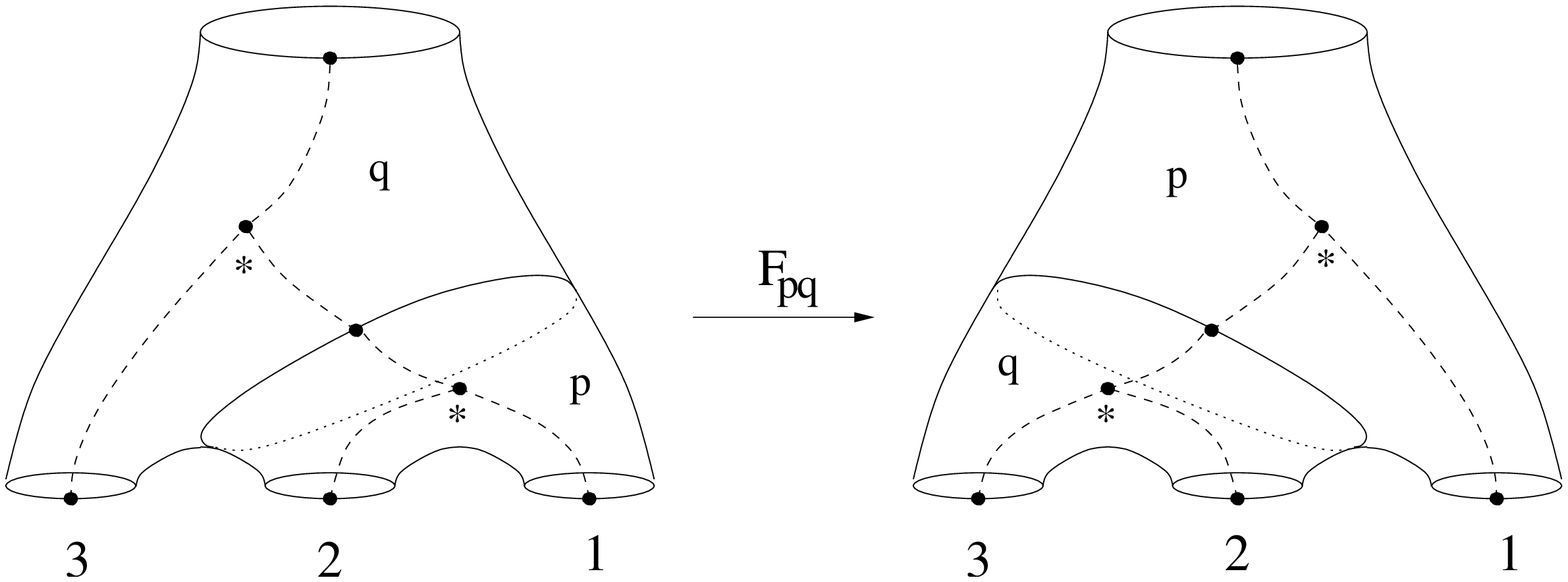}}
\caption{The F-move}\label{fmove}\vspace{.3cm}
\epsfxsize8cm
\centerline{\epsfbox{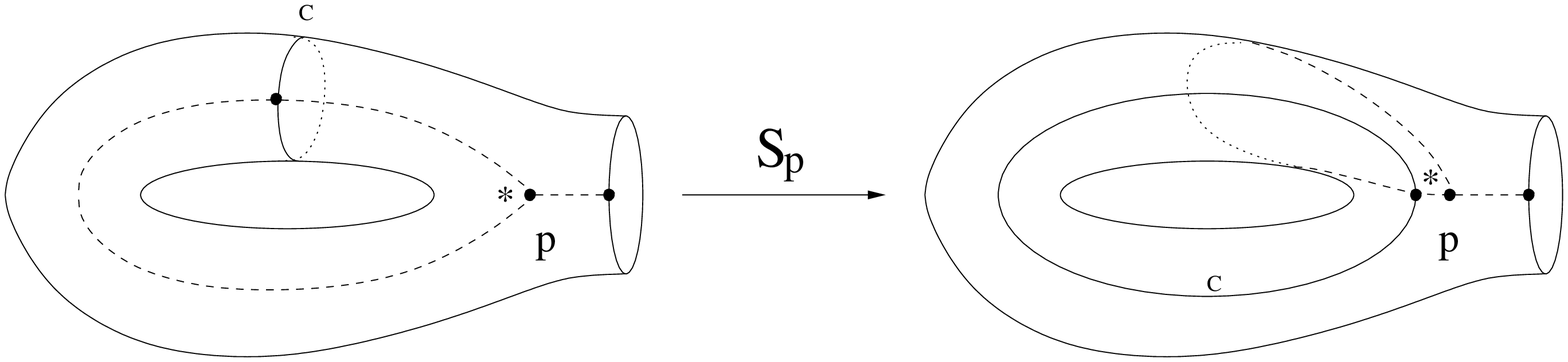}}
\caption{The S-move}\vspace{.3cm}
\label{smove}\end{figure}

The set of edges ${\CM}_{\1}(\Sigma)$ will be given by elementary
moves denoted as $(pq)$,
$Z_p$, $B_p$, $F_{pq}$ and $S_p$. The indices $p,q\in\si_{\0}$ 
specify the relevant trinions within the 
pants decomposition of $\Sigma$ that is determined by $\sigma$. 
The move $(pq)$ will simply be the operation in which the
labels $p$ and $q$ get exchanged.
Graphical representations for the elementary 
moves $Z_p$, $B_p$, $F_{pq}$ and $S_p$ are given in 
Figures \ref{zmove}-\ref{smove}.

\subsection{Relations}\label{MSrels}

The relations of ${\rm M}(\Sigma)$ correspond to 
the faces $\varpi\in\CM_\2(\Sigma)$. 
In the following we will define a set $\CR$ of faces 
which is large enough to make the complex $\CM(\Sigma)$
simply connected. A face $\varpi$
in $\CM_\2(\Sigma)$ may be characterized 
by choosing a chain
$\Gamma_\varpi=
E_{\varpi,n(\varpi)}\circ\dots\circ E_{\varpi,1}$, where
$E_{\varpi,i}\in \CM_{\1}(\Sigma)$ for $i=1,\dots,n(\varpi)$.
In order to simplify notation we
will generically factorize $\Gamma_\varpi$
as $\Gamma_\varpi=\Gamma_\varpi^{\2}\circ\Gamma_\varpi^{\1}$
and write $(\Gamma_\varpi^{\2})^{-1}=\Gamma_\varpi^{\1}$
instead of $\Gamma_\varpi=\id$.

\paragraph{Locality.}

Let us introduce the 
notation ${\rm supp}(m)$ by ${\rm supp}(m)=\{p\}$ if $m=Z_p,B_p,S_p$,
${\rm supp}(m)=\{p,q\}$ if $m=(pq),F_{pq}$ and
${\rm supp}(m_\2\circ m_\1)={\rm supp}(m_\2)\cup{\rm supp}(m_\1)$. 
We then have
\begin{equation}\label{locality}
m_\2\circ m_\1\,=\,m_\1\circ m_\2\quad {\rm whenever}\;\;
{\rm supp}(m_\1)\cap {\rm supp}(m_\2)=\emptyset\,.
\end{equation}
We will list the remaining relations ordered by the topological type of the
surfaces on which the relevant graphs can be drawn.

\paragraph{Relations supported on surfaces of genus zero.}

\begin{figure}[htb]
\epsfxsize11cm
\centerline{\epsfbox{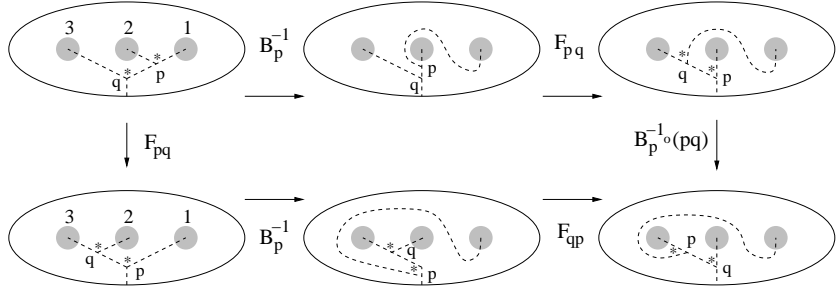}}
\caption{The hexagon relation.}\label{hexafig}
\end{figure}

\begin{figure}[t]
\epsfxsize10cm
\centerline{\epsfbox{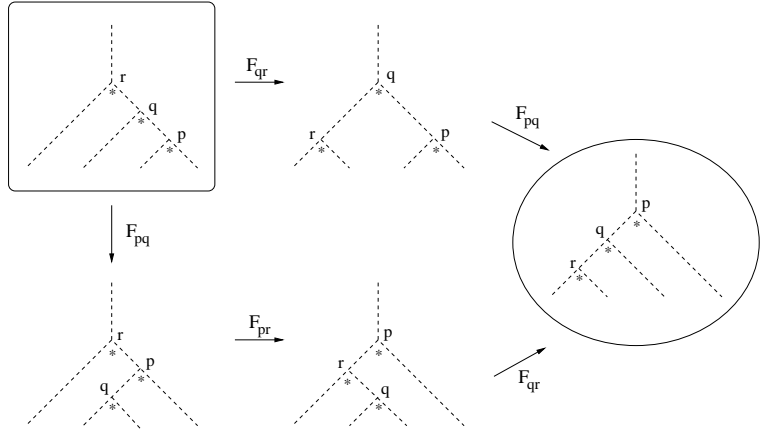}}
\caption{The pentagon identity}
\label{pentafig}
\end{figure}

\newcommand{\LR}{\quad\Leftrightarrow\quad}
\begin{align}\label{zrel}
g=0,~~s=3:\quad & Z_p\circ Z_p\circ Z_p = \id \, .\\[1ex]
g=0,~~s=4:\quad & 
\begin{aligned}\label{hexarel} {}{\rm a)} \quad& 
F_{qp}^{\phantom{1}}\circ B_p\circ F_{pq}^{\phantom{1}}\;=\;
(pq)\circ B_q \circ F_{pq}^{\phantom{1}}\circ B_p\,,\\
{\rm b)}\quad & F_{qp}^{\phantom{1}}
\circ B_p^{-1}\circ F_{pq}^{\phantom{1}}\;=\;
(pq)\circ B_q^{-1} \circ F_{pq}^{\phantom{1}}\circ B_p^{-1},\\
{\rm c)} \quad & A_{pq}\circ A_{qp}\;=\;(pq)\,.
\end{aligned}\\[1ex]
g=0,~~s=5:\quad & F_{qr}\circ F_{pr}\circ F_{pq}\;=\; F_{pq}\circ F_{qr}.
\label{pentarel}\end{align}
We have used the abbreviation
\begin{equation}\label{Adef}
A_{pq}\;\df\; Z_q^{-1}\circ
F_{pq}^{\phantom{1}}\circ Z_q^{-1}\circ Z_p^{\phantom{1}}\,.
\end{equation}
In Figures \ref{hexafig} and \ref{pentafig} we have given
diagrammatic representations for relations \rf{hexarel}, b) and
\rf{pentarel} respectively.
\paragraph{Relations supported on surfaces of genus one.}

\begin{figure}[htb]\epsfxsize10cm
\centerline{\epsfbox{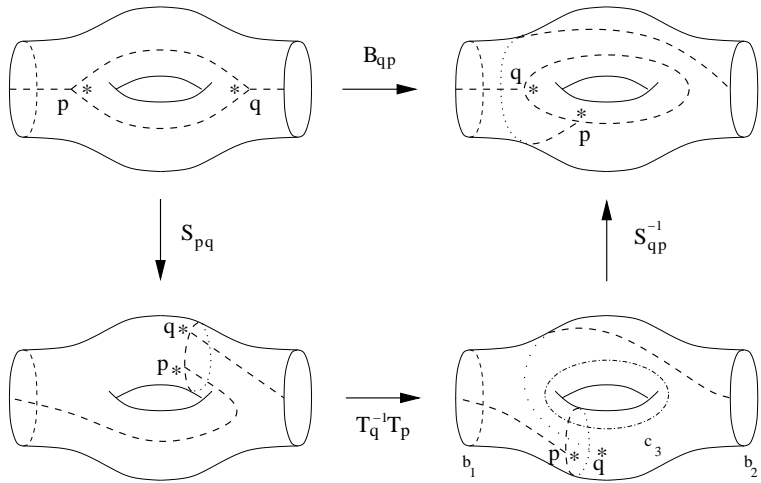}}
\caption{Relation for the two-punctured torus.}
\label{twotor}\end{figure}

In order to write the relations transparently let us introduce the
following composites of the elementary moves. 
\begin{align}
\label{b-composites}g=0,~~s=3:\qquad & \begin{aligned}
{\rm a)}\quad & 
B_p'\;\df\; Z_p^{-1}\circ B_p\circ Z_p^{-1},\\
{\rm b)}\quad & T_{p} \;\df\; Z_p^{-{1}}\circ
B_p\circ Z_p^{\phantom{1}}\circ B_p\, ,
\end{aligned}\\[1ex] \label{bcompdef}
g=0,~~s=4:\qquad &
B_{qp}\,\df\, Z_q^{-1}\circ F_{qp}^{-1}\circ B_q'
\circ F_{pq}^{-1}\circ Z_q^{-1}\circ(pq)
\,,\\[1ex]
g=1,~~s=2:\qquad &
S_{qp}\;\df\; (F_{qp}\circ Z_q)^{-1}\circ S_p\circ (F_{qp}\circ Z_q)\,.
\label{scompdef}\end{align}
It is useful to observe that the move
$T_{p}$,  represents the Dehn twist around the 
boundary component of the trinion $\ft_p$ numbered by $i=\1$ 
in Figure \ref{marking}.
With the help of these definitions we may write the 
relations supported on surfaces of genus one as follows:
\begin{align}\label{onetor:a}g=1,~~s=1:\quad & 
\begin{aligned}{\rm a)}\quad &
S^2_p \;=\; B'_p,\\
{\rm b)}\quad & S_p^{\phantom{1}}\circ  T_p^{\phantom{1}}
\circ S_p^{\phantom{1}}\;=\;T^{-1}_p\circ S_p^{\phantom{1}}\circ T^{-1}_p .
\end{aligned}\\[1ex]
g=1,~~s=2:\quad & 
B_{qp}^{\phantom{1}}\;=\;S_{qp}^{-1}\circ T_{q}^{-1}
T_{p}^{\phantom{1}}
\circ S_{pq}^{\phantom{1}}\,.
\label{twotorrel}\end{align}
Relation \rf{twotorrel} is represented diagrammatically in 
Figure \ref{twotor}.

\begin{thm}\label{BKthm}  
The complex $\CM(\Sigma)$ is connected and simply connected for 
any e-surface $\Sigma$.
\end{thm}

\begin{proof}
The theorem follows easily from \cite{BK}, Theorem 5.1. 
We noted previously that our complex $\CM(\Sigma)$ differs
from the complex 
$\CM^{\rm\sst max}(\Sigma)$ of \cite{BK} in having
a set of vertices which corresponds to decompositions
of $\Sigma$ into connected components with {\it exactly} three holes.
The edges of $\CM^{\rm\sst max}_\1(\Sigma)$ which 
correspond to the F-move of \cite{BK} simply can't appear in 
$\CM_\1(\Sigma)$. Otherwise the set of edges of $\CM(\Sigma)$
coincides with the relevant subset of 
$\CM^{\rm\sst max}_\1(\Sigma)$, with the exception
that our move $B_p$ is a composition of the $B$- and the $Z$-move
of \cite{BK}. To complete the proof it remains to check that
our set of faces is equivalent to the subset of 
$\CM^{\rm\sst max}_\2(\Sigma)$ which involves only
the vertices $\CM_\0(\Sigma)$ of
our smaller complex $\CM(\Sigma)$. This is a useful exercise.
\end{proof}

\begin{defn}\label{CMprimedef} Let ${\CM}'(\Sigma)$ 
be the complex which has the same set of vertices as ${\CM}(\Sigma)$, a set of 
edges given by the moves
\begin{equation}\label{moveslist}
(pq)~,~Z_p~,~B_p~,~B'_p~,~S_p~,~T_{p}~,~
F_{pq}~,~A_{pq}~,~B_{pq}~,~S_{pq}
\end{equation}
defined above, as well as faces given by equations 
\rf{locality}-\rf{twotorrel}. 
\end{defn}




\section{From markings\index{marking} to fat 
graphs\index{fat graphs}}\label{pt->MS}
\setcounter{equation}{0}

A key step in our construction of a stable modular functor 
from the quantized Teichm\"uller spaces will be the definition 
of a distinguished class of fat graphs $\vf_\si$ which
are associated to the elements $\si$ of a 
certain subset 
of the set $\CM_\0(\Sigma)$ of all markings of $\Sigma$. 
\begin{defn}
Let $\CA_\si\subset\CC_\si$ be the set 
of all curves $c$ which are incoming for both adjacent trinions.
We will say that a marking $\si$ is
admissible iff there is no curve $c\in\CC_\si$ which is outgoing for
both adjacent trinions, and if
cutting $\Sigma$ along all curves $c\in\CA_\si$
yields connected components all of which have genus zero.
The set of all admissible markings will be denoted by
$\CM_\0^{\rm ad}(\Sigma)$.
\end{defn}

To each admissible marking $\si$ 
we may naturally associate a fat graph $\vf_{\si}$ on $\Sigma$
by the following construction. 
In order to construct the fat graph $\vf_\si$, 
it will be useful to consider a certain
refinement of the pants decomposition associated to 
$\si$ which is defined as follows. For each  curve
$c\in\CA_{\si}$ let $A_c$ be a small annular neighborhood
of $c$ which contains $c$ in its interior, and which is bounded
by a pair of curves $(c^+,c^-)$ that are isotopic to $c$.
\begin{figure}[t]
\epsfxsize8cm
\centerline{\epsfbox{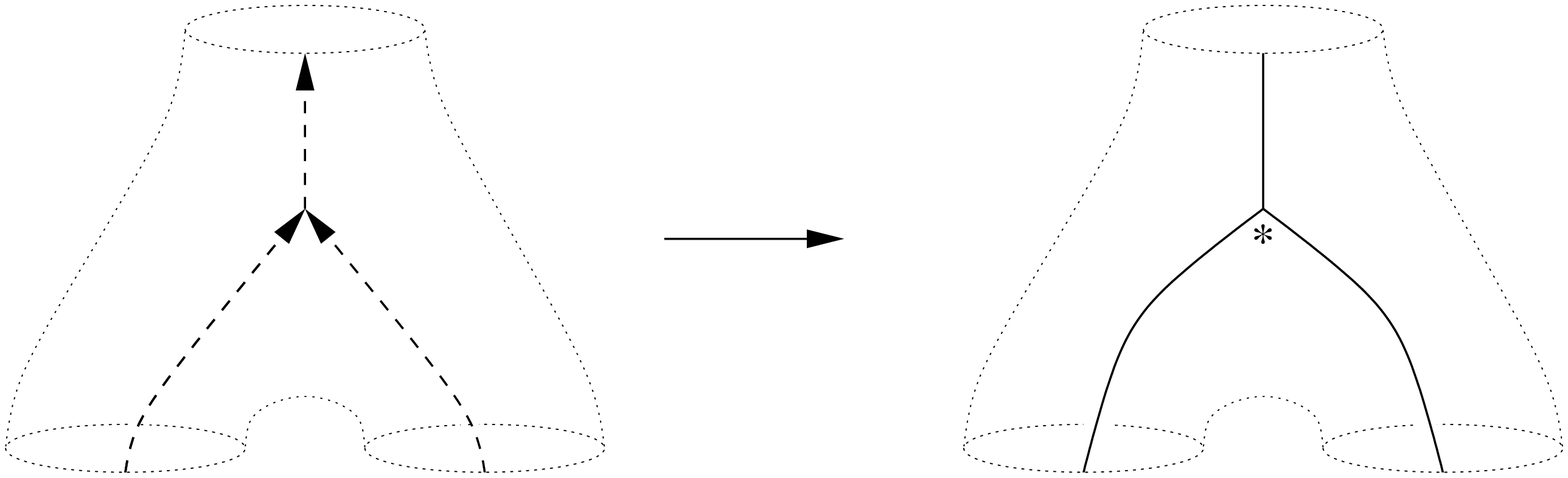}}
\caption{Substitution for a trinion $T$}\label{subst1}
\epsfxsize8cm
\centerline{\epsfbox{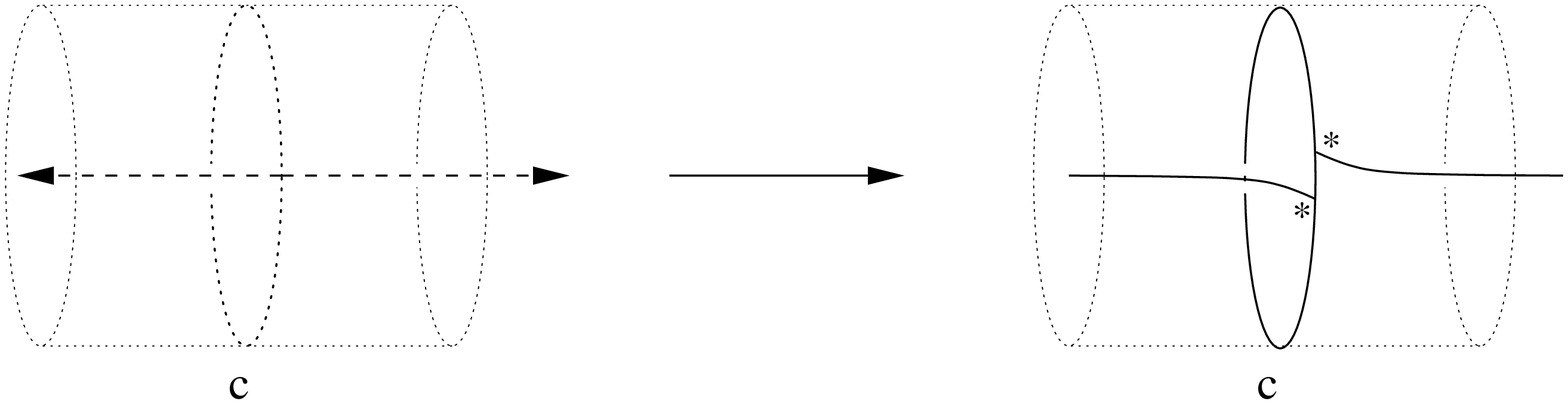}}
\caption{Substitution for an annulus $A$}\label{subst2}
\epsfxsize6cm
\centerline{\epsfbox{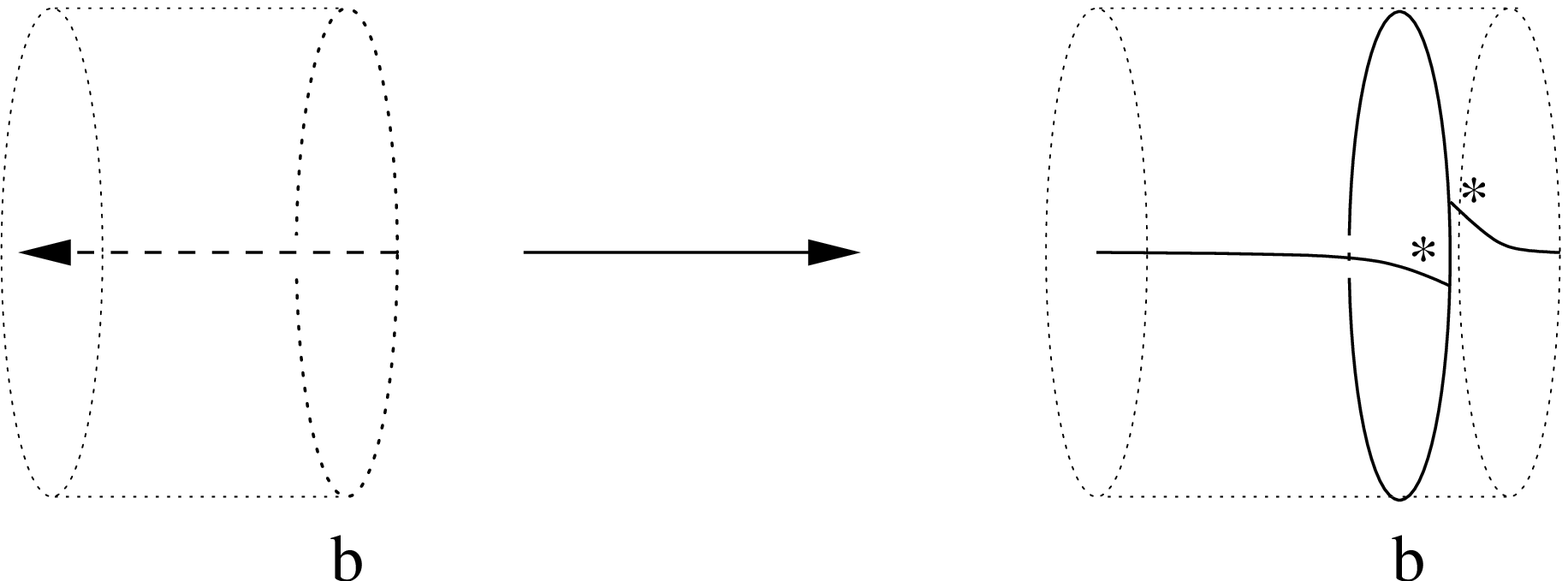}}
\caption{Substitution for a boundary component $B$}\label{subst3}
\end{figure}
If $c\in A(\Sigma)$ represents a boundary component of $\Sigma$
we may similarly consider an annular 
subset $B_c$ of $\Sigma$ bounded by $c^+\equiv c$ 
and another curve $c^-$ isotopic to $c$. If $c\in\CC_\si$ 
belongs to neither
of these two classes we will simply set $c^+\equiv c\equiv c^-$.
By cutting along all such curves $c^\pm$ 
we obtain a decomposition of $\Sigma$ into 
trinions $T_p$, $p\in \si_\0$ and annuli $A_c$, $c\in\CA_{\si}$, 
each equipped with a marking graph. We may then replace the markings 
on each of these connected components 
by fat graphs according to Figures \ref{subst1}-\ref{subst3}. 
The re-gluing of trinions $T_p$, $p\in \si_\0$ and annuli 
$A_c$, $B_c$ to recover the surface $\Sigma$ 
may then be performed in such a way that the fat graphs 
on the connected components glue to a fat graph $\vf_\si$
on $\Sigma$.

\subsection{The complex $\CM^{\rm ad}(\Sigma)$}

It is natural to consider the complex $\CM^{\rm ad}(\Sigma)$
for which the set of vertices $\CM^{\rm ad}_\0(\Sigma)$
is the subset of $\CM_\0(\Sigma)$ which consists of the 
admissible markings, and which has a set of edges 
$\CM^{\rm ad}_\1(\Sigma)$ given by the subset of $\CM_\1'(\Sigma)$
that contains those edges which connect two 
admissible markings.

\begin{propn}\label{Madconn}
The complex $\CM^{\rm ad}(\Sigma)$ is connected and simply connected.
\end{propn}

\begin{proof}
Let us consider two
admissible markings $\si,\si'\in\CM_{\0}(\Sigma)$. 
There exists a path $\varpi$ in $\CM(\Sigma)$ which connects 
$\si$ and $\si'$. This path may be represented
as a chain $C_\varpi=E_{\varpi,n(\varpi)}\circ\dots\circ E_{\varpi,1}$, 
$E_{\varpi,i}\in\CM_i(\Sigma)$ 
composed out of the moves $Z_p$, $F_{pq}$, $B_p$ and $S_p$.

For a given marking $\si\in\CM_{\0}(\Sigma)$ 
let $[\si]$ be the set of all markings 
$\si'$ which differ from $\si$ only in the choice of decoration.
The moves $Z_p$ act transitively on $[\si]$. By inserting 
$Z_p$ moves if necessary we may therefore modify 
$C_\varpi$ to a chain $D_\varpi$ which takes the
form
\[
D_\varpi\,=\,Z_{\varpi,n(\varpi)}\circ
F_{\varpi,n(\varpi)}\circ Z_{\varpi,n(\varpi)-1}\circ
\dots\circ F_{\varpi,1}\circ Z_{\varpi,0},
\]
where $Z_{\varpi,i}$, $i=0,\dots,n(\varpi)$ are chains composed out
of $Z_p$-moves only, and the moves
$F_{\varpi,n(\varpi)}\in\CM_{\1}(\Sigma)$ connect markings
$\tilde{\si}_{\varpi,i}$ and $\si_{\varpi,i}$, 
$i=1,\dots,n(\varpi)$ which are
{\it admissible}. We clearly must have $[\si_{\varpi,i+1}]=
[\tilde{\si}_{\varpi,i}]$, $i=1,\dots,n(\varpi)-1$ and
$[\si_{\varpi,1}]=[\si]$, $[\si_{\varpi,n(\varpi)}]=[\si']$.
Connectedness of $\CM^{\rm ad}(\Sigma)$
would follow if the chains $Z_{\varpi,i}$, $i=0,\dots,n(\varpi)$
are homotopic to chains $Y_{\varpi,i}$
which represent paths in $\CM^{\rm ad}(\Sigma)$.
That this indeed the case follows from the following lemma.

\begin{lem}
Assume that $\si_a,\si_b\in\CM_\0^{\rm ad}(\Sigma)$ satisfy
$\si_b\in[\si_a]$. There then exists a path $\varpi_{ab}$ in 
$\CM^{\rm ad}(\Sigma)$ which connects $\si_a$ and $\si_b$.
\end{lem}

\noindent{\it Sketch of proof}. $\frac{\quad}{}$ 
Let us call  a marking 
$\si$ {\it irreducible}
if it is admissible and if there are no edges $e$ 
in $\CA_{\si}$ such that cutting $\Sigma$ along $c(e,\si)$ yields
two disconnected components. 
A marking $\si$ which is  irreducible has 
only one outgoing external edge.
With the help of Lemma \ref{rotlem}  in Appendix D
it is easy to show that for an admissible 
graph $\si$ there always exists a sequence
of $Z_p$-moves in $\CM^{\rm ad}_\1(\Sigma)$
which transforms $\si$ to a graph $\si'$ that is irreducible.

We may and will therefore assume that $\si_a$ and $\si_b$
are both irreducible. 
Using Lemma \ref{rotlem} again allows us to transform
$\si_b$ to a marking $\si_c$ which is such
that the outgoing external edges of $\si_a$ and $\si_c$
correspond to the same outgoing boundary component of $\Sigma$.

The graphs $\si_a$ and $\si_c$ can finally be connected by
a chain which is composed out of moves $B_p$ and $F_{qp}$ {it only}.
This follows from the connectedness of $\CM(\Sigma_\0)$ in the
case of a surface $\Sigma_\0$ of genus zero \cite{MS}.
In this way
one constructs a sequence of moves that connects $\si_a$ to $\si_b$
\hfill $\square$

It remains to prove that $\CM^{\rm ad}(\Sigma)$ is simply connected
as well. Let us consider any closed path $\varpi$ in $\CM^{\rm ad}(\Sigma)$.
The path $\varpi$ is contractible in $\CM'(\Sigma)$. The deformation
of $\varpi$ to a trivial path may be performed recursively,
face by face. The crucial observation to be made is the following one.

\begin{lem}\label{contractlem}
The paths $\varpi$ which represent the boundaries
of the faces of $\CM'_\2(\Sigma)$ are paths
in $\CM^{\rm ad}(\Sigma)$. 
\end{lem}

\begin{proof}
By direct inspection of the relations \rf{zrel}-\rf{twotorrel}.
\end{proof}

The Lemma \ref{contractlem} implies that deforming a path 
$\varpi\in\CM^{\rm ad}(\Sigma)$ by contracting a face in $\CM'_\2(\Sigma)$
will produce a path $\varpi'$ which still represents a path 
in $\CM^{\rm ad}(\Sigma)$. It follows that $\varpi$ is contractible
in $\CM^{\rm ad}(\Sigma)$.
\end{proof}

\subsection{Separated variables}\label{sepofvar}

We had observed in 
Section \ref{coordholes} that the Fock variables $z_e$, $e\in\vf_\1'$ form
a set of coordinates for the Teichm\"uller spaces of surfaces with holes.
When considering fat graphs $\vf_\si$ associated to a marking
$\si$ it will be useful to replace the 
Fock variables $z_e$, $e\in\vf_\1'$ by an alternative 
set of coordinates $({q}_c,{p}_c)$, $c\in\CC_\si$
for $T_\vf$ which satisfy
\begin{equation}\label{CCRcl}
\begin{aligned}
& \Om_{\vf}({p}_{c_\1}, {p}_{c_\2})\,=0=\,
\Om_{\vf}({q}_{c_\1}, {q}_{c_\2}),\\
&
\Om_{\vf}({p}_{c_\1}, {q}_{c_\2})\,=\,\de_{c_\1,c_\2},
\end{aligned}
\qquad c_\1,c_\2\in\CC_\si.
\end{equation}
These coordinates are constructed as follows.

\begin{figure}[htb]
\epsfxsize3cm
\centerline{\epsfbox{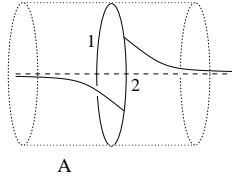}}
\caption{Annulus $A_c$ and fat graph $\vf$ on $A$.}
\label{cycle}\end{figure}

For $c\in\CA_\si$ let $A_c$ be an annular neighborhood of
$c$ such that the part of $\vf_\si$ which is contained in 
$A_c$ is isotopic to the model depicted in Figure \ref{cycle}. 
Let $e_\1$ and $e_\2$ be the
two edges that are entirely contained in $A_c$ with labelling 
defined by Figure \ref{cycle}.
\noindent Out of ${z}_{e_\1}$ and ${z}_{e_\2}$ we may then define
\begin{equation}\label{pqc1}
{q}_c\,\df\,\fr{1}{2}({z}_{e_\1}-{z}_{e_\2}),\qquad
{p}_c\,\df\,-\fr{1}{2}({z}_{e_\1}+{z}_{e_\2})\, .
\end{equation}

\begin{figure}[htb]\epsfxsize3.5cm
\centerline{\epsfbox{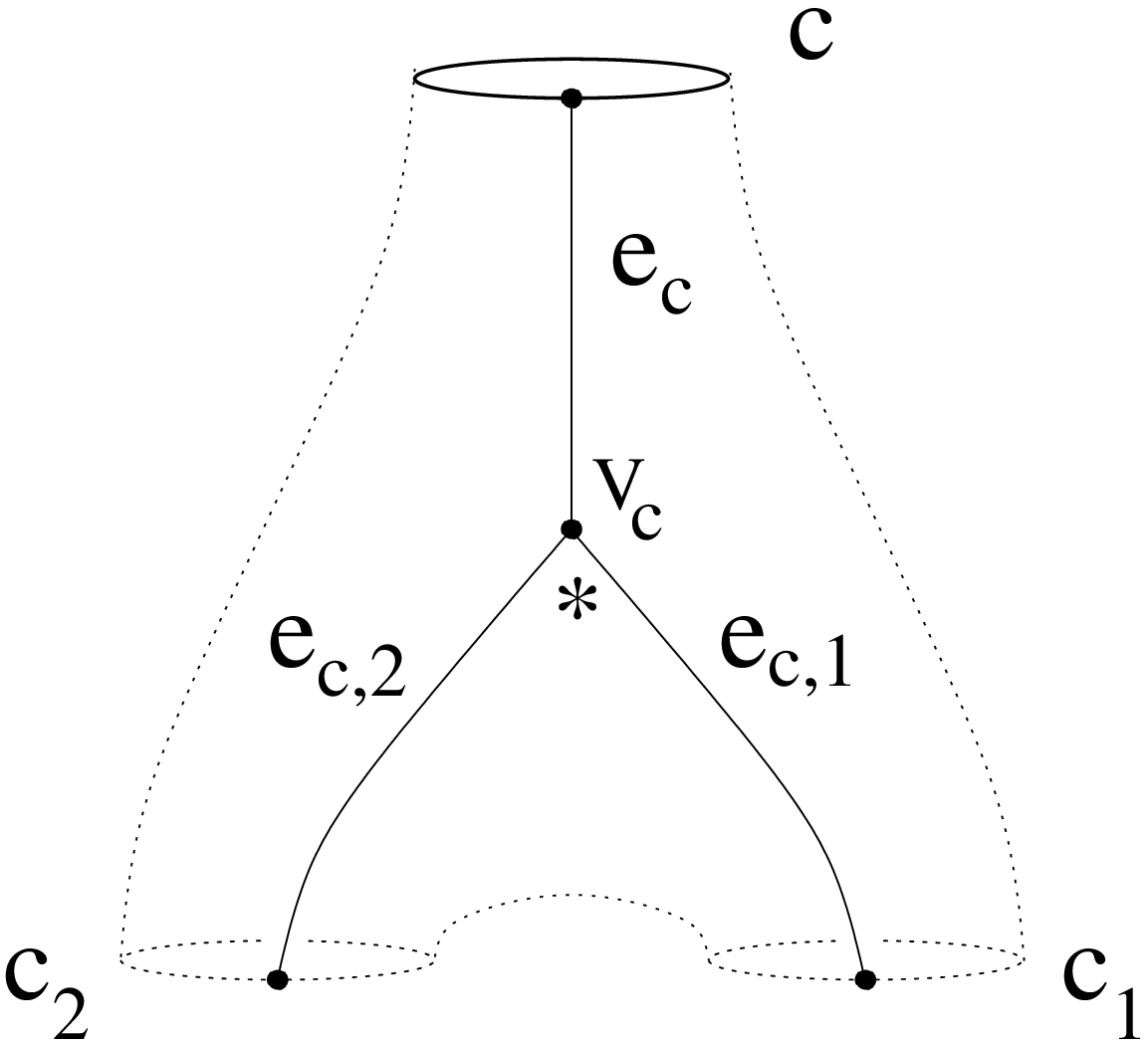}}
\caption{Part of $\vf_\si$ contained in $T_c$}
\label{threefat}\end{figure}

In the case $c\in\CC_\si\setminus\CA_\si$ let us note that 
there is a unique trinion $T_c$ for which the curve $c$ is the outgoing 
boundary component. The part of $\vf_\si$ contained in $T_c$ is 
depicted in Figure \ref{threefat}.

Let $e_{c,\2}$ and $e_{c,\1}$ be the
edges of $\vf_\si$  as 
indicated in Figure \ref{threefat},
and let $c_\ep$, $\ep=1,2$ be the curves
which represent the corresponding boundary components
of $T_p$. We shall then define
\begin{equation}\label{pqc2}
{q}_c\,\df\,y_{{c_\2}},\quad
{p}_c\,\df\,-y_{{c_\1}},
\end{equation}
where
\begin{equation}\label{pqc3}
y_{c_\ep}\,\df\,{z}_{e_{c,\ep}}
+\frac{1}{2}{f}_{c_\ep},\quad\ep=1,2.
\end{equation}

The same construction yields linear combinations 
$(\hat{p}_c,\hat{q}_c)$,
$c\in\CC_\si$ of the Kashaev variables which represent the
coordinates $({p}_c,{q}_c)$ within $W_{\vf_\si}$.
Note that
our construction of the fat graph $\vf_\si$ implies that $T_c$ contains a
unique vertex $v_c\in\vf_{\si,\0}$. 
It turns out that
the variables $\hat{q}_c$, $\hat{p}_c$ have a simple relation
to the Kashaev variables $q_{v_c}$, $p_{v_c}$. 
\begin{lem}
\[
\begin{aligned}
{\rm (i)} \quad& \hat{q}_c=q_{v_c}+h_{c_\2},\qquad 
\hat{p}_c=p_{v_c}-h_{c_\1}.\\
{\rm (ii)} \quad &
\Om_{\vf}(p_{c},q_{c'})\;=\;\de_{c,c'},\qquad
c,c'\in\CC_{\si}.
\end{aligned}
\]
\end{lem}
\begin{proof} 
In order to prove part (i) let us first
look at the vertices $v\in\vf_{\si,\0}$
which appear on the graph geodesic homotopic to $c$. The
contribution to $y_{c_\ep}$ of the edges that are incident to $v$
is $\hat{z}_{e_\1'}+\hat{z}_{e_\2'}+\hat{z}_{e_\3'}$ if 
$e_i'$, $i=1,2,3$ are the 
three edges incident to $v$. 
It then follows from \rf{kafock} that $y_c$ does not depend on both
$p_{v}$ and $q_{v}$. What remains are contributions from the vertex
$v_c$, as well as contributions from the vertices contained in annuli 
$A_{c'}$ which are determined as follows. 
Note that the homology class
$[c]$ can be decomposed as a linear combination of
classes $[c']$ with $c'\in\CA_{\si}$. It is straightforward to
check that each $c'\in\CA_{\si}$ which appears in this decomposition
yields a contribution $-h_{c'}$ to  $\hat{f}_c$.
These contributions sum up to give $-h_{c}$. What remains
is the contribution from the vertex $v_c$. Part (i) of the lemma 
now follows easily by recalling the definition \ref{kafock}.

Part (ii) of the lemma is a trivial consequence of part (i).
\end{proof}

\subsection{Quantized Teichm\"uller spaces for surfaces
with holes}

The coordinates introduced in  the previous subsection
make it straightforward to modify the discussion of the quantization of 
Teichm\"uller spaces from the case of punctured Riemann surfaces 
(Section \ref{Quantsimple}) to the
case of Riemann surfaces with holes.

Bearing in mind that the Fock variables $z_e$, $e\in\vf_\1'$  
are unconstrained in the
present case (see Lemma \ref{TTlem}) leads us to identify the
algebra of functions on the Teichm\"uller spaces with
the algebras of function of the variables
$z_e$, $e\in\vf_\1'$. 
A convenient set of coordinates is 
given by the coordinate functions $(p_c,q_c)$, $c\in\SC_\si$. 
Canonical quantization of the 
Teichm\"uller spaces is therefore straightforward, and leads to 
an algebra of operators with generators $(\spp_c,\sq_c)$, $c\in\CC_\si$,
which is 
irreducibly represented  (in the sense of Remark \ref{irredrem})
on the Hilbert space 
$\CH_z(\si)\simeq L^2(\BR^{3g-3+2s})$.
We will normalize the operators $(\spp_c,\sq_c)$, $c\in\CC_\si$,
such that 
\begin{equation}\label{CCR}
\begin{aligned}
{\rm (i)} \quad &
[\,{\spp}_{c}\, ,\, {\sq}_{c'}\,]\;=\;(2\pi i)^{-1}
\,\de_{cc'},
\quad c,c'\in\CC_\si,\\
{\rm (ii)} \quad &
[\,{\spp}_{c}\, ,\, {\spp}_{c'}\,]\;=\;0\;=\;
[\,{\sq}_{c}\, ,\, {\sq}_{c'}\,].
\end{aligned}
\end{equation}

Quantization of the Kashaev space $W_{\vf_\si}$ 
produces a {\it reducible} 
representation of the algebra \rf{CCR} which is generated by
operators  $(\hat{\sq}_c,\hat{\spp}_c)$ associated to the pairs
of variables $(\hat q_c,\hat p_c)$. It also yields operators
$\sh_c$, $c\in H_1(\Sigma,\BR)$ which represent 
the quantization of the Poisson vector space $H_{\vf_\si}$
with basis $h_c$, $c\in H_1(\Sigma,\BR)$. The algebra generated by the 
$\sh_c$, $c\in H_1(\Sigma,\BR)$
has a center generated by the $\sh_c$, $c\in B(\Sigma)$.
Following the 
discussion at the beginning of Subsection \ref{Qredsubsec} one 
constructs a representation of this algebra on the space
\[
\CH_h(\si)\,\equiv\,\int_{B'(\Sigma)}d\ff \;\CH_{h,\ff}(\si)\,,
\]
where $\CH_{h,\ff}(\si)\simeq L^2(\BR^g)$ 
is an irreducible representation of the 
algebra $i[h_{c_\2},h_{c_\1}]=b^2{\rm I}(c_\2,c_\1)$ for 
$c_\2,c_\1\in H_1(\Sigma_{\rm\sst cl},\BR)$.

The following
Proposition \ref{Qkashdecomp2} describes how the quantized Teichm\"uller spaces
are related to the quantized Kashaev space. 

\begin{propn}\label{Qkashdecomp2}
There exists a unitary operator $\SI_\si$,
\[
\SI_\si:\CK(\vf_\si)\ra \CH_z(\si)\ot\CH_h(\si)
\] 
such that
\begin{equation*}\label{zhmap} \begin{aligned}
&\SI_\si^{}\cdot\hat{\spp}_{c}^{}\cdot\SI_\si^{-1}=
{\spp}_{c}\ot 1\,,\\
&\SI_\si^{}\cdot\hat{\sq}_{c}\cdot\SI_\si^{-1}=
{\sq}_{c}\ot 1\,,
\end{aligned}\quad
{\rm and}
\quad\SI_\si^{}\cdot{\sh}_c^{}\cdot\SI_\si^{-1}=1\ot\sh_c'\,,
\end{equation*}
for any $e\in\si_\1$ and $c\in H_1(\Sigma,\BR)$, respectively.
\end{propn}

Let us recall that a move $m\in[\tau_{m},\si_{m}]$, $m\in\CM_\1(\Sigma)$
is admissible if 
both $\tau_{m}$ and $\si_{m}$ are admissible. 
Given an element $m\in\CM_\1^{\rm ad}(\Sigma)$
it is natural to consider the corresponding fat graphs
$\theta_m\df \vf_{\csi_m}$,
 $\hvf_m\df \vf_{\hsi_m}$ on $\Sigma$, to pick a path 
$\pi_m\in[\cvf_m,\hvf_m]$ and consider the  operator
$\su(m)\equiv\su(\pi_m)$. 
The reduction to the quantized Teichm\"uller spaces 
proceeds as in Subsection \ref{reduction}. By multiplying the
operators $\SI_\si^{}\cdot \su(m)\cdot\SI_\si^{-1}$ 
with suitably chosen operators $\SH(m)$ one 
gets operators $\sv(m)$ which factorize as $\sv(m)=\sw(m)\ot 1$.
The resulting
operators $\sw(m):\CH_z(\hvf_m)\ra \CH_z(\cvf_m)$ will 
then generate a unitary projective representation of the path groupoid of
$\CM^{\rm ad}(\Sigma)$.

\section{Geodesic length operators\index{geodesic length operators}} 
\label{lengthsec}
\setcounter{equation}{0}

Of fundamental importance for us will be to 
define and study quantum analogs of the geodesic length functions
on the Teich\-m\"uller spaces, the geodesic length operators.

\subsection{Overview}

When trying to define operators which represent the geodesic length
functions one has to face the following difficulty:
The
classical expression for $L_{\vf,c}\equiv
2\cosh\frac{1}{2}l_c$ as given by formula \ref{glength}
is a linear combination of monomials
in the variables $e^{\pm\frac{z_e}{2}}$ of a very particular form,
\begin{equation}\label{Lclass}
L_{\vf,c}\,
=\,
\sum_{\tau\in\CF}\,C_{\vf,c}(\tau)\,e^{x(\tau)}\, ,\quad x(\tau)\,\df
\,\sum_{e\in\vf_\1}\,\tau(e)\,z_e,
\end{equation}
where the summation is taken over the space $\CF$ of all maps 
$\vf_\1\ni e\ra \tau(e)\in\frac{1}{2}\BZ$. The coefficients 
$C_{\vf,c}(\tau)$ are positive 
integers, and non-vanishing for a finite number of $\tau\in\CF$
only.

In the quantum case one is interested in the definition 
of length operators $\SL_{\vf,c}$ which should be representable
by expressions similar to \rf{Lclass},
\begin{equation}\label{Lquant}
\SL_{\vf,c}\,
=\,
\sum_{\tau\in\CF}\,C_{\vf,c}^b({\tau})\,e^{\sx(\tau)}\, ,\quad \sx(\tau)\,\df
\,\sum_{e\in\vf_\1}\,\tau(e)\,\sz_e.
\end{equation}

The following properties seem to be indispensable if one wants to
interpret an operator of the general form \rf{Lquant}
as the quantum counterpart 
of the functions $L_{\vf,c}=2\cosh\frac{1}{2}l_c$:
\begin{enumerate}
\item[(a)] {\bf Spectrum:} $\SL_{\vf,c}$ is self-adjoint.
The spectrum of $\SL_{\vf,c}$ is simple and equal
to $[2,\infty)$. This is necessary and sufficient for the existence
of an operator $\sll_{\vf,c}$ - the {\it geodesic length operator} - 
such that 
$\SL_{\vf,c}=2\cosh\frac{1}{2}\sll_c$.
\item[(b)] {\bf Commutativity:} 
\[
\big[\,\SL_{\vf,c}\,,\,\SL_{\vf,c'}\,\big]\,=\,0\quad
{\rm if}\;\; c\cap c'=\emptyset.
\]
\item[(c)] {\bf Mapping class group invariance:}
\[ 
\sa_\mu(\SL_{\vf,c})\,=\,\SL_{\mu.\vf,c},
\quad\sa_\mu\equiv\sa_{[\mu.\vf,\vf]},\quad
\text{for all}\;\;\mu\in{\rm MC}(\Sigma).
\]
\item[(d)] {\bf Classical limit:} 
The coefficients $C_{\vf,c}^b({\tau})$ which appear in \rf{Lquant}
should be  deformations of the 
classical coefficients $C_{\vf,c}^{}({\tau})$ in the sense that
\[ \lim_{b\ra 0}C_{\vf,c}^b({\tau})=C_{\vf,c}^{}({\tau}).\]
\end{enumerate}
Length operators were first defined and studied
in the pioneering work \cite{CF2}. 
It was observed in \cite{CF2} that the necessary deformation of the
coefficients $C_{\vf,c}^b({\tau})$ is indeed nontrivial in general.
However, a full proof that the length operators introduced
in \cite{CF2} fulfil the requirements (a) and (c) does not seem to 
be available yet. We will therefore present an alternative 
approach to this problem,
which will allow us to define length operators that
satisfy (a)-(d) in full generality.

\subsection{Construction of the length operators}

Our construction of the length operators will proceed in two
steps. First, we will construct length
operators $\SL_{\si,\ga}$ in the case that the 
fat graph $\vf$ under consideration equals $\vf_\si$.
This will facilitate the verification of the properties (a)-(d)
formulated above. In order to define the length operators
$\SL_{\vf,\ga}$ in the general case we shall then simply
pick any marking $\si$ such that the given
curve $c$ is contained in the cut system $\CC_\si$, and define
\begin{equation}\label{lengthgenvf}
\SL_{\vf,\ga}\,\equiv\,\sa_{[\vf,\vf_\si]}(\SL_{\si,\ga}).
\end{equation}
Independence of this construction from the choice of $\si$ will follow 
from Theorem \ref{lengthindep} below.
Definition \rf{lengthgenvf} reduces the proof of properties (a), (b)
to the proof of the corresponding statements for 
the length operators $\SL_{\si,\ga}$ which will be given below.
Property (c) follows from $\sa_{[\mu.\vf,\vf]}\circ\sa_{[\vf,\vf_\si]}=
\sa_{[\mu.\vf,\vf_\si]}$.

In order to prepare for our construction of length operators it is useful
us recall the construction of the fat graph $\vf_\si$
in section \ref{pt->MS}.

%

\begin{defn}\label{lengthdef} $\quad$\begin{itemize}
\item[$$] {\bf Case $c\in\CA_\si$:}
Let $A$ be an annular neighborhood of
$c$ such that the part of $\vf_\si$ which is contained in 
$A$ is isotopic to the model depicted in Figure \ref{cycle}. 
We will
then define 
\begin{equation}\label{qlodef}
\SL_{\vf,\ga}\;\df\;e^{-2\pi b\sq_c}+
2\cosh2\pi b\spp_c\;.
\end{equation}
\item[$$] {\bf Case $c\notin\CA_\si$:} If a curve $c$ is not 
contained in $\CA_\si$, it is necessarily the outgoing boundary component 
of a trinion $T_p$ (cf. Figure \ref{subst1}).
Let $c_\ep$, $\ep=1,2$ be the curves
which represent the incoming boundary components
of $T_p$  as 
indicated in Figure \ref{threefat}. 
Given that 
$\SL_{\si;c_i}$, $i=1,2$ are already defined we will define
$\SL_{\si;c}$ by 
\begin{equation}\label{modLlem2}
\SL_{\si,c}\;=\;2\cosh(\sy_{c_\2}+\sy_{c_\1})
+e^{-\sy_{c_{\2}}}\SL_{\si,c_\1}
+e^{\sy_{c_\1}}\SL_{\si,c_\2}
+e^{\sy_{c_\1}-\sy_{c_\2}}\, ,
\end{equation}
where $\sy_{c_{\ep}}$, $\ep=1,2$
are defined as
$\sy_{c_\2}=2\pi b {\sq}_{c}$, 
$\sy_{c_\1}=- 2\pi b {\spp}_{c}$.
\end{itemize}\end{defn}
It is easy to see that this recursively defines length operators
for all remaining $c\in\CC_\si$.

\begin{rem} 
Let us note that \rf{CCR} implies that $[\sy_{c_\ep},
\SL_{\si,c_{\ep'}}]
=0$ for
$\ep,\ep'\in\{1,2\}$. We therefore do not have an 
issue of operator ordering in \rf{modLlem2}.
\end{rem}

\begin{propn}\label{lengthopcomm}
The length operators  $\SL_{\si,c}$, 
$\SL_{\si,c'}$ associated to different curves $c,c'\in\CC_\si$
commute with each other. 
\end{propn}
\begin{proof}
Let us recall that cutting the surface $\Sigma$ along all of the
curves $c\in\CA_\si$, yields a set of connected
components which all have genus zero.
For a given curve
$c\in\CC_{\si}\setminus\CA_\si$ let $\Sigma_c$ be the connected component 
which contains $c$. Cutting $\Sigma_c$ along $c$ produces two connected
components. The component which has $c$ as its outgoing
boundary component will be denoted $\Sigma'_c$.
It follows from Definition \ref{lengthdef}
that $\SL_{\si,c}$ is an operator function of the operators
${\spp}_{d}$ and ${\sq}_d$, where $d\in\CC_{\si}$ 
is contained in $\Sigma'_c$.
The claim therefore follows immediately 
from \rf{CCR} if $\Sigma_{c}'$ and $\Sigma_{c'}'$ 
are disjoint. 

Otherwise we have the situation that one 
of $\Sigma_{c}'$, $\Sigma_{c'}'$, say $\Sigma_{c'}'$ is a subsurface of
the other. The crucial point to observe is that
the resulting expression for $\SL_{\si,c}$ depends on the variables
${\spp}_{d}$ and ${\sq}_d$  
associated to the subsurface $\Sigma_{c'}'$
exclusively via $\SL_{\si,c'}$. The claim therefore 
again follows from  \rf{CCR}.
\end{proof}

The following theorem expresses the consistency of our
definition with the automorphisms induced by a change of the
marking $\si$. 
\begin{thm}\label{lengthindep} For a given curve $\ga$
let $\si_i$, $i=1,2$ be markings such that
$\ga$ is contained in both cut systems
$\CC_{\si_\1}$ and $\CC_{\si_\2}$. We then have
\begin{equation}\label{Inveqn}
\sa_{[\vf_{\si_\2},\vf_{\si_\1}]}(\SL_{\si_\1,\ga})\,=\,\SL_{\si_\2,\ga.}
\end{equation}
\end{thm}

{\it On the proof of Theorem \ref{lengthindep}:} 
The description of the modular groupoid
in terms of generators and relations given in section \ref{modgroup} reduces 
the proof of Theorem \ref{lengthindep} to the case that
$\si_\2$ and $\si_\1$ are connected by one of the elementary moves $m$
defined in section \ref{MS-gens}. In order to reduce the proof
of  Theorem \ref{lengthindep} to a finite number of verifications
one would need to have simple standard choices for 
the paths $\pi_m\in[\vf_{\si_\2},\vf_{\si_\1}]$ for all 
elementary moves $m$. Existence of such 
standard paths $\pi_m$ turns out to be nontrivial, though. 
The task to find suitable paths $\pi_m$ is particularly simple for
a subclass of moves $m$ which is defined as follows.
\begin{defn}\label{locreal} $\quad$
\begin{itemize}
\item[(i)] Let $\pi\in[\vf',\vf]$ be a path 
in the complex ${\CP t}(\Sigma)$ which 
is described by
a sequence $S_{\pi}\df(\vf'\equiv
\vf_n~,~\dots~,~\vf_1\equiv
\vf)$ of fat graphs such that consecutive elements
of $S_{\pi}$  are connected
by edges  in ${\CP t}_\1(\Sigma)$. We will say that
$\pi\in[\vf',\vf]$ is realized locally in 
a subsurface $\Sigma'\hookrightarrow\Sigma$ if 
the restrictions of $\vf_i$ to $\Sigma\setminus\Sigma'$
coincide for all $i=1,\dots,n$.
\item[(ii)]
We will say that the move $m=[\tau_{m},\si_{m}]\in
\CM_\1'(\Sigma)$
can be realized locally
if there exists a path $\pi\in[\theta_m,\hvf_m]$ in
the complex ${\CP t}(\Sigma)$ that is realized locally in 
$\Sigma_m$ in the sense of (i).
\end{itemize}
\end{defn}

For moves $m=[\tau_{m},\si_{m}]\in\CM_\1(\Sigma)'$  
that can be realized locally we may 
choose essentially the {\it same}
path $\pi_m\in[\theta_m,\hvf_m]$ for all surfaces 
$\Sigma$ into which $\Sigma_m$ can be embedded. 
It is then crucial to observe the
following fact
\begin{propn}\label{locrelfact}
If $m\in\CM_\1'(\Sigma)$ is admissible, but can not be realized locally,
there always exists a path $\varpi_m$ which is (i) homotopic to
$m$ within $\CM^{\rm ad}(\Sigma)$, and (ii) takes the
form 
\begin{equation}
\varpi_m\,=\,Y_m^{\phantom{1}}\circ m\circ Y_m^{-1},
\end{equation}
where $Y_m$ is a chain composed out of $Z_p$-moves and $F_{pq}$-moves
which can all be realized locally.
\end{propn}
The proof of Proposition \ref{locrelfact} is given in Appendix 
\ref{locrelfactpf}. 
It therefore suffices to prove Theorem \ref{lengthindep} in the case
that $\si_\2$ and $\si_\1$ are connected by any elementary move $m$
that can be realized locally. This amounts to
a finite number
of verifications which can be carried out by straightforward, but 
tedious calculations. Some details are given in the Appendix 
\ref{lengthindeppf}. 
\hfill $\square$

\subsection{Spectrum} 

\begin{thm} \label{Lspecthm}
The spectrum of $\SL_{\si,c}$ is simple and equal to $[2,\infty)$.
\end{thm}
\begin{proof}
To begin with, let us consider the following simple model for the length
operators:
\begin{equation}
\SL\,\df\,2\cosh 2\pi b\spp+e^{-2\pi b\sq},
\end{equation}
where $\spp$, $\sq$ are operators on 
$L^2(\BR)$ which satisfy the commutation
relations $[\spp,\sq]=(2\pi i)^{-1}$.

A basic fact is that $\SL$ is self-adjoint.
Indeed, being a sum of two {\it positive} 
self-adjoint operators $\SL$ is self-adjoint on the 
intersection of the domains of the summands. The
main spectral properties of this operator are summarized in 
the following proposition. 

\begin{propn} \label{kathm}$\;\frac{\quad}{}$
\cite{Ka4} $\;\frac{\quad}{}$
\begin{itemize}
\item[(i)] We have $
{\rm Spec}(\SL)=(2,\infty)$.
\item[(ii)] The spectrum of $\SL$ in $L^2(\BR)$ is simple.
\end{itemize}
\end{propn}

Validity of Theorem \ref{Lspecthm}
in the case $c\in\CA_\si$ is a direct consequence
of Proposition \ref{kathm}. It remains to treat the case $c\notin\CA_\si$.
We will keep 
the notations introduced in Definition \ref{lengthdef}.
The main ingredient will be an operator $\SC_{\si,c}:\CH(\si)\ra\CH(\si)$ 
which maps all length operators $\SL_{\si,c}$
to the simple standard form $\SL_{\si,c}^{\rm st}$,
\begin{equation}
\SL_{\si,c}^{\rm st}\;=\; 
2\cosh 2\pi b{\spp}_{c}+ e^{-2\pi b{\sq}_{c}}
\end{equation}
in the sense that the following commutation relations are satisfied:
\begin{equation}\label{Ltostandard'}
\SC_{\si,c}\cdot\SL_{\si,c}^{}\;=\;
\SL_{\si,c}^{\rm st}
\cdot\SC_{\si,c}.
\end{equation}

\begin{defn}\label{Cdef}
Let the unitary operator $\SC_{\si,c}$ be defined by
\begin{equation}\label{SCdef}
\SC_{\si,c}^{-1}\;=\;e_b\big({\sq}_{c}-\mss_{{\2}}\big)
\frac{s_b\big(\mss_{{\1}}-{\spp}_{c}\big)}
     {s_b\big(\mss_{{\1}}+{\spp}_{c}\big)}
e^{2\pi i \mss_{{\2}}{\sq}_{c}}\,,
\end{equation}
where
\[
\mss_\ep\;=\;(2\pi b)^{-1} \,{\rm arcosh}\fr{1}{2}
\SL_{\si,c_\ep}\,\quad \ep=1,2 ,
\] 
and the special function $s_b(x)$ is a close relative of 
$e_b(x)$ defined in the Appendix A.
\end{defn}

\begin{lem}\label{SClem}
The unitary operator $\SC_{\si,p}$ satisfies \rf{Ltostandard'}.
\end{lem}
\begin{proof} 
Proposition \ref{SClem} is proven by means of 
a direct calculation using the 
explicit form of $\SL_{\si,p}$ given in \rf{modLlem2} and
the functional equations \rf{sb_feq} and
\rf{eb_feq}.
\end{proof}
The proof of Theorem \ref{Lspecthm} is thereby reduced to Proposition
\ref{kathm}.
\end{proof}

\subsection{Relation with the Dehn twist\index{Dehn twist} generator}

To round off the picture we shall now discuss, following \cite{Ka3,Ka4},
the relation between the length operators $\SL_{\vf,c}$ and the 
operator $\SD_{\vf,c}$ which represents the Dehn twist $D_{\ga}$ along $c$.

\begin{figure}[htb]
\epsfxsize3cm
\centerline{\epsfbox{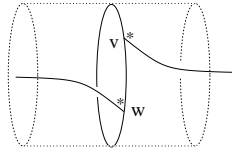}}
\caption{Annulus $A_c$ and fat graph $\vf$ on $A$.}
\label{cycle2}\end{figure}

A closed curve $c$ will be called a curve of simple type
if the the connected components of $\Sigma\setminus c$ all have more than 
one boundary component. It is not hard to see (using the construction in 
Section \ref{pt->MS}, for example) that for curves of 
simple type there always exists a fat graph $\vf$ and an annular 
neighborhood $A_c$ of $c$ in which $\vf$ takes the form depicted in
Figure \ref{cycle2}.

\begin{lem}\label{MCgenpropn}
Let $\Sigma$ be a surface with genus $g$ and $s$ holes.
The pure mapping class group ${\rm MC}(\Sigma)$ is generated by
the Dehn twists along geodesics of simple type.
\end{lem}

\begin{proof}
For $g\geq 2$ it is known that the Dehn twists along non-separating 
closed curves suffice to generate the mapping class group ${\rm MC}(\Sigma)$
\cite{G}. Closed curves which are non-separating are always of simple type.

For $g=0$ any closed curve is of simple type. In the remaining 
case $g=1$ one may note that the only closed curve $c$ which is 
not of simple type is the one which separates a one-holed
torus from the rest of $\Sigma$. It is then well-known
that the Dehn twists along $a$- and $b$-cycles 
of the one-holed torus generate the Dehn-twist along $c$. 
If we supplement these generators by the Dehn twists along the 
remaining closed curves (which are all of simple type), we get
a complete system of generators for ${\rm MC}(\Sigma)$ \cite{G}.
\end{proof}

We may therefore use a fat graph
which in an annular neighborhood of $c$ takes the 
form depicted in Figure \ref{cycle2}.
It is easy to see
that the action of the Dehn twist $D_{\ga}$
on the fat graph $\vf$ can be undone by a single flip $\omega_{vw}$.
As the representative $\SF_{\vf,\ga}$ for the Dehn twist $D_{\ga}$
on $\CK(\vf)$ 
we may therefore choose $\SF_{\vf,\ga}=\ST_{vw}$.

The operator  $\sh_{\vf,c}$
associated to the homology cycle $c$ is 
$\sh_{\vf,c}\df\frac{1}{2}(\spp_v+\sq_w)$. It is not very difficult to
verify that the operator 
\begin{equation}
\SD_{\vf,\ga}\;=\;\zeta^{-6}\exp(2\pi i \sh_{\vf,\ga}^2)\,\SF_{\vf,\ga}\,.
\end{equation}
commutes with all $\sh_c$, $c\in H_1(\Sigma,\BR)$.
The prefactor $\zeta^{-6}$, $\zeta=e^{\pi i b^2/3}$ 
was inserted to define a convenient 
normalization.
We then have the following result.
\begin{propn} \label{Dehnev}\cite{Ka4} \;\,$\frac{\quad}{}$
$\SD_{\ga}$ coincides with the following function of the length operator
$\sll_{\ga}$:
\[ \SD_{\ga}\;=\;\zeta^{-6}\exp\biggl( i \frac{\sll_{\ga}^2}{8\pi b^2}\biggr).
\]
\end{propn}
This should be compared with the classical result that
the geodesic length functions are the Hamiltonian generators
of the Fenchel-Nielsen twist flow, which reproduces the 
Dehn twist for a twist angle of $2\pi$.

\section{Passage to the length representation}\label{lengthrep}
\setcounter{equation}{0}

\subsection{The length representation}

Our aim is to define a representation for $\CH(\Sigma)$ which
is such that the 
length operators associated to a cut system
are all realized as multiplication operators.
For a surface $\Sigma$ and a marking 
$\sigma$ on $\Sigma$ let 
\begin{equation}\label{HMS}
\CH_{\rm L}^{}(\si)\;\df\;
\CK_{\rm\sst sp}^{\ot\si_\1},
\end{equation}
where $\CK_{\rm\sst sp}\simeq L^2(\BR_+,\eta)$, and
$\eta$ is the spectral measure of the 
operator $\sll_{\rm A}=2{\rm arcosh}\frac{1}{2}\SL_{\rm A}$.
The numbering of the edges associated to $\si$ defines
canonical isomorphisms
\begin{equation}
\CH_{\rm L}^{}(\si)
\;\simeq\; \CH_{\rm L}^{}(\Sigma)\,\df\,
L^2({\mathfrak L},d\eta_{\mathfrak L}),
\end{equation}
where ${\mathfrak L}\simeq \BR_+^{3g-3+2s}$ and 
$d\eta_{\mathfrak L}=d\eta(l_1)\land\dots\land d\eta(l_{3g-3+2s})$ 
is the corresponding product measure.
For $e\in\si_{\1}$ and $f:\BR^+\ra \BC$ let us define the
multiplication operator $\sm_e[f]$ by
\begin{equation}
\big(\sm_e[f]\Psi\big)(\la_{\si})\;\df\;f(l_e)\Psi(\la_{\si}).
\end{equation}
The operator $\sm_e[f]$ will represent a bounded operator 
$\sm_e[f]:\CH_{\rm L}^{}(\si)\ra \CH_{\rm L}^{}(\si)$
iff $f\in L^{\infty}(\BR^+,d\eta)$.

For the rest of this section let us 
consider an  admissible marking $\si$ and the
associated fat graph $\vf_{\si}$.
The Definition \ref{lengthdef} yields a family of mutually
commuting length operators $\SL_{\si,c}$ associated to
the curves $c$ which constitute the cut sytem $\CC_\si$. 
It follows from the spectral theorem\footnote{See Appendix B 
for the precise statement.} 
for the family $\{\SL_e;e\in\si_{{\1}}\}$ of self-adjoint operators
that 
there exists a unitary operator 
\[
\SD(\si)\;:\; \CH_{z}({\si})\longrightarrow
\CH_{\rm L}^{}(\si)
\]
that diagonalizes the set of length operators $\SL_{\si,e}$, 
$e\in \si_{{\1}}$ in the sense that
\begin{equation}\label{diag}
\SD(\si)\cdot\SL_{\si,e}\;=\;\sm_e[2\cosh\fr{1}{2}l]\cdot \SD(\si).
\end{equation}

Our aim for the rest of this section will be to give a recursive 
construction for such an operator $\SD(\si)$ in terms of
operators $\SC_{p}$ which are associated to the vertices
$p\in\si_\0$ of $\si$.

\subsection{Construction of the operator $\SD(\si)$}\label{Cprimedef}

The main ingredient will be an operator $\SC_\si:\CH(\si)\ra\CH(\si)$ 
which maps all length operators $\SL_{\si,c}$
to the simple standard form $\SL_{c}^{\rm st}$,
\begin{equation}
\SL_{\si,c}^{\rm st}\;=\; 
e^{-2\pi b{\sq}_c}+2\cosh 2\pi b{\spp}_c.
\end{equation}
The operators $\SC_{\si,c}$ that were defined in 
Definition \ref{Cdef} solve this task locally for every curve
$c\in\CC_{\si}\setminus\CA_\si$. If the trinion $T_p$ has $c_p$ as the 
outgoing boundary component we will define
\begin{equation}
\SC_{\si,p}\,\equiv\,\SC_{\si,c_p}.
\end{equation}
We are now ready to define 
\begin{equation}\label{SCdef2}
\SC_{\si}\;\df\;\prod_{p\in\si\!{}_{\0}}\,\SC_{\si,p}\,
.\end{equation}
Let us note that we do not have to specify the order in
which the operators $\SC_{\si,p}$ appear 
thanks to the following lemma.
\begin{lem}
We have $\SC_{\si,p}\SC_{\si,q}=\SC_{\si,q}\SC_{\si,p}$.
\end{lem} 
\begin{proof} This  
follows from Proposition \ref{lengthopcomm} and
Definition \ref{lengthdef}, keeping in mind 
\rf{CCR}.
\end{proof}

In order to construct the sought-for operator $\SD(\si)$ it now 
remains to map the length operators $\SL_{\si,e}^{\rm st}$
to multiplication operators. 
Proposition \ref{kathm} ensures existence of an operator 
\[
\sd_{\si,e}:L^2(\BR)\ra \CK_{\rm\sst sp}\quad{\rm such\;\, that}\quad
\sd_{\si,e}\cdot\SL_e^{\rm st}=\sm_e[2\cosh\fr{1}{2}l]\cdot\sd_{\si,e}.
\]
Out of $\SE_{\si}$ we may finally define the operator 
$\SD(\si)$ as
\begin{equation}\label{SDdef}
\SD(\si)\;\df\;\;\sd_\si
\cdot\SC_{\si}^{}\,, 
\qquad\sd_\si\,\df\,
\bigotimes_{e\in \si_{{1}}}
\sd_{\si,e}\;,
\end{equation}
where $\SI_\si^{}$ is the operator introduced in Proposition \ref{Qkashdecomp}.
It is straightforward to 
verify that the operator $\SD(\si)$ indeed satisfies the 
desired property \rf{diag}.

\section{Realization of ${\rm M}(\Sigma)$\index{modular groupoid}}
\label{Mconstr}
\setcounter{equation}{0}

\subsection{Two constructions for the generators}

\renewcommand{\csi}{\tau}
\renewcommand{\hsi}{\si}
\renewcommand{\cvf}{\theta}
\renewcommand{\hvf}{\vf}
Our aim is to define operators $\SU(m)$ associated to the 
edges $m\in\CM_\1(\Sigma)$. 
We will give two constructions for these operators, each of which makes
certain properties manifest. The proof of the equivalence of these
two constructions
will be the main difficulty that we will have to deal with.

\paragraph{The first construction.}

Let us recall that a move $m\in[\tau_{m},\si_{m}]$, $m\in\CM_\1(\Sigma)$
is admissible if 
both $\tau_{m}$ and $\si_{m}$ are admissible. 
Given an element $m\in\CM_\1^{\rm ad}(\Sigma)$
it is natural to consider the corresponding fat graphs
$\theta_m\df \vf_{\csi_m}$,
 $\hvf_m\df \vf_{\hsi_m}$ on $\Sigma$, to pick a path 
$\pi_m\in[\cvf_m,\hvf_m]$ and define operators
$\widetilde{\SU}(m):\CH_{\rm L}^{}(\hsi_m)\ra\CH_{\rm L}^{}(\csi_m)$ as
\begin{equation}\label{SUdef2}
\widetilde{\SU}(m)\;\df\;\SD(\csi_m)\cdot\su(\pi_m)\cdot 
\SD(\hsi_m)^{\dagger}\;.
\end{equation}

\paragraph{The second construction.}

We note that for all $m\in\CM_\1(\Sigma)$ 
the markings $\csi_m$ and $\hsi_m$ will coincide outside of 
a subsurface $\Sigma_m\hookrightarrow\Sigma$.
Let us therefore consider the restrictions
$\csi'_m$ and 
$\hsi'_m$ of $\csi_m$ and $\hsi_m$ to $\Sigma_m$
respectively. Admissibility of $\csi'_m$ and 
$\hsi'_m$ is obvious for $m\in\CM_\1(\Sigma)$,
allowing us to use the first construction in order to define
an operator 
\[
\SU'(m):
\CH_{\rm L}(\hsi_m')\ra\CH_{\rm L}(\csi_m')\,.
\]
Out of $\SU'(m)$ we may then construct the sought-for operator
$\SU(m)$ by acting with $\SU'(m)$ non-trivially only
on those tensor factors of $\CH_{\rm L}(\hsi)=
\CK_{\rm\sst sp}^{\ot\hsi_\1}$ which correspond to the subsurface 
$\Sigma_m\hookrightarrow\Sigma$. More precisely,
let $\RE_m\subset \si_\1$ be the set of edges in $\si_\1$ 
which have nontrivial intersection with $\Sigma_m$.
Out of $\SU'(m)$ let us then define the operator
$\SU(m):\CK(\hsi_m)\ra \CK(\csi_m)$ 
by applying definition \rf{legnotation} to the case 
$\SO\equiv\SU'(m)$ and $\BJ\equiv\RE_m$.

\paragraph{Comparison.}

The crucial difference between $\widetilde{\SU}(m)$ and $\SU(m)$
is that the latter is manifestly acting {\it locally} in $\CH_{\rm L}(\si)$, 
in the sense that it acts only on the tensor factor of $\CH_{\rm L}(\si)$
which corresponds to the subsurface $\Sigma_m$. This is not obvious
in the case of $\widetilde{\SU}(m)$. 

\medskip

The length operators $\SL_{\si,c}$ associated to the boundary components
$c\in A(\Sigma)$ form a commutative family of operators. The 
joint spectral decomposition for this family of operators
leads us to represent $\CH_{\rm L}(\Sigma)$ as
\begin{equation}\label{bdspecdecomp}
\CH_{\rm L}(\Sigma)\,\simeq\,\int_{\FL}d\eta_{\Sigma}^{}(\fc)\;
\CH_{\rm L}(\Sigma,\fc),
\end{equation}
where the integration is extended over the set $\FL\simeq \BR_+^{s}$ of 
all colorings $\fc$ of the boundary by
elements of $\BR_+$. It follows from Theorem \ref{lengthindep}
that the operators $\SU(m)$, $\widetilde{\SU}(m)$ 
commute (up to permutations
of the boundary components) with the length operators
$\SL_{\si,c}$, $c\in A(\Sigma)$. Within the representation
\rf{bdspecdecomp} one may therefore\footnote{According
to Proposition \ref{decomppropn} in Appendix B}
represent  $\SU(m)$ and $\widetilde{\SU}(m)$ 
by families of operators
$(\SU(m,\fc))_{\fc\in \FL}$ and $(\widetilde{\SU}(m,\fc))_{\fc\in \FL}$,
respectively. 

\subsection{The main result} \label{prooflocalprop1}

The following theorem is the main result of this paper.

\begin{thm}\label{mainthm} 
The operators $\SU(m,\fc)$, $\fc\in \FL$, 
generate a tower of projective 
unitary representations
of the modular groupoids ${\rm M}(\Sigma)$.
\end{thm}
The proof of Theorem \ref{mainthm} will take up the rest of   this 
subsection.


To begin with, let us note that 
the necessary structure \rf{CHtriniondecomp}
of the Hilbert spaces follows trivially from our definition 
of the length representation in Section \ref{lengthrep}, where
in the present case we simply have $\CH(S_\3,\fc_\3)\simeq \BC$.
It is furthermore
clear that the operators $\widetilde{\SU}(m,\fc)$, $\fc\in \FL$
generate a unitary projective representation of the modular
groupoid ${\rm M}(\Sigma)$ for each surface $\Sigma$ within the
considered class. Let us finally note that the 
naturality properties formulated
in Subsection \ref{towerdef} clearly
hold for the operators ${\SU}(m,\fc)$, $\fc\in \FL$. 
Our main task is therefore to show that
$\widetilde{\SU}(m,\fc)={\SU}(m,\fc)$, as will be 
established in Proposition
\ref{localprop1} below.

\begin{propn}
\label{localprop1} 

For all $m\in\CM_\1^{\rm ad}(\Sigma)$ there exists 
a path $\pi_m\in[\cvf_m,\hvf_m]$ such that we have
\[\widetilde{\SU}(m,\fc)\,=\,\SU(m,\fc)\,.
\]
\end{propn}

\begin{proof}
When we compare the respective definitions of $\SU(m)$
and $\widetilde{\SU}(m)$, we observe that there are 
two main discrepancies that 
we need to deal with.
First, it is not always true that the path $\pi_m$
can be realized locally in the sense of Definition \ref{locreal}. 
It may therefore
not be clear a priori why there should exist a simple
relation between $\SU(m)$ and $\widetilde{\SU}(m)$.

Second, we may observe that the definition of the operators $\SU(m)$
and $\widetilde{\SU}(m)$ involves products of operators
$\SC_{\si,p}$, where the set of vertices $p$ that the product is
extended over is generically much smaller in the case of $\SU(m)$.
This means that most of the factors $\SC_{\si,p}$ must 
ultimately cancel each other 
in the expression for $\widetilde{\SU}(m)$. 
The first step will be to prove Proposition \ref{localprop1}
in the case that $m$ can be realized locally.

\begin{lem}\label{locallem1}
We have
\begin{equation}\label{localact}
\widetilde{\SU}(m)\,=\,\SU(m)\,
\end{equation}
whenever $m=[\tau_{m},\si_{m}]\in\CM_\1'(\Sigma)$  
can be realized locally.
\end{lem}

\begin{proof}
It follows from  \rf{SCdef} and \rf{CCR} that $\SC_{\si,p}$ 
can be represented as a function of the following 
operators 
\[
\SC_{\si,p}=\SC_{\si,p}
\big(\,\spp_{v_p},\sq_{v_p}\,;\,
\SL_{c_\2},\SL_{c_\1}\,\big)\,,
\]
where $v_p$ is the vertex of $\vf_\si$ contained in the trinion $T_p$ 
and $c_{\ep}\equiv c_{\ep}(p)$, $\ep=1,2$ are the curves which represent
the outgoing and incoming boundary components of $T_p$
respectively.
%
Let us recall that the length operators
$\SL_{\csi_m,c_{\ep}}$ and $\SL_{\hsi_m,c_{\ep}}$
$\ep=1,2$ satisfy \rf{Inveqn}. 
These observations imply that 
\begin{equation}\label{SCcommrel}
\su(\pi_m)\cdot \SC_{\hsi_m,p}  \;=\;
\SC_{\csi_m,p}\cdot \su(\pi_m),
\end{equation}
whenever the operators $\spp_{\si,p}$ and $\sq_{\si,p}$ commute with 
$\su(\pi_m)$. Our task is therefore to determine the set of
all $p\in\si_\0$ for which this is the case.
The condition that $m$ can be realized locally 
implies that \rf{SCcommrel} will hold unless
$p$ is located within $\Sigma_m$. The claim
now follows straightforwardly from these observations.\end{proof}

In order to treat the general case let us recall that
Proposition \ref{locrelfact} implies that $\pi_m$ may be chosen as $\pi_{m}=
y_m\circ \hat{\pi}_m\circ y_m^{-1}$, where $y_m$ is
uniquely defined by the factorization of $Y_m$ into elementary moves that
can be realized locally, and $\hat{\pi}_m$ is the (fixed) path which
was chosen to represent $m$ in the case that $m$ can be realized locally.
This leads to the following representation for 
$\widetilde{\SU}(m)$:
\begin{equation}\label{Vtildefactor}\begin{aligned}
\widetilde{\SU}(m)& \,=\,\SD(\csi_m)\cdot\su(\pi_m)\cdot 
\SD(\hsi_m)^{\dagger}\\
& \,=\,\SD(\csi_m)\cdot\su(y_m)\cdot\su(\hat{\pi}_m)\cdot\su(y_m^{-1})\cdot 
\SD(\hsi_m)^{\dagger}\,.
\end{aligned}\end{equation}
By using Lemma \ref{locallem1} one may deduce
from \rf{Vtildefactor} that the following holds:
\begin{equation}\begin{aligned}
\widetilde{\SU}(m) & \,=\,\widetilde{\SU}(Y_m^{\phantom{1}})\cdot
\SU(m)\cdot\widetilde{\SU}(Y_m^{-1})\\
& \,=\,{\SU}(Y_m^{\phantom{1}})\cdot
\SU(m)\cdot{\SU}(Y_m)^{\dagger}\,.
\end{aligned}\end{equation}
It remains to observe that 
\begin{lem}
We have $\SU(m_\1)\cdot\SU(m_\2)\cdot\SU(m_\1)^{\dagger}=\SU(m_\2)$ whenever 
${\rm supp}(m_\1)\cap {\rm supp}(m_\2)=\emptyset$.
\end{lem}
Recall that we had
defined the  
notation ${\rm supp}(m)$ by ${\rm supp}(m)=\{p\}$ if $m=Z_p,B_p,S_p$
and ${\rm supp}(m)=\{p,q\}$ if $m=(pq),F_{pq}$.
\begin{proof}
Let us factorize $\CH_{\rm L}(\Sigma)$ as
$\CH_{\rm L}(\Sigma)\,=\,\CK_{\rm\sst sp}^{\ot\beta_\1}\ot
\CK_{\rm\sst sp}^{\ot\si_\1\setminus\beta_\1}
$, where $\beta_\1\subset\si_\1$ is the set of all edges which end in 
boundary components of $\Sigma$. 
The representation \ref{bdspecdecomp} may be rewritten as
\newcommand{\bc}{{\mathbf c}}
\begin{equation}\label{CHMrepr}
\CH_{\rm L}(\Sigma)\,\simeq\,
\int_{\BR_+^s}^{\oplus}d\eta_s(\bc)\;\CH_{\rm L}(\Sigma,\bc)\,,
\end{equation}
where $\CH_{\rm L}(\Sigma,\bc)\simeq 
\CK_{\rm\sst sp}^{\ot\si_\1\setminus\beta_\1}$ for all $\bc\in \BR_+^s$.
Let $\SU(m,\bc)$, $\bc\in\BR_+^{s(m)}$ be the  unitary
operators 
on $\CH_{\rm L}(\Sigma_m,\bc)$
which represent the operators $\SU(m)$ in the representation 
\rf{CHMrepr}.
Within this representation it becomes almost 
trivial to complete the proof of Proposition \ref{localprop1}.
Let $\be_{\2\1}=\be_\2\cup\be_\1$, where $\be_\jmath\subset \si_\1$, 
$\jmath=\1,\2$ 
are the sets of edges which correspond to boundary components of 
$\Sigma_{m_\jmath}$. Let $\de_\jmath\subset \si_\1$, 
$\jmath=\1,\2$ be the 
sets of edges that are fully contained in the interior of 
$\Sigma_{m_\jmath}$ respectively. Let finally 
$\si_\1'=\si_\1\setminus(\be_{\2\1}\cup\de_\2\cup\de_\1)$
We may then factorize $\CH_{\rm L}(\Sigma)$ in the following way:
\begin{equation}
\CH_{\rm L}(\Sigma)\,=\, \CK_{\rm\sst sp}^{\ot\beta_{\2\1}}\ot
\CK_{\rm\sst sp}^{\ot\de_\2}\ot
\CK_{\rm\sst sp}^{\ot\de_\1}\ot\CK_{\rm\sst sp}^{\ot\si_\1'}\;,
\end{equation}
where we define $\CK_{\rm\sst sp}^{\ot\de}=\BC$ if $\de=\emptyset$.
This may be rewritten as
\begin{equation}\label{factorrep21}
\CH_{\rm L}(\Sigma)\,\simeq\,
\int_{\BR_+^{s_{\2\1}}} d\mu(\bc_{\2\1})\;
\CH_{\rm L}(\Sigma,\bc_{\2\1}),
\end{equation}
where $\CH_{\rm L}(\Sigma,\bc_{\2\1})\simeq \CK_{\rm\sst sp}^{\ot\de_\2}\ot
\CK_{\rm\sst sp}^{\ot\de_\1}\ot\CK_{\rm\sst sp}^{\ot\si_\1'}$ for all
$\bc_{\2\1}\in \BR_+^{s_{\2\1}}$, $s_{\2\1}={\rm card}(\be_{\2\1})$.
In the representation \rf{factorrep21}
we may represent $\SU(m_\jmath)$, $\jmath=1,2$
by families of operators $\SU(m_\jmath,\bc_{\2\1})$ which take the
form
 $\SU(m_\2,\bc_{\2\1})\simeq \SU_\2(m_\2,\bc_{\2\1})\ot\id\ot\id$ and
$\SU(m_\1,\bc_{\2\1})\simeq \id\ot\SU_\1(m_\1,\bc_{\2\1})\ot\id$
respectively. The Lemma follows easily from these observations.
\end{proof}
It follows from the Lemma that
$\SU(Y_m^{\phantom{1}})\cdot
\SU(m)\cdot\SU(Y_m)^{\dagger}
\,=\,\SU(m)$, which completes the proof of
Proposition \ref{localprop1}. 
\end{proof}


Taken together our results show that the quantization of 
Teichm\"uller spaces yields a tower of unitary projective 
representations of the modular groupoid in the sense of 
\S\ref{towerdef}. 

\subsection{Concluding remarks}

Let us 
recall that the mapping class group is a group of symmetries
for both complexes $\CP t(\Sigma)$ and $\CM(\Sigma)$. Having
projective unitary representations of the associated 
groupoids ${\rm Pt}(\Sigma)$ and ${\rm MC}(\Sigma)$ will therefore 
induce corresponding representations of the mapping class group
by means of the construction in \S\ref{symmG}. It is not hard to see
It follows quite easily from Proposition \ref{Madconn}
that these two representations are equivalent.



Finally
it is clearly of interest to calculate the phases
in the relations $\su(\pi_{\varpi})=\zeta^{\nu(\varpi)}$. 

\begin{conj}
There exists a definition for the operators  $\SU(m)$, $m\in\CM_\1(\Sigma)$ such that
the phases $\zeta^{\hat{\nu}(\varpi)}$ which appear in the 
relations $\su(\pi_{\varpi})=\zeta^{\nu(\varpi)}$ are trivial for all
but one $\varpi\in \CM_\2(\Sigma)$, which can be chosen
as the relation \rf{onetor:a},a). The phase which appears in the 
relation \rf{onetor:a},a) is given as 
\[
\zeta^{\nu(\varpi)}=e^{\frac{\pi i}{2}c_{\rm L}}, \quad c_{\rm L}\equiv
1+6(b+b^{-1})^2.
\]
\end{conj}

We now believe to have a proof of this conjecture. Details
will appear elsewhere. Our conjecture is also 
strongly supported by the calculation in \cite{Ka2}
which establishes a similar result for the realization
of the {\it pure} 
mapping class group on the quantized Teichm\"uller spaces.



\newpage

\appendix

\section{The special functions $e_b(x)$ and $s_b(x)$}
\setcounter{equation}{0}\label{specapp}


The function $s_b(x)$ may be defined with the help of the 
following integral representation.
\begin{equation}
\log s_b(x)\;=\;\frac{1}{i}\int\limits_0^{\infty}\frac{dt}{t}
\biggl(\frac{\sin 2xt}{2\sinh bt\sinh b^{-1}t}-
\frac{x}{t}\biggl)\;\;.
\end{equation}
This function, or close relatives of it like 
\begin{equation}\label{ebdef}
e_b(x)\;=\;e^{\frac{\pi i}{2}x^2}\,e^{-\frac{\pi i}{24}(2-Q^2)}s_b(x)\;,
\end{equation}
have appeared in the literature under various names like 
``Quantum Dilogarithm'' \cite{FK1}, ``Hyperbolic G-function''
\cite{Ru}, ``Quantum Exponential Function'' \cite{W} 
and ``Double Sine Function'', we refer to 
the appendix of \cite{KLS} for
a useful collection of properties of $s_b(x)$ and further references.
The most important properties for our purposes are 
\begin{align}
{}& \text{(i) Functional equation:} \quad
s_b\big(x-i\fr{b}{2}\big)\;=\;2\cosh
\pi b x\;
s_b\big(x+i\fr{b}{2}\big). \label{sb_feq}\\
{}& \text{(ii) Analyticity:}\quad
s_b(x)\;\text{is meromorphic,}\nonumber\\ 
{}& \hspace{2.7cm}\text{poles:}\;\,  
x=c_b+i(nb-mb^{-1}), n,m\in\BZ^{\geq 0}.\\
{}& \hspace{2.7cm}\text{zeros:}\;\,  
x=-c_b-i(nb-mb^{-1}), n,m\in\BZ^{\geq 0}.\nonumber \\
{}& \text{(iii) Self-duality:}\quad s_b(x)=s_{1/b}(x).  \\
{}& \text{(iv) Inversion relation:}\quad s_b(x)s_b(-x)\;=\;1.\\
{}& \text{(v) Unitarity:} \quad \overline{s_b(x)}\;=\;1/s_b(\bar{x}).\\
{}& \text{(vi) Residue:} \quad {\rm res}_{x=c_b}s_b(x)
=e^{-\frac{\pi i}{12}(1-4c_b^2)}(2\pi i)^{-1}\label{sbRes}.
\end{align}
The function $e_b(x)$ clearly has very similar properties as 
$s_b(x)$. We list the properties that
get modified compared to $s_b(x)$ below.
\begin{align}
{}& \text{(i)' Functional equation:} \quad
e_b\big(x-i\fr{b}{2}\big)\;=\;(1+e^{2\pi b x})\;
e_b\big(x+i\fr{b}{2}\big). \label{eb_feq}\\
{}& \text{(iv)' Inversion relation:}\quad e_b(x)e_b(-x)\;=\;e^{\pi ix^2}
e^{-\frac{\pi i}{6}(1+2c_b^2)}.\\
{}& \text{(vi)' Residue:} \quad {\rm res}_{x=c_b}s_b(x)
=(2\pi i)^{-1}\label{ebRes}.
\end{align}

Among the most remarkable properties satisfied by the function 
$e_b(x)$ is the so-called pentagonal relation \rf{basicpentagon}
which underlies the validity of the pentagonal relation
\rf{pentrel}, 
\begin{equation}\label{basicpentagon}
e_b(\spp)\cdot e_b(\sq)\;=\;e_b(\sq)\cdot e_b(\spp+\sq)\cdot e_b(\spp),
\end{equation}
Relation \rf{basicpentagon} is valid if 
$\sq$ and $\spp$ represent $[\spp,\sq]=(2\pi i)^{-1}$ on $L^2(\BR)$.
Equation \rf{basicpentagon} in turn is equivalent to the 
following property of the function $E_b(x)\df e_b(-\frac{x}{2\pi b})$:
\begin{equation}\label{q-exp}
E_b(\SU)\cdot E_b(\SV)\,=\,E_b(\SU+\SV)\,,
\end{equation}
where $\SU=e^{2\pi b \sq}$, $\SV=e^{2\pi b \spp}$.
Proofs of \rf{basicpentagon} and \rf{q-exp} can be found in \cite{W,FKV,BT,Vo}.

\section{Operator-theoretical background}\label{opback}

\setcounter{equation}{0}

Let $\CH$ be a separable Hilbert space. The algebra of bounded operators
in $\CH$ will be denoted by $\BFB(\CH)$. We will only need some of the
most basic notions and results from functional analysis and 
the theory of operator
algebras as summarized e.g. in \cite[Chapter 14]{Wa}.
For the reader's convenience  and to fix some notations we shall
formulate the results that we need below.

\begin{thm} \label{specthm} 
$\;\frac{\quad}{}$ \cite[Theorem 14.8.14]{Wa} $\;\frac{\quad}{}$\\
Let $\BFC$ be a commutative von Neumann subalgebra of $\BFB(\CH)$,
where $\CH$ is a separable Hilbert space. Then there exists a compact, 
separable Hausdorff space $X$, a Radon measure $\mu$ on $X$, a measurable
family of Hilbert spaces $\{\CH_x\}_{x\in X}$, and a unitary bijection
\begin{equation}
\SU\,:\,\CH\,\to\,\CH_{\BFC}^{}\,\df\,\int_Xd\mu(s)\;\CH_s
\end{equation} 
such that 
\begin{equation}
\SU\cdot\BFC\cdot\SU^{\dagger}\;=\;\{\,\sm_f\,;\,f\in L^{\infty}(x,d\mu)\,\}\,,
\end{equation}
where $\sm_f$ is the multiplication operator defined by 
\begin{equation}
\sm_f\,:\,\CH_{\BFC}^{}\,\ni\,\{\Psi_x\}_{x\in X}\;\ra\;
\{f(x)\Psi_x\}_{x\in X}\,\in\,\CH_{\BFC}^{}.
\end{equation}
\end{thm}
If $\BFC$ is a commutative von Neumann subalgebra of $\BFB(\CH)$ we will
call
an operator $\SO$ on $\CH$ $\BFC$-decomposable if there
exists a family of bounded operators $\SO_x$,
defined on $\CH_x$ for $\mu$-almost every $x\in X$, such that
\begin{equation}
\SU\cdot\SO\cdot\SU^{\dagger}\;=\;\int_Xd\mu(s)\;\SO_x\;.
\end{equation} 
Let $\BFD_\BFC^{}$ be the algebra of all $\BFC$-decomposable operators.
\begin{propn} \label{decomppropn}
$\;\frac{\quad}{}$ \cite[Proposition 14.8.8]{Wa} $\;\frac{\quad}{}$
$\BFD_\BFC^{}$ is the commutant of $\BFC$ within $\BFB(\CH)$.
\end{propn}

The unbounded operators that we will have to deal with
will all be self-adjoint. We will freely use standard functional 
calculus for self-adjoint operators. 
When we say that two unbounded self-adjoint operators commute,
\[
\,[\,\SA\, ,\,\SB\,]\;=\;0,
\]
we will always mean commutativity of the spectral projections.
Let $\CF=\{\SA_\imath\,;\,\imath\in\CI\}$ be a family of commuting
self-adjoint operators defined on dense domains in a
separable Hilbert space
$\CH$. Standard functional 
calculus for commuting self-adjoint operators associates to 
$\CF$ a commutative von Neumann subalgebra $\BFC_\CF$ of $\BFB(\CH)$.
Theorem \ref{specthm} applied to $\BFC_\CF$ yields the existence
of a common spectral decomposition for the family $\CF$,
where $X$ represents the one-point compactification of the
spectrum of $\CF$.\\[1.5ex]
{\bf Tensor product notation}\\[1.5ex]
For a given finite set $I$ we will often use the notation
$\CH^{\otimes\RI}$ instead of
$\CH^{\otimes{\rm card}(\RI)}$, where ${\rm card}(\RI)$ is the number
of elements in $\RI$. In order to avoid fixing a numbering of the
elements of $\RI$ we 
will find it useful to employ the
following ``leg-numbering'' notation. To a given a subset $\RJ\subset \RI$
we may associate the canonical permutation operators 
\[
\SP_{\RJ}~:~\CH^{\otimes \RI}\ra \CH^{\ot \RJ}
\ot \CH^{\otimes \RI\setminus \RJ}\,.
\]
To an operator 
$\SO\in\BFB(\CH^{\ot \RJ})$ we may then associate an 
operator $\SO_{\RJ}\in\BFB(\CH^{\ot \RI})$ via
\begin{equation}\label{legnotation}
\SO_{\RJ}\df\SP_{\RJ}^{-1}\cdot(\SO\ot\id)\cdot\SP_{\RJ}^{\phantom{1}}.
\end{equation}
If $\RJ=\{i_1,i_2,\dots\}$ we will sometimes 
write $\SO_{i_1i_2\dots}$ instead of $\SO_{\RJ}$. 
We will also abbreviate $\SO_{\RJ_1\cup\RJ_2\cup\dots}$ to 
$\SO_{\RJ_1\RJ_2\dots}$.

\section{On the proof of Theorem \ref{Ptrelcompl}}\label{Ptrelproof}
\renewcommand{\hvf}{\tilde{\vf}}
\renewcommand{\tpi}{\tilde{\pi}}

A similar statement is known \cite{CF} for the closely related 
groupoid $\widetilde{\rm Pt}(\Sigma)$ whose elements are the moves
$[\hvf_\2,\hvf_\1]$ between 
fat graphs $\hvf_\2$, $\hvf_\1$ which do not have 
the decoration introduced in Subsection 5.1, but which  have
a numbering of the edges instead. 

\vspace{1ex}
\noindent
{\bf Theorem 2'.} {\it The groupoid $\widetilde{\rm Pt}(\Sigma)$ is the path
groupoid of the complex $\widetilde{\CP t}(\Sigma)$ which has vertices $\hvf$
and the following generators and relations:\\[.5ex]
{\it Generators:}

(i) Flips $ F_e$ along the edges $e\in\hvf_\1$ 
(see Figure \rf{flip}).

(ii) The exchanges $(ef)$ of the numbers assigned to edges $e$ and $f$.

\noindent{\it Relations:}

(a) There is no vertex that both edges $e$ and $f$ are incident to: 
$ F_f\circ F_e= F_e\circ F_f$.

(b) The edges $e$ and $f$ are incident to the same vertex: 
$ F_e\circ F_f\circ F_e\circ F_f\circ F_e=(ef).$

(c) $(ef)\circ F_e= F_f\circ(ef)$

(d) $ F_e\circ F_e={\rm id}$.}\\[1ex]
It was observed in \cite{P3,CF,CP} that
Theorem 2' follows quite easily from the fact \cite{P1,Ko} 
that the fat graphs
$\hvf$ can be used to label cells in certain cell decompositions
\cite{Ha,P1} of $\CT(\Sigma)\ti\BR_+^s$. This allows one to 
associate a path $\hat{\pi}$ in  $\CT(\Sigma)\ti\BR_+^s$ to each path
$\tilde{\pi}\in[\hvf_\2,\hvf_\1]$ in $\widetilde{\rm Pt}(\Sigma)$. Each 
closed path $\hat{\pi}$ in $\CT(\Sigma)\ti\BR_+^s$
can be deformed into small circles around the codimension two faces
of the cell decompositions from \cite{Ha,P1}. The latter are
easily identified with the relations listed in Theorem 2'. This 
implies that no relations other than those listed
in Theorem 2' are needed to contract a closed 
path $\tilde{\pi}\in[\hvf,\hvf]$ to the identity.

It remains to show that Theorem 2 follows from Theorem 2'.
To this aim it is important to observe that $\widetilde{\rm Pt}(\Sigma)$
imbeds into ${\rm Pt}(\Sigma)$ by means of the following construction:
Given the numbered fat graph $\hvf$ let us define
a decorated fat graph $\vf$ according to the following
rule: Pick any numbering of the vertices of $\hvf$.
For each vertex $v$ let us distinguish among the edges which
emanate from $v$ the one with the smallest number assigned to it.
For each move $ F_e$ let us choose a lift $\phi_e$
of the form $\phi_e=\de'_e\circ\om_{u_ev_e}\circ\de_e$, where
the vertices $u_e$ and $v_e$ represent are the ends of $e$ and
$\de'_e$, $\de_e^{}$ change the decoration only. 
It is then crucial to check that the relations of ${\rm Pt}(\Sigma)$ 
ensure that the image of $\widetilde{\rm Pt}(\Sigma)$ within 
${\rm Pt}(\Sigma)$ is simply connected. This follows if the 
images of the relations listed in Theorem 2' are contractible within 
${\rm Pt}(\Sigma)$. We have drawn a particular example for such an
image in Figure \ref{penta} below. It is easy to verify that the
image of the corresponding path in ${\rm Pt}(\Sigma)$ is closed. More 
generally one needs to consider all relations obtained by changes 
of the numbering of the edges. These relations are obtained 
by inserting $\eta^{-1}\circ\eta$ between the moves 
which occur in the image of the relation under consideration, 
where $\eta$ represents the change of decoration induced by a
change of numbering. In a similar way one may convince oneself
that all the relations listed in Theorem 2' are mapped to closed 
paths in ${\rm Pt}(\Sigma)$, which completes the proof that
$\widetilde{\rm Pt}(\Sigma)$ embeds into ${\rm Pt}(\Sigma)$.
\begin{figure}[htb]
\epsfxsize8cm
\centerline{\epsfbox{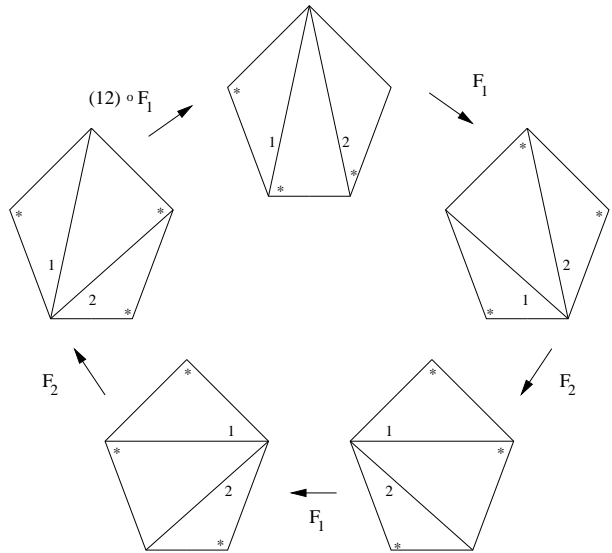}}
\caption{Image of the pentagon in ${\rm Pt}(\Sigma)$.}\label{penta}
\end{figure}

It remains to show that each closed path in ${\rm Pt}(\Sigma)$ is homotopic 
to a closed path in the image of $\widetilde{\rm Pt}(\Sigma)$.
Let us consider a closed path $\pi\in[\vf,\vf]\in{\rm Pt}(\Sigma)$.
Write $\pi=\om_n\circ\dots\circ\om_1$,
where each $\om_k$ is of the form $\om_k=\de'_k\circ\om_{u_kv_k}\circ 
\de_k^{}$, 
with $\de_k$, $\de'_k$ being composed out of moves $\rho_w$ and permutations
only. Choose any numbering for the edges of $\vf$ and denote 
the resulting numbered fat graph by $\hvf$. 
The path $\pi$ then defines a path $\tilde{\pi}$ in 
$\widetilde{\rm Pt}(\Sigma)$ by means of the following construction:
Substitute the moves
$\om_r$ for $k=1,\dots, n$  by the 
corresponding flips along the edges which connect $u_k$ and $v_k$,
and then multiply the result by the necessary permutations of the numbers
assigned to the edges of $\hvf$. $\tilde{\pi}$ has the
form $\tilde{\pi}=\Phi_n\circ\dots\circ\Phi_1$, where 
$\Phi_k\in[\hvf_{k+1},\hvf_k]$ factorize into a single flip $F_k$ times
a change of numbering $P_k$, $\Phi_k=P_k\circ F_k$.
 
The path $\tpi$ is mapped to a path $\pi'$ in ${\rm Pt}(\Sigma)$
by means of the construction above. 
$\pi'$ has the form $\pi'=\phi_n'\circ\dots\circ\phi_1'$
with $\phi_k'\in[\vf_{k+1},\vf_k]$. 
On the other hand let us note that the path
$\pi''=\phi_n''\circ\dots\circ\phi_1''$ with 
$\phi_k''\equiv\eta_{k+1}^{-1}\circ\om_k\circ\eta_k^{}$
is clearly homotopic to $\pi$ for all changes of decorations $\eta_k$
which satisfy $\eta_{n+1}=\eta_1$. 
For suitable choice of the $\eta_k$ one gets
$\phi_k''\in[\vf_{k+1},\vf_k]$. It is then  easy to see that 
$\phi_k'$ and $\phi_k''$ are homotopic.
\hfill$\square$

\section{Proof of Proposition \ref{locrelfact}}\label{locrelfactpf}

One may easily convince oneself that the moves $Z_p$
and $B_p$ can generically {\it not} be realized locally. 
There are simple sufficient criteria for a 
move $m$ to be realizable locally. 

\begin{lem} $\quad$ \label{locrelcrit}
\begin{itemize}
\item[(i)]
A move $m=[\tau_{m},\si_{m}]\in\CM_\1'(\Sigma)$ 
can always be realized locally if all but one of 
the boundary components of $\Sigma_m$
are contained in $\{c_e;e\in\CA_{\si}\}$.
\item[(ii)] The moves $B_p$, $B_p'$, $T_p$ can be realized locally if
the curve $c_{p,\1}$ which represents the boundary component
of $\ft_p$ assigned number $\1$ in Figure \ref{bmove} is contained
in $\CA_\si$.
\item[(iii)] The moves $F_{pq}$, $A_{pq}$, $S_p$, $S_{pq}$ and
$B_{pq}$ can always be realized locally.
\end{itemize}
\end{lem} 
\begin{proof}
Straightforward verifications.
\end{proof}


\begin{propn}\label{lrdeform} {\bf ( $\equiv$ Proposition 10)}
If $m\in\CM_\1'(\Sigma)$ is admissible, but can not be realized locally,
there always exists a path $\varpi_m$ which is (i) homotopic to
$m$ within $\CM^{\rm ad}_\1(\Sigma)$, and (ii) takes the
form 
\begin{equation}
\varpi_m\,=\,Y_m^{\phantom{1}}\circ m\circ Y_m^{-1},
\end{equation}
where $Y_m$ is a chain composed out of $Z_p$-moves and $F_{pq}$-moves
which can all be realized locally.
\end{propn}

\begin{proof} 
The claim will follow easily from another
auxiliary result that we will formulate as a separate lemma.
As a preparation, let us recall that cutting a surface $\Sigma$ 
along the curves $c(\si,e)$, $e\in\CA_{\si}$ produces  
a surface $\Sigma^{\dagger}$, the connected components of which
all have genus zero. Let us call a marking $\si$ {\it irreducible}
if the set of connected components of $\Sigma^{\dagger}$ has only
one element. It is easy to see that a marking $\si$ is 
irreducible iff it has precisely one outgoing external edge.

\begin{lem}\label{rotlem}
Let $\si$ be a marking on a surface $\Sigma$ 
that is irreducible. 
For any chosen boundary component $b$ of $\Sigma$
there exists a chain of 
$Z_p$-moves and $F_{pq}$-moves that (i) preserves admissibility
in each step, (ii) consists only of moves that can be realized
locally, and (iii)
transforms $\si$ to a irreducible marking $\si'$ whose
outgoing external edge ends in the chosen component $b$.
\end{lem}

\noindent{\it Sketch of proof}. 
In order to check the following arguments 
it may be useful to think of an irreducible marking 
as being obtained from 
the corresponding marking $\si^\dagger$ 
on the surface $\Sigma^{\dagger}$ of genus
zero by identifying the appropriate boundary circles.
To each incoming boundary component of $\Sigma^{\dagger}$ there
corresponds a unique element of $\CA_\si$.

The following claim is easy to verify.
By means of $F_{pq}$-moves one may transform $\si$ to a
marking $\tilde{\si}$ which has the following two properties.
First, for each vertex $p\in\tilde{\si}_\0$ at least one of the
edges that are incident at $p$ is contained in $\CA_\si$. 
Second, there is no edge $e\in\si_\1$ which connects a vertex $p$ to
itself. 
These two properties
insure that for each $p\in\tilde{\si}_\0$ either $Z_p$ or $Z_p^{-1}$
preserves the admissibility of the marking. 
It is then not very hard to construct a sequence of $Z_p$-moves
that (i) can all be realized locally and (ii) which 
transform $\tilde{\si}$ to a marking 
$\tilde{\si}'$ whose
outgoing external edge ends in the chosen boundary component $b$. 

By means of a chain of $F_{pq}$-moves one may finally 
transform $\tilde{\si}'$ back to a marking $\si'$ that differs from 
the original marking only by the desired change of decoration.
\hfill$\square$

{\it End of proof of Lemma \ref{lrdeform}}. 
Lemma \ref{locrelcrit} allows us to restrict
attention to the cases $m=Z_p,B_p,B_p',T_p$.
We need to transform the original marking 
to another one which has the property 
that both incoming boundary components
of the trinion $\ft_p$ are contained in $\CA_\si$. 
Given a curve $c\in\CC_\si\setminus\CA_\si$ there is a unique subsurface
$\Sigma_c\hookrightarrow\Sigma$ with marking $\si_c$ such 
that (i) $c$ is the unique outgoing curve in the boundary of $\Sigma_c$, and
(ii) the incoming boundary curves of $\Sigma_c$
are contained in $\CA_\si$. 
In order to infer the existence of the chain $Y_p$ it 
clearly suffices to apply Lemma \ref{rotlem} to the subsurfaces
$\Sigma_{c_\1}$ and  $\Sigma_{c_\2}$ respectively, where $c_{\1}$ 
and $c_{\2}$ are the incoming boundary components 
of the trinion with label $p$.
\end{proof}

Lemma \ref{lrdeform} implies that $\pi_m$ may be chosen as $\pi_{m}=
y_m\circ \hat{\pi}_m\circ y_m^{-1}$, where $y_m$ is
uniquely defined by the factorization of $Y_m$ into elementary moves that
can be realized locally, and $\hat{\pi}_m$ is the (fixed) path which
was chosen to represent $m$ in the case that $m$ can be realized locally.

\section{On the proof of Theorem \ref{lengthindep}} \label{lengthindeppf}

\subsection{A technical preliminary}

To begin with, let us present a technical result that
facilitates the explicit computations. 
The operator functions we will be interested in  are
of the form 
\begin{equation}\label{posexpops}
\SO\,\df\,\SO\big(\{\sz_\imath;\imath\in\vf_{{\1}}\}\big)\,=\,
\sum_{\tau\in\CF}\,C_\tau\,e^{\sx(\tau)}\, ,\qquad \sx(\tau)\,\df
\,\sum_{e\in\vf_\1}\,\tau(e)\,\sz_e,
\end{equation}
where the summation is taken over the space $\CF$ of all maps 
$\vf_\1\ni e\ra \tau(e)\in\frac{1}{2}\BZ$, and the coefficients
$C_\tau$ are assumed to be non-vanishing for a finite number of $\tau\in\CF$
only, in which case we assume $C_\tau\in\BR_+$. These operators
are densely defined and positive self-adjoint  due to the self-adjointness of 
$\sx(\tau)$.
The cone generated by operator functions of the form
 \rf{posexpops} will be denoted $\BFC_+(\vf)$.

\begin{figure}[t]\epsfxsize8cm
\centerline{\epsfbox{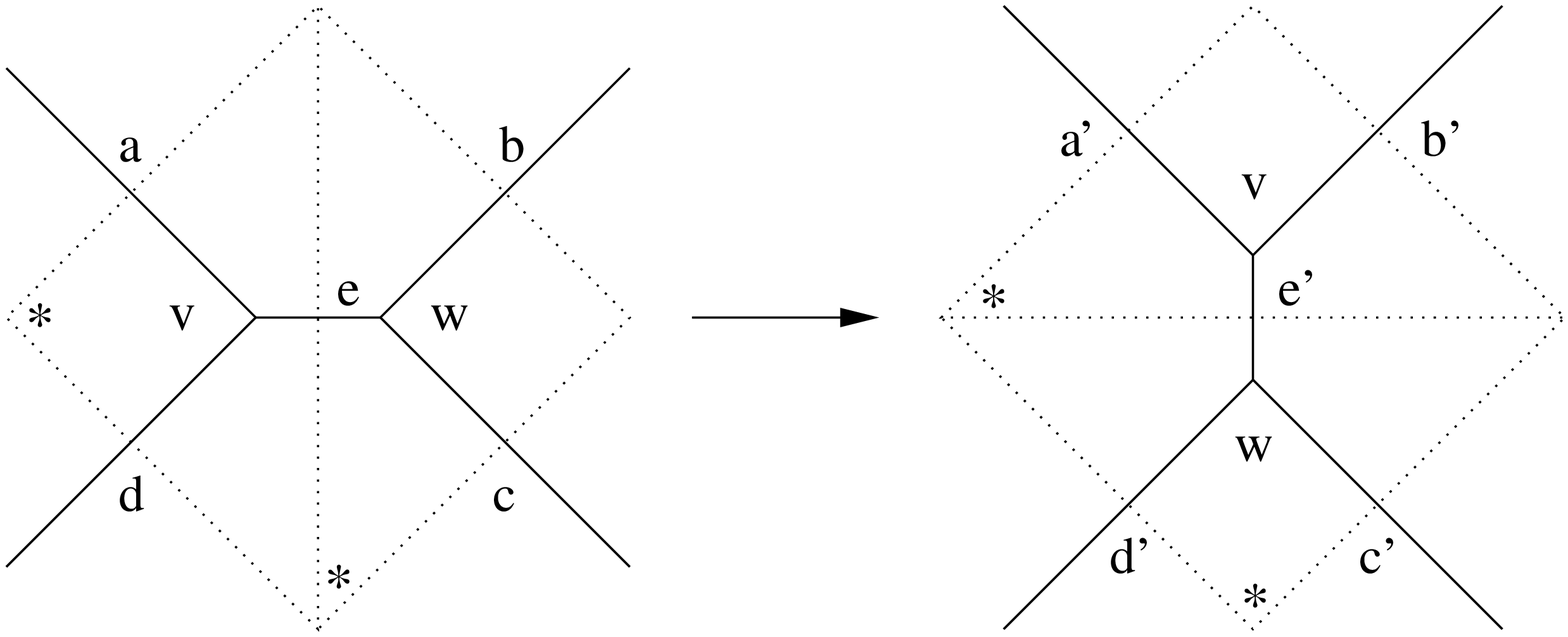}}
\caption{Labelling of the edges which are relevant for the 
description of a flip.}
\label{flip+}\end{figure}

\begin{propn}\label{Fockflip}
Let $\vf$, $\vf'$ be two fat graphs that are related by 
$\rho\df\omega_{vw}\in {\CP t}(\Sigma)$,
and let us adopt the labelling of the 
relevant edges given by Figure \ref{flip+}. 
For each $\SO\in\BFC_+(\vf)$ one has
\begin{equation}\label{SOfocktrf}
{} \sa_{\rho}(\SO)\,\equiv\,\su(\pi)\cdot\SO\cdot\su(\pi)^{-1}\,
=\,\SE_b^{}(\sz_{e'})\cdot\SO'\cdot
\left(\SE_b^{}(\sz_{e'})\right)^{-1},
\end{equation}
where $\SE_b(\sz)\df e_b\big(-\fr{\sz}{2\pi b}\big)$, and $\SO'$ is related 
to $\SO$ via
\begin{equation}\label{soprimedef}\begin{aligned}
\SO'\;\df\;&
\SO\big(\,\sz_{a'},\sz_{b'}+\sz_{e'},\sz_{c'},\sz_{d'}+\sz_{e'},-\sz_{e'},
\big\{\sz_\imath;\imath\in\vf_{{\1}}\setminus\{{a',b',c',d',e'\}}\big\}\,\big)
\;\;{\rm if}\\
\SO\;=\;
& \SO\big(\,\sz_a,\sz_b,\sz_c,\sz_d,\sz_e,
\big\{\sz_\imath;\imath\in\vf_{{\1}}\setminus\{{a,b,c,d,e\}}\big\}\,\big).
\end{aligned}
\end{equation}
\end{propn}

\begin{proof}
Given the explicit expression for $\ST_{vw}$ we only have to verify that
\[
e^{-2\pi i \spp_v\sq_w}\cdot\SO\cdot e^{-2\pi i \spp_v\sq_w}\;=\;\SO',
\]
with $\SO$ and $\SO'$ being related by \rf{soprimedef}. 
Keeping in mind \rf{kafock} we notice that the 
operators $\sz_f$ may be represented in the form $\sz_f=\sz_{f,v}+\sz_{f,w}$, 
where $v,w\in\vf_\0$ are the vertices connected by the edge $f\in\vf_\0$.
Let us label the vertices in $\vf_\0$ and $\vf'_\0$ such that the
edges $a$ and $d$ are incident to the 
vertices $v_a$ and $v_d$ besides to the vertex $v$ 
respectively, 
and similarly for the vertices $v_b$, $v_c$ and $w$. 
This leads to
the following representations for the relevant Fock variables:
\begin{equation}
\begin{aligned}
\sz_e = & 2\pi b(\sq_v-\spp_v+\spp_w)\, ,\\
\sz_a = & 2\pi b\spp_v+\sz_{a,v_a}\, ,\\
\sz_b = & 2\pi b(\sq_w-\spp_w)+\sz_{b,v_b}\, ,\\
\sz_c = & 2\pi b(-\sq_w)+\sz_{c,v_c}\, ,\\
\sz_d =  & 2\pi b(-\sq_v)+\sz_{d,v_d}\, ,
\end{aligned}\qquad
\begin{aligned}
\sz_{e'} = & 2\pi b(\sq_w'-\sq_v'-\spp_w')\, ,\\
\sz_{a'} = & 2\pi b\spp_v'+\sz_{a',v_{a}}\, ,\\
\sz_{b'} = & 2\pi b(\sq_v'-\spp_v')+\sz_{b',v_{b}}\, ,\\
\sz_{c'} = & 2\pi b(-\sq_w')+\sz_{c',v_{c}}\, ,\\
\sz_{d'} =  & 2\pi b\spp_w'+\sz_{d',v_{d}}\, ,
\end{aligned}\end{equation}
To complete the proof of Proposition 
\ref{Fockflip} is now the matter of a straightforward calculation.
\end{proof}

\begin{rem}
This yields in particular the formulae
\begin{equation*}\label{flipfovar}\begin{aligned}
e^{-\sz_{a'}}=& e^{-\frac{1}{2}\sz_a}(1+e^{-\sz_e})e^{-\frac{1}{2}\sz_a},\\
e^{+\sz_{d'}}=& e^{+\frac{1}{2}\sz_d}(1+e^{+\sz_e})e^{+\frac{1}{2}\sz_a},
\end{aligned} 
\qquad \sz_{e'}=-\sz_e,\qquad\begin{aligned}
e^{-\sz_{b'}}=& e^{+\frac{1}{2}\sz_b}(1+e^{+\sz_e})e^{+\frac{1}{2}\sz_b},\\
e^{+\sz_{c'}}=& e^{-\frac{1}{2}\sz_c}(1+e^{-\sz_e})e^{-\frac{1}{2}\sz_c}.
\end{aligned}
\end{equation*}
It is quite obvious that these transformations reduce to their
classical counterparts \rf{flipfovarclass} in the
limit $b\ra 0$. 
\end{rem}

As an example let us consider the monomials 
$
\SM_{n_an_c}^{n_bn_d}\;\df\; e^{\frac{1}{2}(n_a\sz_a+n_b\sz_b+n_c\sz_c
+n_d\sz_d)},
$
where $n_{\flat}$, $\flat\in\{a,b,c,d\}$ are restricted by
the requirement that $2N\df n_a+n_c-n_b-n_d$ must be even. 
One then finds certain simplifications on the right hand side of 
\rf{SOfocktrf}:
\begin{equation}\label{Monosimp}\begin{aligned}
\SE_b^{}(\sz_e')\cdot \SM_{n_an_c}^{n_bn_d} \cdot
\left(\SE_b^{}(\sz_e')\right)^{-1}\;=\;&
\sqrt{\SM_{n_an_c}^{n_bn_d}}\;
\frac{\SE_b^{}(\sz_e'+\pi ib^2 N)}{\SE_b^{}(\sz_e'-\pi ib^2 N)}\;
\sqrt{\SM_{n_an_c}^{n_bn_d}}\\
\;=\;&\sqrt{\SM_{n_an_c}^{n_bn_d}}\left[
\prod_{m=-\frac{N-1}{2}}^{\frac{N-1}{2}}(1+q^{2m} e^{-\sz_e'})
\right]
\sqrt{\SM_{n_an_c}^{n_bn_d}}\;.
\end{aligned}
\end{equation}

\subsection{Invariance of length operators}

In the main text we have explained how to reduce the proof
of Theorem \ref{lengthindep} to the following Lemma:

\begin{lem} $\quad$ \label{localprop2}
Assume that $m\in\CM_\1'(\Sigma)$ can be realized locally.
We then have the relation
\begin{equation}\label{inveqn}
\su(\pi_m)\cdot \SL_{\hsi_m,c} \;=\;
\SL_{\csi_m,c}\cdot\su(\pi_m)
\end{equation}
for all curves $c$ such that $c\in\CC_{\hsi_m}$ and $c\in\CC_{\csi_m}$.
\end{lem}

{\it On the proof of Lemma \ref{localprop2}.}
It is easy to see that our recursive Definition \ref{lengthdef} of the
length operators $\SL_{\si,c}$ reduces the task of proving 
Lemma \ref{localprop2} 
to the curves that represent a boundary component
of $\Sigma_m$. Let us write $\SL_{\si,e}\equiv\SL_{\si,c(\si,e)}$.

In the cases $m=F_{pq}, S_p,T_p,Z_p$ 
we will need to verify Lemma \ref{localprop2}
by direct calculations. 
Let us begin with the case $m=F_{pq}$. 
The claim is again trivial for the length operators
associated to the incoming boundary components
of $\Sigma_m$. Let us therefore focus on the length operators
$\SL_{\hsi_m,f}$, $\SL_{\csi_m,f}$ assigned to the outgoing boundary component
of $\Sigma_m$.
\begin{figure}[b]
\epsfxsize10cm
\centerline{\epsfbox{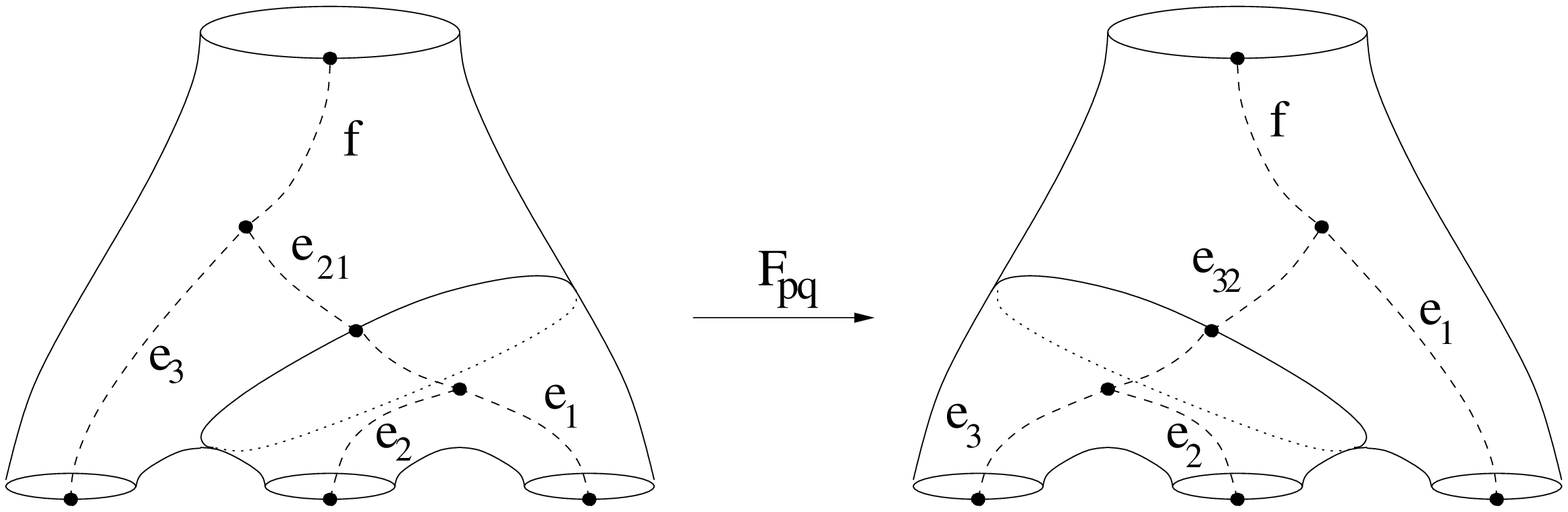}}
\caption{Labelling of the edges of $\si_m$ for $m=F_{pq}$.}\label{fusedges}
\end{figure}
The Definition \ref{lengthdef} yields  the following expressions
for the length operators $\SL_{\hsi_m,f}$ and $\SL_{\csi_m,f}$
respectively.
\begin{equation}
\begin{aligned}
\SL_{\hsi_m,f} \,=\,& e^{\sz_{\2\1}+\sy_\3+\sy_\2+\sy_\1}+
                   e^{\sz_{\2\1}-\sy_\3+\sy_\2+\sy_\1}+
e^{-\sz_{\2\1}-\sy_\3-\sy_\2-\sy_\1}\\
&+e^{\sy_\1-\sy_\2-\sy_\3}+e^{\sy_\1+\sy_\2-\sy_\3}+
  e^{-\sy_\1-\sy_\2-\sy_\3}\\ 
&+e^{\sz_{\2\1}+\sy_\2+\sy_\1}\SL_{\3}+
e^{\sy_\1-\sy_\3}\SL_\2+
e^{-\sy_\2-\sy_\3}\SL_\1\,,\\
\SL_{\csi_m,f} \,=\,& e^{+\sz_{\3\2}'+\sy_\3'+\sy_\2'+\sy_\1'}+
                   e^{-\sz_{\3\2}'-\sy_\3'-\sy_\2'-\sy_\1'}+
e^{-\sz_{\3\2}'-\sy_\3'-\sy_\2'+\sy_\1'}\\
&+e^{\sy_\1'+\sy_\2'+\sy_\3'}+e^{\sy_\1'-\sy_\2'-\sy_\3'}+
  e^{\sy_\1'+\sy_\2'-\sy_\3'}\\ 
&+e^{\sy_\2'+\sy_\1'}\SL_{\3}'+
e^{\sy_\1'-\sy_\3'}\SL_\2'+
e^{-\sz_{\3\2}'-\sy_\2'-\sy_\3'}\SL_\1'\,.
\end{aligned}
\end{equation}
We have used the abbreviations  $\sy_\iota\df\sy_{\vf_\hsi,e_\iota}$, 
$\sy_\iota'\df\sy_{\vf_\csi,e_\iota}$,
$\SL_{\iota}\df\SL_{\hsi,e_\iota}$ and $\SL_{\iota}'\df\SL_{\csi,e_\iota}$ 
for $\iota\in\{\1,\2,\3\}$, as well as 
$\sz_\iota\df\sz_{\vf_\hsi,e_\iota}$, 
$\sz_\iota'\df\sz_{\vf_\csi,e_\iota}$
for $\iota\in\{\3\2,\2\1\}$.
The labelling of the edges is the one introduced in Figure
\ref{fusedges}. According to Proposition \ref{Fockflip} we need
to calculate
\begin{equation}\label{calc}
\begin{aligned}
E_b(\sz_{\3\2}')\cdot \Big(& e^{\sy_\3'+\sy_\2'+\sy_\1'+\sz_{\3\2}^{}}+
e^{\sy_\1'-\sy_\2'-\sy_\3'}+\\
 +&e^{\sy_\2'+\sy_\1'}\SL_{\3}'+
e^{\sy_\1'-\sy_\3'}\SL_\2'+
e^{-\sz_{\3\2}'-\sy_\2'-\sy_\3'}\SL_\3'+\\
+& e^{\frac{1}{2}(\sy_\1'+\sy_\2'-\sy_\3')}(1+e^{-\sz_{\3\2}'})
   e^{\frac{1}{2}(\sy_\1'+\sy_\2'-\sy_\3')}+\\
+& e^{-\frac{1}{2}(\sy_\1'+\sy_\2'+\sy_\3'+\sz_{\3\2}')}(1+e^{-\sz_{\3\2}'})
   e^{-\frac{1}{2}(\sy_\1'+\sy_\2'+\sy_\3'+\sz_{\3\2}')}\Big) \cdot E_b(\sz_{\3\2}')^{\dagger}\,.
\end{aligned}
\end{equation}
We finally need to apply equation \rf{Monosimp}. The terms in the first 
line of \rf{calc} have $N=1$, those in the second line $N=0$, and all
other terms have $N=-1$. Straightforward application
of equation \rf{Monosimp} shows that the expression given 
in \rf{calc} equals $\SL_{\csi_m,f}$, as claimed.

\begin{figure}[t]
\epsfxsize10cm
\centerline{\epsfbox{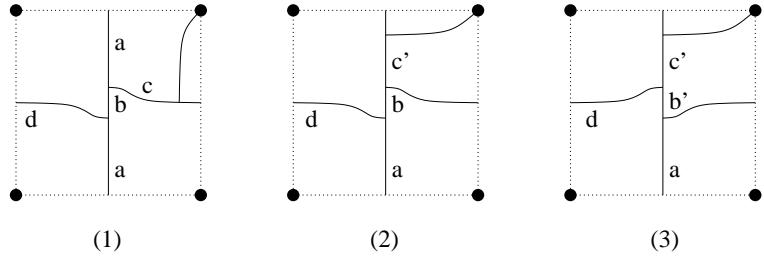}}
\caption{Fat graphs $\vf_i$ on the one-holed torus.}\label{smovegraphs}
\end{figure}

The next case we will consider is $m=S_p$. It clearly
suffices to restrict attention to the case that
$\Sigma_m=\Sigma_{\1,\1}$, the one-holed torus. On $\Sigma_m$ let us consider 
the fat graphs $\vf_i$, $i=\1,\2,\3$ depicted in 
Figure \ref{smovegraphs}. 

We have
$\vf_\1=\hvf_m'$ and $\vf_\3=\cvf_m'$ 
respectively. The sequence of fat graph $(\vf_\3,\vf_\2,\vf_\1)$ 
defines the path $\pi_m\in [\cvf_m',\hvf_m']$ which we will use.
The only relevant curve is the curve $\be$ which represents the
boundary of $\Sigma_m$.
Definition \ref{lengthdef} 
yields the following expressions for the length operators
$\SL_{\vf_\1,\be}$ and $\SL_{\vf_\3,\be}$ respectively:
\begin{equation}
\begin{aligned}
\SL_{\vf_\1,\be}\;=\;& e^{\sz_d-\sz_c}+2\cosh(\sz_c+\sz_d+\sz_a+\sz_b)\\
&+e^{\sz_d}(1+e^{\sz_b}+e^{\sz_b+\sz_a})
+e^{-\sz_c}(1+e^{-\sz_a}+e^{-\sz_a-\sz_b}),\\
\SL_{\vf_\3,\be}\;=\;& e^{\sz_{c'}-\sz_a}+
2\cosh(\sz_a+\sz_d+\sz_{c'}+\sz_{b'})\\
&+e^{\sz_{c'}}(1+e^{\sz_{b'}}+e^{\sz_{b'}+\sz_d})
+e^{-\sz_a}(1+e^{-\sz_d}+e^{-\sz_d-\sz_{b'}})\,.
\end{aligned}
\end{equation}
The labelling of the relevant edges is the one introduced in
Figure \ref{smovegraphs}.  
With the help of Proposition \ref{Fockflip} one may calculate
\begin{equation}\label{interstep}
\begin{aligned}
\SL_{\vf_\2,\be}\,\df\,
\sa_{[\vf_\2,\vf_\1]}\big(\SL_{\vf_\1;\be}\big)\;=\;&
e^{\sz_{c'}-\sz_a}+
2\cosh(\sz_a+\sz_d+\sz_{c'}+\sz_{b})\\
+&e^{\sz_{c'}}(1+e^{\sz_{d}}+e^{\sz_{d}+\sz_b})
+e^{-\sz_a}(1+e^{-\sz_b}+e^{-\sz_b-\sz_{b}})\,.
\end{aligned}
\end{equation}
One may then apply Proposition \ref{Fockflip} once more to calculate
$\sa_{[\vf_\3,\vf_\2]}\big(\SL_{\vf_\2,\be}\big)$.
The result is $\SL_{\vf_\3,\be}$ as claimed.

We have furthermore verified by direct calculations 
that Lemma \ref{localprop2} 
holds in the cases $m=T_p,Z_p,B_p,B_p'$ respectively. These 
calculations 
proceed along very similar lines as in the previous two cases,
which motivates us to omit the details.

In the remaining cases $m=A_{pq}$, $m=B_{pq}$ and $m=S_{pq}$
one may use the factorization of the chains $\pi_m$ which follows 
from the definitions \rf{Adef}, \rf{bcompdef} and \rf{scompdef}
respectively in order to reduce the proof to the cases
where the result is already established.
This completes our discussion of the 
proof of Lemma \ref{localprop2}.\hfill $\square$



\end{document}